\documentclass{article}

\pdfoutput=1

\usepackage{amsmath}
\usepackage{paralist}
\usepackage{graphicx}
\usepackage{epstopdf}
\usepackage[colorlinks=true]{hyperref}
  \hypersetup{urlcolor=blue, citecolor=red}

\newtheorem{proposition}{Proposition}
\title{A Generalized Bidiagonal-Tikhonov Method Applied To Differential Phase Contrast Tomography}

\author{}

\usepackage{amsfonts}
\usepackage{amssymb}
\usepackage{bm}
\usepackage{mathrsfs}
\usepackage{algorithm}
\usepackage{algorithmic}
\usepackage{tikz}
\usepackage{enumerate}
\usepackage{subfig}

\graphicspath{{Images/}}
\usetikzlibrary{arrows}
\newcommand{\mbbR}{\mathbb{R}}
\newcommand{\mbbC}{\mathbb{C}}
\newcommand{\bfx}{\bm{x}}
\newcommand{\bft}{\bm{\theta}}
\DeclareMathOperator{\radon}{\mathscr{R}}
\DeclareMathOperator{\diff}{\ensuremath{\mathrm{d}}}
\DeclareMathOperator*{\argmin}{arg\, min}
\DeclareMathOperator{\spn}{span}
\DeclareMathOperator{\Null}{Null}

\newcommand\numberthis{\addtocounter{equation}{1}\tag{\theequation}}% Number 1 equation in align*.
\newcommand{\hypref}[2]{\hyperref[#2]{#1 \ref*{#2}}}

\begin{document}

\maketitle

% Enter the first author's name and address:
\centerline{\scshape Nick Schenkels}
\medskip
{\footnotesize
% please put the address of the first author
 \centerline{Department of Mathematics and Computer Science,}
   \centerline{University of Antwerp, Belgium}
	 \centerline{Middelheimlaan 1, M.G.328, 2020, Antwerpen, Belgi\"e}
} % Do not forget to end the {\footnotesize by the sign }

\medskip

\centerline{\scshape Wim {van Aarle}}
\medskip
{\footnotesize
 % please put the address of the second  and third author
 \centerline{Department of Physics, Vision Lab}
   \centerline{University of Antwerp, Belgium}
   \centerline{Universiteitsplein 1, D.N.114, 2610, Wilrijk, Belgi\"e}
}

\medskip

\centerline{\scshape Jan Sijbers}
\medskip
{\footnotesize
 % please put the address of the second  and third author
 \centerline{Department of Physics, Vision Lab}
   \centerline{University of Antwerp, Belgium}
   \centerline{Universiteitsplein 1, D.N.113, 2610, Wilrijk, Belgi\"e}
}

\medskip

\centerline{\scshape Wim Vanroose}
\medskip
{\footnotesize
 % please put the address of the second  and third author
 \centerline{Department of Mathematics and Computer Science,}
   \centerline{University of Antwerp, Belgium}
	 \centerline{Middelheimlaan 1, M.G.210, 2020, Antwerpen, Belgi\"e}
}

\bigskip

\begin{abstract}
Phase contrast tomography is an alternative to classic absorption
contrast tomography that leads to higher contrast reconstructions in
many applications. We review how phase contrast data can be acquired
by using a combination of phase and absorption gratings. 
Using algebraic reconstruction techniques the object can be reconstructed
from the measured data. In order to solve the resulting linear system we
propose the Generalized Bidiagonal Tikhonov (GBiT) method, an
adaptation of the generalized Arnoldi-Tikhonov method that uses the
bidiagonal decomposition of the matrix instead of the Arnoldi
decomposition. We also study the effect of the finite difference
operator in the model by examining the reconstructions with either a
forward difference or a central difference approximation. We validate
our conclusions with simulated and experimental data.
\end{abstract}

%%%%%%%%%%%%%%%%%%%%%%%% SECTION: Introduction %%%%%%%%%%%%%%%%%%%%%%%%
\section{Introduction}
Since their discovery by Wilhelm R\"ontgen in 1895, x-rays have become an invaluable tool in medical imaging, material testing, security screenings, etc. In classic x-ray computed tomography (CT), the absorption of the x-rays by the material is imaged. However, very different materials can have similar absorption effects, making them hard to differentiate in the reconstructed image. An alternative tomographic technique is phase contrast computed tomography (PC-CT). This method is based on the phase shift that occurs in the wave front when the x-rays pass through the object of interest. This shift is more sensitive to small changes in the material -- especially for soft tissues or samples with low atomic numbers -- this technique can be used to achieve a higher contrast in the reconstructions.

The reconstruction methods for absorption contrast CT (AB-CT) are well established in literature and can be divided in two groups. First there are analytical techniques, such as filtered back projection (FBP), which are based on the Radon transform and its relation to the Fourier transform as expressed by the Fourier Slice theorem \cite{Natterer2001}. The second type  are algebraic reconstruction methods, which are based on the discretization of the Radon transform $R\in\mbbR^{m\times n}$ and solve a linear system
\begin{equation}\label{eq:linsys}
Rx = b,
\end{equation}
using some direct or iterative method.  Here, $x\in\mbbR^n$ represents a rasterized version of the attenuation coefficients of the object and $b\in\mbbR^m$ the measured x-ray intensity data.  In practice, $b$ is subject to noise and measurement errors, which can result in various types of reconstruction artefacts.  Moreover, constraints on radiation dose, scanning time and scanning geometry can make \eqref{eq:linsys} an underdetermined system, leading to non-unique solutions.  Various regularization methods have been developed in order to find a good approximation of the image $x$.  These methods typically include the exploitation of prior knowledge about the scanned object.  One widely used approach is Tikhonov regularization, which solves:
\[
x_\lambda=\argmin_{x\in\mathcal{K}}\left\{\left\|Rx-b\right\|^2+\lambda\left\|L(x-x_0)\right\|^2\right\},
\]
where $\left\|\cdot\right\|$ denotes the 2-norm and $\mathcal{K}$ is a subspace of the data-space $\mbbR^n$. Prior knowledge of the solution is incorporated via the initial guess $x_0$ and the operator $L$, which can be used to place constraints on the solution $x_\lambda$. The choice of the regularization parameter $\lambda$ plays a key role here, since it determines the quality of the reconstruction. We come back to this in \hypref{section}{sec:artik}.

\subsection*{Differential phase contrast tomography}\label{sec:dpcct}
In traditional absorption contrast CT (AC-CT), the intensity of the x-rays is measured immediately after they pass through the object of interest. However, this doesn't provide any information about the occurred phase shift, information which can lead to a significantly enhanced contrast between two materials with similar absorption coefficients. Various different approaches have been developed in order to acquire data that can be linked to the phase shift: set-ups with a crystal interferometer \cite{Bonse1965, Momose1995}, propagation based imaging \cite{Bronnikov2002, Kostenko2013, Snigirev1995}, analyser based imaging \cite{Ingal1995, Davis1995} or differential phase contrast imaging \cite{Momose2003, Pfeiffer2006}.

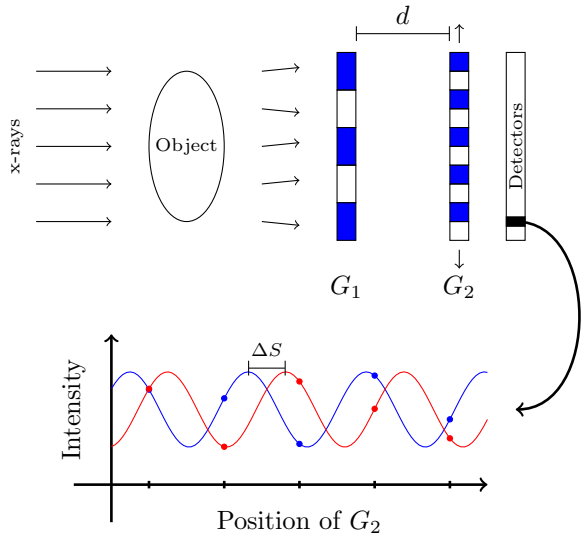
\begin{figure}
\centering
\begin{tikzpicture}[scale = 0.5]

% x-rays in:
\draw[->] (0,4) -- (2,4);
\draw[->] (0,3) -- (2,3);
\draw[->] (0,2) -- (2,2);
\draw[->] (0,1) -- (2,1);
\draw[->] (0,0) -- (2,0);
\node[rotate=90] at (-0.5,2) {\scriptsize x-rays};

% Ellips:
\draw  (4,2) ellipse (1 and 2);
\node at (4,2) {\scriptsize Object};

% xrays uit:
\draw[->] (6,4) -- (7,4.1);
\draw[->] (6,3) -- (7,2.9);
\draw[->] (6,2) -- (7,2.1);
\draw[->] (6,1) -- (7,1.1);
\draw[->] (6,0) -- (7,-0.1);

% Grating 1:
\draw [fill = blue]  (8,4.5) rectangle (8.5,3.5);
\draw  (8,3.5) rectangle (8.5,2.5);
\draw [fill = blue] (8,2.5) rectangle (8.5,1.5);
\draw  (8,1.5) rectangle (8.5,0.5);
\draw [fill = blue] (8,0.5) rectangle (8.5,-0.5);
\node at (8.25,-1.75) {$G_1$};

% Grating 2:
\draw [fill = blue] (11,4.5) rectangle (11.5,4);
\draw  (11,4) rectangle (11.5,3.5);
\draw [fill = blue] (11,3.5) rectangle (11.5,3);
\draw  (11,3) rectangle (11.5,2.5);
\draw [fill = blue] (11,2.5) rectangle (11.5,2);
\draw  (11,2) rectangle (11.5,1.5);
\draw [fill = blue] (11,1.5) rectangle (11.5,1);
\draw  (11,1) rectangle (11.5,0.5);
\draw [fill = blue] (11,0.5) rectangle (11.5,0);
\draw  (11,0) rectangle (11.5,-0.5);
\draw [->] (11.25, 4.75) -- (11.25,5.25);
\draw [->] (11.25,-0.75) -- (11.25,-1.25);
\node at (11.25,-1.75) {$G_2$};

% Detectors:
\draw  (12.5,4.5) rectangle (13,-0.5);
\node[rotate=90] at (12.75,2) {\scriptsize Detectors};
\draw[fill = black] (12.5, -0.125) rectangle (13, 0.125);
\draw[->, line width = 1, color = black] (12.75,0) .. controls (15, 0) and (15, -5) .. (12.75, -5);

% Distance d:
\draw[|-|] (8.5, 5) -- (11, 5);
\node at (9.75, 5.5) {$d$};

% Begin shift ----------------------------------------------------------------------------------
\begin{scope}[shift = {(-14, -10)}]

\scalebox{1}{
% Axis:
\draw[->, line width = 1] (16,2) -- (16,7);
\draw[->, line width = 1] (15,3) -- (26,3);
\draw[line width = 1] (17, 3.1) -- (17, 2.9);
\draw[line width = 1] (19, 3.1) -- (19, 2.9);
\draw[line width = 1] (21, 3.1) -- (21, 2.9);
\draw[line width = 1] (23, 3.1) -- (23, 2.9);
\draw[line width = 1] (25, 3.1) -- (25, 2.9);
\node at (21,2) {Position of $G_2$};
\node[rotate = 90] at (15,5) {Intensity};

% Sine wave 1:
\draw[domain=16:26,samples=100, blue] plot(\x,{sin(2*\x r) + 5});
\draw[blue, fill] (17, 5.5291) circle (2pt);
\draw[blue, fill] (19, 5.2964) circle (2pt);
\draw[blue, fill] (21, 4.0835) circle (2pt);
\draw[blue, fill] (23, 5.9018) circle (2pt);
\draw[blue, fill] (25, 4.7376) circle (2pt);

% Sine wave 2:
\draw[domain=16:26,samples=100, red] plot(\x,{sin(2*\x r- 2 r) + 5});
\draw[red, fill] (17, 5.5514) circle (2pt);
\draw[red, fill] (19, 4.0082) circle (2pt);
\draw[red, fill] (21, 5.7451) circle (2pt);
\draw[red, fill] (23, 5.0177) circle (2pt);
\draw[red, fill] (25, 4.2317) circle (2pt);

% Intensity Shift:
\draw[|-|] (19.6344, 6.1) -- (20.6344, 6.1);
\node at (20.1344, 6.6) {\scriptsize $\Delta S$};
}

\end{scope}
% End shift ------------------------------------------------------------------------------------

\end{tikzpicture}
\caption{Overview of the DPC set-up. The scanned object causes a phase shift in the wave front. The angular deviations this causes in the x-rays are proportional to the differential phase shift. By using a phase grating ($G_1$) and an absorption grating ($G_2$) with a phase-stepping procedure, it is possible to measure these small angles.}
\label{fig:dpcsetup}
\end{figure}

Differential phase contrast CT (DPC-CT) is based on the small angular deviations $\alpha$ that appear in the wave front after going through the scanned object, which in turn can be linked to the derivative of the phase shift $\Phi$. Recovery of these refraction angles based on the intensity measurements can be done using a reverse projection method \cite{Zhu2010}, a Fourier based algorithm \cite{Takeda1982} or a phase stepping approach. In the latter, two gratings are placed between the object and the detector array, see \hypref{Fig.}{fig:dpcsetup}. The first grating $G_1$ is a phase grating that, due to the Talbot effect and the Fresnel diffraction principle, causes an interference pattern at the Talbot distance $d$ \cite{Born1999}. In order to accurately measure this interference pattern, a second absorption grating $G_2$ is used. For different positions of $G_2$, measurements of the intensity of the outgoing x-rays are taken with and without an object present in the detector. By comparing the two intensity curves (position of $G_2$ versus the intensity measurement), the shift $\Delta S$ caused by the object can be measured. Since $\Delta S$ is related to $\alpha$, the refraction angles can be recovered. For a more detailed overview of this set-up and the phase stepping approach we refer to \cite{Momose2006, Pfeiffer2006}.

The reconstruction methods for DPC-CT are closely related to those for AB-CT. Analytically, the presence of a derivative of the Radon transformation in the FBP theorem results in a special complex valued filter. Algebraically on the other hand, the model can, for example, be represented using a discrete version for the derivative of the Radon transform \cite{nilchian2013} or by adding a finite difference operator $D$ as proposed in \cite{Fu2015, Fu2013, Fu2014}:
\[
DRx = b.
\]
Here, the resulting linear system is solved using a Kaczmarz method or algebraic reconstruction technique (ART).

In this paper we follow a similar algebraic approach, but we propose to use an adaptation of the generalized Arnoldi-Tikhonov method to solve the linear system. This is a Krylov subspace method that uses Tikhonov regularization and an secant update step based on the discrepancy principle in order to find a suitable value for the regularization parameter $\lambda$. This method was first introduced in \cite{Gazzola2014_2, Gazzola2014} where it was derived for square matrices. As we will illustrate further on, square matrices rarely appear in tomographic reconstructions -- absorption or phase contrast -- but the method can be adjusted using a bidiagonal decomposition of the matrix instead of the Arnoldi algorithm. Note that this generalization is independent of the tomographic setting we use as a test case. We will also take a closer look at how the choice of finite difference operator $D$ effects the reconstruction. We do this by comparing the forward difference scheme used in \cite{Fu2015, Fu2013, Fu2014} with a central difference approximation and study the effects they both have on the reconstruction.

\subsection*{Outline}
The outline of this paper is as follows. In \hypref{section}{sec:dpc} we give a brief review of absorption and phase contrast tomography and describe how they can be solved analytically and algebraically. In \hypref{section}{sec:artik} we discuss the Arnoldi-Tikhonov method and explain how it can be adjusted for non-square matrices. In \hypref{section}{sec:findiff} we examine how the finite difference operator  affects the solution by looking at a simple test problem. Finally, in \hypref{section}{sec:numapp}, we look at some numerical experiments with simulated and experimental data.

%%%%%%%%%%%%%%%%%%%%%%%% SECTION: Phase Contrast %%%%%%%%%%%%%%%%%%%%%%%%
\section{Review of the forward and inverse problem}\label{sec:dpc}
The propagation of a wave $u(\bfx)$ through a medium with no inside source is described by the homogeneous Helmholtz equation
\[
\left(\Delta + k(\bfx)\right)u(\bfx) = 0
\]
with $\bfx\in\mbbR^2$ and $k(\bfx)\in\mbbC$ the position dependent wave number. If $\theta$ is the direction of the incoming plane wave $u_{in}(\bfx) = E_0e^{ik_0\bfx\cdot\bft}$, with $\bft = (\cos(\theta), \sin(\theta))^T$ and amplitude $E_0$, then $k(\bfx) = k_0n(\bfx)$. Here,
\[
n(\bfx) = 1 - \delta(\bfx) + i\beta(\bfx)\in\mbbC
\]
is the refractive index of the material and $\beta, \delta:\mbbR^2\longrightarrow\mbbR$. Assuming that the wave number is constant, i.e. $k(\bfx) = k = k_0(1 - \delta + i\beta)$, then we can write the total wave as
\begin{equation}\label{eq:totalwave}
\begin{aligned}
u(\bfx) &= E_0e^{ik\bfx\cdot\bft}\\
				&= \Big(E_0e^{ik_0\bfx\cdot\bft}\Big)\Big(e^{-ik_0\delta\bfx\cdot\bft}\Big)\Big(e^{-k_0\beta\bfx\cdot\bft}\Big)\\
				&= \Big(\text{incoming wave}\Big)\Big(\text{phase shift}\Big)\Big(\text{amplitude change}\Big).
\end{aligned}
\end{equation}
This equation illustrates the effect of the medium on the incoming wave: the imaginary part $\beta$ of the refractive index causes a change in the amplitude of the wave, whereas the real part $\delta$ leads to a phase shift in the wave front. From \eqref{eq:totalwave} it follows that the intensity of the beam at $\bfx$ is given by
\begin{align*}
I(\bfx) &= \Big(\text{amplitude of the wave}\Big)^2\\
				&= \Big(E_0e^{-k_0\beta\bfx\cdot\bft}\Big)^2\\
				&= I_0e^{-2k_0\beta\bfx\cdot\bft},
\end{align*}
with $I_0 = E_0^2$ the intensity of the incoming wave. If the beam travels from $\bfx$ to $\bfx + \Delta\bfx$ the change in intensity can be approximated by
\begin{align*}
I(\bfx + \Delta\bfx) &\approx I(\bfx) + I'(\bfx)\Delta\bfx\cdot\bft\\
										 &= I(\bfx) - 2k_0\beta I(\bfx)\Delta\bfx\cdot\bft\\
\Rightarrow\ dI(\bfx) &= -2k_0\beta I(\bfx)d\bfx\cdot\bft.\numberthis\label{eq:diffint}
\end{align*}
Similarly, if $\Phi(\bfx)$ is the phase shift of the wave, then the phase shift in the wave propagating from $\bfx$ to $\bfx + \Delta\bfx$ is given by
\begin{align*}
&\Phi(\bfx + \Delta\bfx) - \Phi(\bfx) = k_0\delta\Delta\bfx\\
\Rightarrow\ &d\Phi(\bfx) = k_0\delta d\bfx\cdot\bft.\numberthis\label{eq:diffphase}
\end{align*}
Note that this is the phase shift relative to the incoming wave. To find the intensity and total phase shift at the detector for a non-constant wave number $k(\bfx)$, we have to integrate \eqref{eq:diffint} and \eqref{eq:diffphase} along the entire path of the beam $L(\theta, t) = \{\bfx\in\mbbR^2\mid\bfx\cdot\bft = t\}$, which yields:
\begin{equation}\label{eq:integrals}
\left\{\begin{aligned}
I(\theta, t) &= I_0e^{-2k_0\int_{L(\theta, t)}\beta(\bfx)d\bfx}\\
\Phi(\theta, t) &= \int_{L(\theta, t)}k_0\delta(\bfx)d\bfx
\end{aligned}\right..
\end{equation}
It should be noted that this is a very crude, linear approximation of the Helmholtz equation we started from. However, \eqref{eq:integrals} is much easier to solve and can still be used to reconstruct high quality x-ray images.

% ------------------------- SUBSECTION: Continuous level -------------------------
\subsection{Analytical Reconstruction}\label{sec:anrec}
Analytically, AB-CT is based on the filtered back projection theorem, which states that \cite{Natterer2001}:
\begin{equation}\label{eq:FBP}
\beta(\bfx) = \mathscr{B}\left(\mathscr{F}^{-1}|r|\mathscr{F}\mathscr{R}\beta\right)(\bfx).
\end{equation}
Here, $\mathscr{F}$ is the Fourier transform, $\radon$ and $\mathscr{B}$ the Radon transformation and its inverse and $r$ the radius of $\bfx$ in polar coordinates. Since it follows from \eqref{eq:integrals} that 
\[
\radon\beta(\theta, t) = \int_{L(\theta, t)}{\beta(\bfx)d\bfx} = -\frac{1}{2k_0}\ln\frac{I(\theta, t)}{I_0}
\]
it is possible to reconstruct $\beta$ via \eqref{eq:FBP} using the measurements of $I(\theta, t)$.

Differential phase contrast tomography, however, is based on the small angular deviations $\alpha(\theta, t)$ that occur in the wave front.  By applying a phase stepping technique as described in \hypref{section}{sec:dpcct}, the shift in the intensity measurements caused by the object can be determined. Denoting this shift as $\Delta S(\theta, t)$, they can be linked to the deviation angles in the x-rays by
\[
\Delta S(\theta, t) = d\alpha(\theta, t).
\]
Since these angles are very small, they can be approximated by
\[
%\alpha(\theta, t)=\arctan\left(\frac{1}{k_0}\frac{\partial\Phi(\theta, t)}{\partial t}\right)\approx\frac{1}{k_0}\frac{\partial\Phi(\theta, t)}{\partial t}.
\alpha(\theta, t)=\frac{1}{k_0}\frac{\partial\Phi(\theta, t)}{\partial t}.
\]
Substituting \eqref{eq:integrals} for $\Phi(\theta, t)$, we find that:
\begin{equation}\label{eq:alpha}
\alpha(\theta, t)=\frac{\partial\int_{L(\theta, t)}{\delta(x,y)ds}}{\partial t}=\frac{\partial}{\partial t}\radon\delta(\theta, t).
\end{equation}
This differs from AB-CT because here we have the derivative of the Radon transform. However, the Fourier transform property for derivations states that
\[
\mathscr{F}\left(\frac{df}{dx}\right)(r)=2\pi ri\mathscr{F}f(r),
\]
so it follows that
\[
\mathscr{F}\alpha(\theta, r)=\mathscr{F}\left(\frac{\partial}{\partial t}\radon\delta\right)(\theta, r)=2\pi ri\mathscr{F}\mathscr{R}\delta(\theta, r).
\]
This means that $\delta$ can be reconstructed similarly to \eqref{eq:FBP}, but with a different complex filter and measurements of $\alpha(\theta, t)$:
\begin{equation}\label{eq:FBPphase}
\delta(\bfx) = \mathscr{B}\left(\mathscr{F}^{-1}\frac{|r|}{2\pi ri}\mathscr{F}\alpha\right)(\bfx).
\end{equation}

% ------------------------- SUBSECTION: Discrete level  -------------------------
\subsection{Algebraic Reconstruction}\label{sec:disrec}
In both absorption and differential phase contrast tomography the measurements are related to the Radon transformation, a set of line integrals over the domain. Using a parallel beam projection model, see \hypref{Fig.}{fig:pbp} these can be discretized in the following way:
\begin{itemize}
\item Assume that there are $k$ detector elements and  that projections are made under $l$ different angles $\theta_1, \ldots, \theta_l\in[0, \pi[$. This means there are $m = kl$ different line integrals or measurements, which can be written as a vector $b=(b_i)_{i=1}^m\in\mathbb{R}^m$. Each component $b_i$ will be the measurement corresponding to a specific detector element and projection angle.
\item The domain is considered to be a grid of $n_x$ by $n_y$ pixels in the $x$ and $y$ direction, which are written column wise as a vector $x=(x_j)_{j=1}^n\in\mathbb{R}^n$ with $n=n_xn_y$.
\item For a measurement $b_i$, every pixel $x_j$ is given a weight $R=(r_{ij})_{ij}\in\mathbb{R}^{m\times n}$ such that
\[
b_i=\sum_{j=1}^n{r_{ij}x_j}\approx\text{ the corresponding line itegral}.
\]
The weights $r_{ij}$ can be chosen in different ways, see \cite{Joseph1982}, but here we take them to be the length of the line segment of projection $b_i$ going through pixel $x_j$, see \hypref{Fig.}{fig:pbp}. For our numerical experiments in \hypref{sections}{sec:findiff} and \ref{sec:numapp} we generate the matrix $R$ on the fly using the ASTRA toolbox \cite{Palenstijn2011, Aarle2015}.
\end{itemize}

The resulting linear system $Rx = b$, for $x\in\mbbR^n, b\in\mbbR^m$ and $R\in\mbbR^{m\times n}$, is a discrete version of the Radon transformation and can be used to solve the absorption based CT problem. However, in DPC-CT we use measurements of the derivative of the Radon transform. This results in a new linear system $DRx = b$, where $D$ is a finite difference matrix.

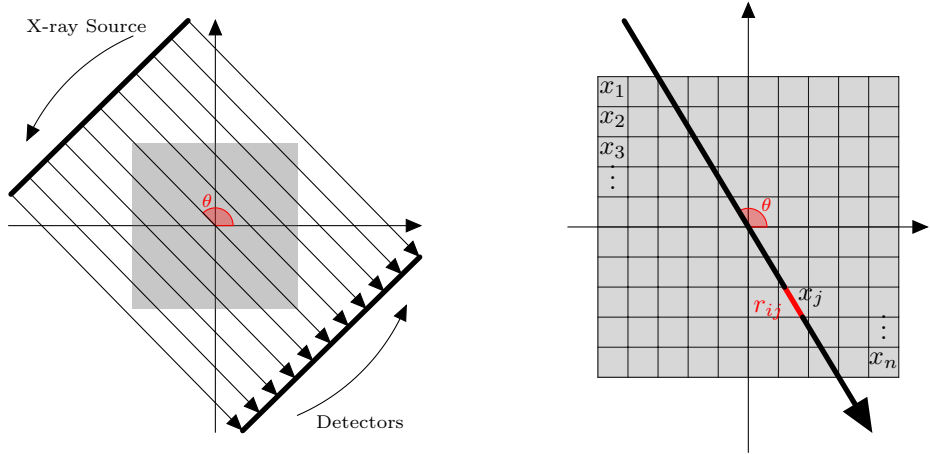
\begin{figure}
\begin{minipage}{0.45\linewidth}
\centering
%Tikz image: Parallel beam projection scheme
\definecolor{ffqqqq}{rgb}{1,0,0}
\definecolor{wqwqwq}{rgb}{0.38,0.38,0.38}
\begin{tikzpicture}[line cap=round,line join=round,>=triangle 45,x=1.0cm,y=1.0cm,scale=0.55]
\draw[->,color=black] (-5,0) -- (5,0);
\draw[->,color=black] (0,-5) -- (0,5);
\clip(-5.0,-6.0) rectangle (7.0,6.0);
\fill[color=wqwqwq,fill=wqwqwq,fill opacity=0.35] (-2,2) -- (-2,-2) -- (2,-2) -- (2,2) -- cycle;
\draw [shift={(0,0)},color=ffqqqq,fill=ffqqqq,fill opacity=0.35] (0,0) -- (0:0.43) arc (0:134.41:0.43) -- cycle;
\draw [line width=2pt] (-4.94,0.76)-- (-0.66,4.96);
\draw [line width=2pt] (4.94,-0.76)-- (0.66,-4.96);
\draw [shift={(0,0)}] plot[domain=1.98:2.71,variable=\t]({1*5*cos(\t r)+0*5*sin(\t r)},{0*5*cos(\t r)+1*5*sin(\t r)});
\draw [shift={(0,0)}] plot[domain=5.12:5.86,variable=\t]({1*5*cos(\t r)+0*5*sin(\t r)},{0*5*cos(\t r)+1*5*sin(\t r)});
\draw [->] (-4.94,0.76) -- (0.66,-4.96);
\draw [->] (-4.51,1.18) -- (1.08,-4.54);
\draw [->] (-4.09,1.6) -- (1.51,-4.12);
\draw [->] (-3.66,2.02) -- (1.94,-3.7);
\draw [->] (-3.23,2.44) -- (2.37,-3.28);
\draw [->] (-2.8,2.86) -- (2.8,-2.86);
\draw [->] (-2.37,3.28) -- (3.23,-2.44);
\draw [->] (-1.94,3.7) -- (3.66,-2.02);
\draw [->] (-1.51,4.12) -- (4.09,-1.6);
\draw [->] (-1.08,4.54) -- (4.51,-1.18);
\draw [->] (-0.66,4.96) -- (4.94,-0.76);
\draw [->] (-4.55,2.07) -- (-4.57,2.03);
\draw [->] (4.55,-2.07) -- (4.63,-1.9);
\begin{scriptsize}
\draw[color=black] (-3.11,4.8) node {X-ray Source};
\draw[color=ffqqqq] (-0.2,0.61) node {$\theta$};
\draw[color=black] (3.50,-4.75) node {Detectors};
\end{scriptsize}
\end{tikzpicture}
\end{minipage}\hspace*{\fill}%
\begin{minipage}{0.45\linewidth}
\centering
%TikZ image: determining weights
$$
\definecolor{ffqqqq}{rgb}{1,0,0}
\definecolor{uuuuuu}{rgb}{0.27,0.27,0.27}
\definecolor{wqwqwq}{rgb}{0.38,0.38,0.38}
\begin{tikzpicture}[line cap=round,line join=round,>=triangle 45,x=1.0cm,y=1.0cm,scale=0.2]
\draw[step=2.0cm] (-10,-10) grid (10,10);
\draw[->,color=black] (-12,0) -- (12,0);
\foreach \x in {-12,-10,-8,-6,-4,-2,2,4,6,8,10}
\draw[shift={(\x,0)},color=black] (0pt,-2pt);
\draw[->,color=black] (0,-15) -- (0,15);
\foreach \y in {-14,-12,-10,-8,-6,-4,-2,2,4,6,8,10,12,14}
\draw[shift={(0,\y)},color=black] (-2pt,0pt);
\clip(-12,-15) rectangle (12,15);
\fill[color=wqwqwq,fill=wqwqwq,fill opacity=0.25] (-10,10) -- (-10,-10) -- (10,-10) -- (10,10) -- cycle;
\draw [shift={(0,0)},color=ffqqqq,fill=ffqqqq,fill opacity=0.35] (0,0) -- (0:1.24) arc (0:121.06:1.24) -- cycle;
\draw (-9.0,8.0) node[anchor=south] {$x_1$};
\draw (-9.0,6.0) node[anchor=south] {$x_2$};
\draw (-9.0,4.0) node[anchor=south] {$x_3$};
\draw (-9.0,1.9) node[anchor=south] {$\vdots$};
\draw (9.0,.-8.1) node[anchor=south] {$\vdots$};
\draw (9.0,-10.0) node[anchor=south] {$x_n$};
\draw [line width=2pt,color=ffqqqq] (2.41,-4)-- (3.61,-6);
\draw [line width=2pt] (-8.25,13.71)-- (2.41,-4);
\draw [->,line width=2pt] (3.61,-6) -- (8.25,-13.71);
\draw (3.0,-5.5) node[anchor=east] {$\color{red}{r_{ij}}$};
\draw (2.7,-4.7) node[anchor=west] {$x_j$};
\begin{scriptsize}
\draw [fill=uuuuuu] (0,0) circle (1.5pt);
\draw[color=ffqqqq] (1.2,1.5) node {$\theta$};
\end{scriptsize}
\draw (-10, -10) -- (-10, 10);
\draw (-10,-10) -- (10, -10);
\end{tikzpicture}
$$
\end{minipage}
\caption{Overview of parallel beam projection \& definition of the weigts.}
\label{fig:pbp}
\end{figure}

Two natural choices for $D$ are the forward or the central finite difference approximation. If we set
\[
q(\theta, t) = \radon\delta(\theta, t) = \int_{L(\theta, t)}{\delta(\bfx)\diff\bfx},
\]
and take $h$ to be the distance between the detector cells, the forward and central difference scheme are respectively given by:
\begin{align*}
\frac{\partial q(\theta, t)}{\partial t}&=\frac{q(\theta, t+h)-q(\theta, t)}{h}+\mathcal{O}(h)\\
\intertext{or}
\frac{\partial q(\theta, t)}{\partial t}&=\frac{q(\theta, t+h)-q(\theta, t-h)}{2h}+\mathcal{O}(h^2).
\end{align*}
So the discrete values of $\alpha$ are approximated by one of the following schemes:
\begin{align*}
\alpha(\theta_i, t_j)&=\left.\frac{\partial q(\theta_i, t)}{\partial t}\right|_{t = t_j}=\frac{q(\theta_i, t_{j+1})-q(\theta_i, t_j)}{h}\\
\intertext{or}
\alpha(\theta_i, t_j)&=\left.\frac{\partial q(\theta_i, t)}{\partial t}\right|_{t=t_j}=\frac{q(\theta_i, t_{j+1})-q(\theta_i, t_{j-1})}{2h}.
\end{align*}
The matrices $D$ corresponding with these approximations are given below. Note that both are constructed using a Kronecker product because the derivative in \eqref{eq:alpha} is taken with respect to the detector space and independently of the projection angle. In \hypref{section}{sec:findiff} we will take a closer look at how each matrix affects the reconstruction.

\begin{itemize}
\item Forward difference scheme:

\[
D_f=I_l\otimes\underbrace{\begin{pmatrix}-1& 1\\
													   &-1& 1\\
														 &  &-1 & 1\\
														 &  &   &\ddots&\ddots\\
														 &      &  &  &-1& 1\\
														 &      &  &  &  &-1
						\end{pmatrix}}_{k\times k}\in\mbbR^{m\times m}.
\]
\item Central difference scheme:
\[
D_c=I_l\otimes\frac{1}{2}\underbrace{\begin{pmatrix}0& 1\\
													 -1& 0& 1\\
													   &-1& 0& 1\\
														 &  &\ddots&\ddots&\ddots\\
														 &  &      &-1& 0& 1\\
														 &  &      &  &-1& 0
						\end{pmatrix}}_{k\times k}\in\mbbR^{m\times m}.
\]
\end{itemize}

%%%%%%%%%%%%%%%%%%%%%%%% SECTION: Iterative Methods %%%%%%%%%%%%%%%%%%%%%%%%
\section{Iterative reconstruction with Krylov subspaces}\label{sec:artik}
One of the problems when dealing with inverse problems is that the forward problem still contains some model and approximation errors and that the right hand side $b$ typically contains small measurement errors, i.e.
\[
Ax = b = \widetilde{b} + e
\]
where $\widetilde{b} = A\widetilde{x}$ for the exact original image $\widetilde{x}$. When using an iterative method to solve these problems, the first few iterations usually seem to converge very nicely to $\widetilde{x}$. However, as the iterations continue, the noise starts to have an increasing effect on the approximations $x_k$ and the error with respect to $\widetilde{x}$ increases again. This effect is called semi-convergence \cite{hansen2010}. Therefore, some form of regularization method has to be performed in order to acquire a good approximation for $\widetilde{x}$. An often used method for large scale linear problems is Tikhonov regularization, which instead of solving $Ax = b$, solves the problem
\begin{equation}\label{eq:TikReg}
x_\lambda=\argmin_{x\in\mathcal{K}}\left\{\left\|Ax-b\right\|^2+\lambda\left\|L(x-x_0)\right\|^2\right\}
\end{equation}
for some subspace $\mathcal{K}\subseteq\mbbR^n$. Here, we will take $\mathcal{K}$ to be the Krylov subspace spanned by the matrix $A^TA\in\mbbR^{m\times n}$ and $A^Tr_0$ for the initial residue $r_0 = b - Ax_0\in\mbbR^m$. This method can be derived from the singular value decomposition of $A$ \cite{hansen2010}. Intuitively, the idea is to balance solving the original problem -- by minimizing the fidelity term $\left\|Ax-b\right\|^2$ -- and suppressing the noise  -- by minimizing the penalty term $\left\|L(x-x_0)\right\|^2$. This is done by choosing an appropriate value for the regularization parameter $\lambda$. If prior knowledge of the solution is available, the choice of a regularization matrix $L$ that incorporates these conditions in the solution can greatly improve the convergence of the method. Various variations exist to this standard formulation, e.g. replacing the penalizing term by $\left\|x\right\|^2_1$ or with the total variation operator $\text{TV}(x)$. However, the difficulty lies in finding an optimal value for $\lambda$.  If it is chosen too small, the method will still exhibit semi-convergence, and, if it is chosen too large, the solution will tend to be over-smooth.

There exist a number of different methods to find a value for $\lambda$ that is in some way suitable: the L-curve \cite{Calvetti2000, Calvetti2004}, the L-ribbon \cite{Calvetti1999}, generalized cross validation \cite{Chung2008, Golub1979}, the discrepancy principle \cite{Siltanen2012}, \ldots. In \cite{Gazzola2014_2, Gazzola2014} a generalized Arnoldi-Tikhonov method is proposed to solve \eqref{eq:TikReg} that incorporates an iterative scheme based on the discrepancy principle to find an optimal value for $\lambda$. However, since this method is based on the Arnoldi decomposition of a matrix, it only works for square matrices $A$. This will not be true for general inverse problem, e.g. in tomographic reconstructions, where the matrix $A$ is non square due to the number of detectors, projection angles and pixels (see section \ref{sec:disrec}). We therefore suggest an adaptation of this algorithm that uses a bidiagonal decomposition of the matrix $A$ instead of the Arnoldi decomposition. Since both decompositions are closely connected, the Tikhonov problem \eqref{eq:TikReg} can be solved in a similar way.

% ------------ SUBSECTION  ALGORITHM  ---------------------
\subsection{Generalized Bidiagonal-Tikhonov}
In \cite{golub1965}, it is shown that a matrix $A\in\mbbR^{m\times n}$ can be decomposed as
\begin{equation}\label{eq:BidiagDecomp}
A = U\begin{pmatrix}B\\0\end{pmatrix}V^T,
\end{equation}
with $U\in\mathbb{R}^{m\times m}$ and $V\in\mathbb{R}^{n\times n}$ orthonormal, $0\in\mathbb{R}^{(m-n)\times n}$ and $B\in\mathbb{R}^{(n+1)\times n}$ an lower bidiagonal matrix. We now denote the columns of $U$ and $V$ respectively as $u_1, \ldots, u_m\in\mathbb{R}^m$ and $v_1, \ldots, v_n\in\mathbb{R}^n$ and the matrices of the first $k$ columns as $U_k=(u_1, \ldots, u_k)$, $V_k=(v_1, \ldots, v_k)$ and $B_{k+1, k}$. If
\[
B=\begin{pmatrix}\alpha_1\\ \beta_2 & \alpha_2\\  & \beta_3 & \ddots\\  &  & \ddots & \alpha_{n}\\  &  &  & \beta_{n+1}\end{pmatrix}\in\mathbb{R}^{(n+1)\times n},
\]
then \eqref{eq:BidiagDecomp} implies the following relations:
\[
AV=U_nB\qquad\qquad A^TU_n=VB^T.
\]
Setting $\beta_1v_0=\alpha_{n+1}v_{n+1}=0$, this means that for all $k = 1, 2, \ldots, n$:
\begin{align*}
A^Tu_k&=\beta_kv_{k-1}+\alpha_kv_k\\
  Av_k&=\alpha_ku_k+\beta_{k+1}u_{k+1}.
\end{align*}
Starting with a given vector $u_1\in\mathbb{R}^m$ and $\left\|u_1\right\|=1$ we get the following recursion by Paige and Saunders called Bidiag1 \cite{Paige1982_2, Paige1982}:
\begin{equation}\label{eq:Bidiag1}
\left\{\begin{aligned}
r_k&=A^Tu_k-\beta_kv_{k-1},&\ \alpha_k&=\left\|r_k\right\|,&\ v_k&=\frac{1}{\alpha_k}r_k\\
p_k&=Av_k-\alpha_ku_k,&\ \beta_{k+1}&=\left\|p_k\right\|,&\ u_{k+1}&=\frac{1}{\beta_{k+1}}p_k
\end{aligned}\right.
\end{equation}
This process stops when either $\alpha_k = 0$ or $\beta_k = 0$. The former can only happen when $\text{rank}(A^T) < m$ and the latter implies that
\[
AV_{k-1}=U_{k-1}B_{k-1, k-1}\qquad\qquad A^TU_{k-1}=V_{k-1}B_{k-1, k-1}^T,
\]
meaning that the singular values of $B_{k, k-1}$ are a subset of the singular values of $A$ \cite{bjorck1996}. After $k$ steps of Bidiag1 with initial vector $u_1 = r_0 = b - Ax_0$, orthonormal matrices $V_k\in\mathbb{R}^{n\times k}$ and $U_{k+1}\in\mathbb{R}^{m\times (k+1)}$ and a lower bidiagonal matrix $B_{k+1,k}\in\mathbb{R}^{(k+1)\times k}$ have been constructed that satisfy
\[
AV_k=U_{k+1}B_{k+1, k}.
\]
Where the generalized Arnoldi-Tikhonov algorithm was intrinsically connected to the GMRES algorithm -- which is based on the Arnoldi decomposition -- the method we describe here will be linked to the LSQR algorithm -- which is based on the bidiagonal decomposition. In each iteration, the current approximation is given by
\[
x_{k, \lambda}\ =\ \argmin_{x_k\in\mathcal{K}_k}\left\{\left\|Ax_k-b\right\|^2+\lambda\left\|L(x_k-x_0)\right\|^2\right\},
\]
where $x_{k, \lambda}$ lies in the Krylov subspace $\mathcal{K}_k = \mathcal{K}_k(A^TA, A^Tr_0) = \spn{V_k}$. This means that the iterations can be written as $x_{k, \lambda} = x_0 + V_ky_{k, \lambda}$ for some $y_{k, \lambda}\in\mbbR^n$. It's therefore equivalent in each iteration to solve
\begin{align*}
y_{k, \lambda}\ =\ &\argmin_{y_k\in\mathbb{R}^k}\left\{\left\|A(x_0+V_ky_k)-b\right\|^2+\lambda\left\|LV_ky_k\right\|^2\right\}\\
=\ &\argmin_{y_k\in\mathbb{R}^k}\left\{\left\|AV_ky_k-r_0\right\|^2+\lambda\left\|LV_ky_k\right\|^2\right\}\\
=\ &\argmin_{y_k\in\mathbb{R}^k}\left\{\left\|U_{k+1}B_{k+1,k}y_k-r_0\right\|^2+\lambda\left\|LV_ky_k\right\|^2\right\}\\
=\ &\argmin_{y_k\in\mathbb{R}^k}\left\{\left\|B_{k+1, k}y_k-\|r_0\|e_1\right\|^2+\lambda\left\|LV_ky_k\right\|^2\right\}\\
=\ &\argmin_{y_k\in\mathbb{R}^k}\left\|\begin{pmatrix}B_{k+1, k}\\ \sqrt{\lambda}LV_k\end{pmatrix}y_k-\begin{pmatrix}\|r_0\|e_1\\ 0\end{pmatrix}\right\|^2\numberthis\label{eq:RegSys}.
\end{align*}
Note that in the next to last equation we used the orthogonality of the columns of $U_{k+1}$ and the fact that $u_1 = r_0 / \|r_0\|$ and if $L$ is the identity matrix, \eqref{eq:RegSys} can written without the matrix $V_k$. Furthermore, if the regularization term $\left\|L(x-x_0)\right\|^2$ where dropped, this would simply be the LSQR algorithm .

% ------------ SUBSECTION  UPDATE STEP  ---------------------
\subsection{Secant update step}
We now describe how the parameter $\lambda$ can be updated iteratively by using a approximation for the discrepancy in each iteration. This method was suggested in \cite{Gazzola2014_2, Gazzola2014} for the generalized Arnoldi-Tikhonov method. Due to the resemblance with the bidiagonal decomposition, analogous derivations can be made. The idea behind the discrepancy principle is to choose the regularization parameter $\lambda$ such that
\begin{equation}\label{eq:dp}
\|b - Ax_\lambda\|=\eta\varepsilon,
\end{equation}
with the error norm $\varepsilon = \|e\| = \|b - \widetilde{b}\|$ and a small tolerance value $\eta\gtrsim 1$. We now define the discrepancy as the residual in an iteration as a function of $\lambda$
\[
\phi_k(\lambda)=\left\|b-Ax_{k,\lambda}\right\|,
\]
and say that the discrepancy principle is satisfied when $\phi_k(\lambda)\leq\eta\varepsilon$. In order to find a good approximation for $\lambda$ in iteration $k$, we can perform one step of the secant method in order to approximate a solution for $\phi_k(\lambda) = \eta\varepsilon$. Assuming we have found a value $\lambda_{k-1}$ in the previous iteration, this can be done using the initial points $(0, \phi_k(0))$ and $(\lambda_{k-1}, \phi_k(\lambda_{k-1})$, see \hypref{Fig.}{fig:secant}. Using these points, the update formula for $\lambda$ in each iteration with the secant method is given by:
\begin{equation}\label{eq:RegParaTemp}
\lambda_k=\frac{\eta\varepsilon-\phi_k(0)}{\phi_k(\lambda_{k-1})-\phi_k(0)}\lambda_{k-1}.
\end{equation}
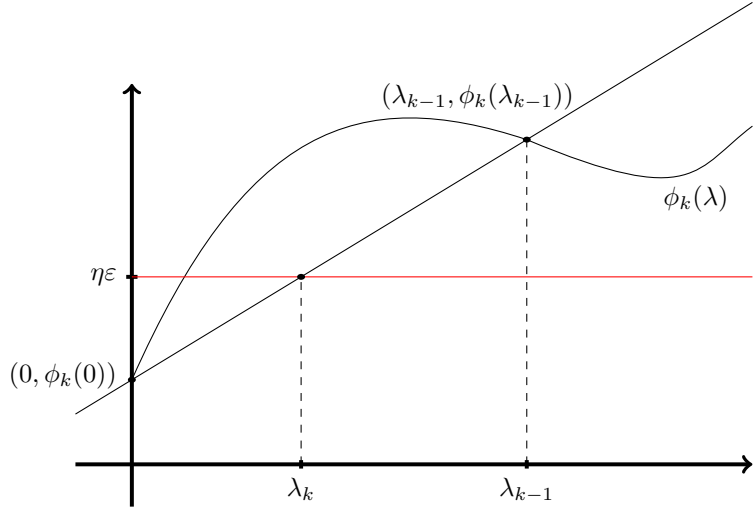
\begin{figure}
\centering
\begin{tikzpicture}[yscale = 0.75, scale = 0.75]

% Axis:
\draw[->, line width = 1.5] (0,-1) -- (0,9);
\draw[->, line width = 1.5] (-1,0) -- (11, 0);

% Ticks:
\draw[line width = 1.5] (-0.1, 4.4343) -- (0.1, 4.4343);
\draw[line width = 1.5] (3, -0.1) -- (3, 0.1);
\draw[line width = 1.5] (7, -0.1) -- (7, 0.1);

% Nodes:
\node[below] at (3, -0.1) {$\lambda_k$};
\node[below] at (7, -0.1) {$\lambda_{k-1}$};
\node at (10, 6.3) {$\phi_k(\lambda)$};
\node[left] at (-0.1, 4.4343) {$\eta\varepsilon$};
\node[left] at (8, 8.7) {$(\lambda_{k-1}, \phi_k(\lambda_{k-1}))$};
\node[left] at (-0.1, 2.1) {$(0, \phi_k(0))$};
%\node at (10, 9.5) {$r(\lambda)$};

% Discrepancy curve:
\draw (0,2) .. controls (2,8) and (4,9) .. (7, 7.68);
\draw (7, 7.68) .. controls (10,6) and (10,7) .. (11, 8);

% Lines:
\draw[color = red] (0,4.4343) -- (11, 4.4343);
\draw (-1, 1.1886) -- (11, 10.9257);
\draw[dashed] (3, 0) -- (3, 4.4343);
\draw[dashed] (7, 0) -- (7, 7.68);

% Dots:
\fill (7,7.68) circle (2pt);
\fill (0, 2) circle (2pt);
\fill (3, 4.4343) circle (2pt);

\end{tikzpicture}
\caption{Overview of the secant update step.}
\label{fig:secant}
\end{figure}
Therefore, in order to perform the secant update step, \eqref{eq:RegSys} has to be solved twice. Once with $\lambda = \lambda_{k-1}$, the current best estimate, and once with $\lambda = 0$. This means that $\phi_k(0)$ is simply the residual in the $k$-th iteration of the non-regularized LSQR method and can be found by solving
\begin{equation}\label{eq:LSQRsol}
\min_{y\in\mathbb{R}^k}\left\|c-B_{k+1,k}y\right\|,
\end{equation}
with $c=\left\|r_0\right\|e_1$. Because in the first few iterations the LSQR residual will be bigger than the tolerance value, i.e. $\phi_k(0)>\eta\varepsilon$, \eqref{eq:RegParaTemp} would give a negative value for the regularization parameter, therefore \eqref{eq:RegParaTemp} is replaced by
\begin{equation}\label{eq:RegPara}
\lambda_k=\left|\frac{\eta\varepsilon-\phi_k(0)}{\phi_k(\lambda_{k-1})-\phi_k(0)}\right|\lambda_{k-1}
\end{equation}
This is also why in the first few iterations the value for $\lambda$ typically becomes very large. Only after $\phi_k(0)$ becomes smaller than $\eta\varepsilon$, the value of the parameter will improve. This is why the initial value of $\lambda_0$ has little to no effect on the final solution. In our numerical experiments we will therefore take $\lambda_0 = 1$. We solve two linear systems in each iteration, namely \eqref{eq:RegSys} and \eqref{eq:LSQRsol}, but these are two sparse systems due to the bidiagonal structure of the matrix $B_{k+1, k}$ that is used in both systems. However, when the LSQR algorithm needs many iterations for its residual to go below the threshold value $\eta\varepsilon$, and hence for the regularization parameter to converge, these linear systems become larger and more compute time is required.

An overview of this method is given in \hypref{Algorithm}{alg:GBiT}. Here, an extra condition is added to the stop criterion in order to ensure the convergence of $\lambda_k$. This is done by demanding that the discrepancy principle is satisfied a number of times before quitting.

\begin{algorithm}
\caption{Generalized Bidiagonal-Tikhonov (GBiT)}
\label{alg:GBiT}
\begin{algorithmic}[1]
\STATE{\textbf{Input:} $A, b, x_0, \eta, \epsilon, \lambda_0, L$ and maxcounter}
\STATE{counter $= 0$}
\FOR{$k = 1, \ldots,$ maxIter}
	\STATE{\begin{enumerate}
	\item Update the bidiagonal decomposition by \eqref{eq:Bidiag1} by calculating $v_k$ and $u_{k+1}$ with reorthogonalization.
	\item Solve \eqref{eq:LSQRsol} and calculate $\phi_k(0)$.
	\item Solve \eqref{eq:RegSys} and calculate $\phi_k(\lambda_{k-1})$.
	\item Calculate $\lambda_k$ by \eqref{eq:RegPara}.
	\end{enumerate}}
	\IF{$\phi_k(\lambda_{k-1}) < \eta\varepsilon$}
		\STATE{counter$++$}
		\IF{maxcounter $<$ counter}
			\STATE{\textbf{break}}
		\ENDIF
	\ENDIF
\ENDFOR
\STATE{$x_k = x_0 + V_ky_k$}
\end{algorithmic}
\end{algorithm}

% ------------ SUBSECTION  UPDATE STEP  ---------------------
\subsection{Alternative update scheme}\label{sec:altup}
A downside of this iterative scheme for the regularization parameter, and of the discrepancy principle in general, is that a fairly good estimate of the noise norm $\varepsilon$ is necessary. If it too small, the method will still exhibit semi-convergence, and if it is too big, the solution will be over-smoothed. This can be solved by changing the update formula \eqref{eq:RegPara} to
\begin{equation}\label{eq:RegPara2}
\lambda_k=\left|\frac{\eta\phi_{k-1}(0)-\phi_k(0)}{\phi_k(\lambda_{k-1})-\phi_k(0)}\right|\lambda_{k-1}.
\end{equation}
as was suggested in \cite{Gazzola2014_3}. The motivation for this adjustment is that the non-regularized residual of LSQR tends to stabilizes slightly below the the noise norm, i.e. $\phi_{k-1}(0)\approx\varepsilon$. The only problem with \eqref{eq:RegPara2} is that it may prevent the discrepancy principle from being satisfied. The reason for this is that when the approximation for the regularization parameter starts to stabilize, the following will hold:
\[\begin{aligned}
\lambda_k\approx\lambda_{k-1}&\ \Rightarrow\ \phi_k(\lambda_{k-1})-\phi_k(0)\approx\eta\phi_{k-1}(0)-\phi_k(0)\\
														 &\ \Rightarrow\ \phi_k(\lambda_{k-1})\approx\eta\phi_{k-1}(0)
\end{aligned}\]
Therefore, we propose to take a slightly bigger tolerance value for the discrepancy principle than for the update step. For example, by taking its value $1\%$ bigger, the condition on line 5 of \hypref{Algorithm}{alg:GBiT} becomes
\[
\phi_k(\lambda_{k-1}) < 1.01\eta\phi_{k-1}(0).
\]

%%%%%%%%%%%%%%%%%%%%%%%% SECTION: Finite difference approximation %%%%%%%%%%%%%%%%%%%%%%%%
\section{Finite difference approximation}\label{sec:findiff}
In \hypref{section}{sec:disrec}, we introduced the finite difference operators $D_f$ and $D_c$. In this section, we will examine how these affect the reconstruction by studying a test problem.  The following proposition summarizes some properties of the $D_f$ and $D_c$ which can easily be verified using induction and the properties of the kronecker product.

\begin{proposition}\label{prop1}
\begin{enumerate}[(i)]
\item For $D_f\in\mbbR^{m\times m}$ the following holds:
\[
\det(D_f) = (-1)^m\qquad\text{and}\qquad\Null(D_f) = \{0\}.
\]
\item For $D_c^m\in\mbbR^{m\times m}$ the following holds:
\[
\det(D_c^m) = \left\{\begin{aligned} &\frac{1}{2^m}&\quad &m\text{ even}\\
																		 &\ 0&\quad &m\text{ odd.}
																		\end{aligned}\right.
\]
and
\[
\Null(D_c^m) = \left\{\begin{aligned} &\ \{0\}&\quad &m\text{ even}\\
																		 &\left\langle(1, 0, 1, 0, \ldots, 1, 0)^T\right\rangle&\quad &m\text{ odd}
																		\end{aligned}\right.
																		\]
\end{enumerate}
\end{proposition}
\noindent Note that the inverse of $D_f$ is given by the kronecker product of $I_l$ and a backsubstitution matrix:
\[
D_f^{-1} = I_l\otimes\underbrace{\begin{pmatrix} -1 & -1 & \ldots & -1\\
																	   & -1 & \ldots & -1\\
																		 &    & \ddots & \vdots\\
																		 &    &    &-1
						\end{pmatrix}}_{k\times k}\in\mbbR^{m\times m}.
\]

% ------------ SUBSECTION: Test problem ---------------------
\subsection{Test problem}\label{sec:testp}
In order to study the effect of the finite difference operator on the reconstruction, we take the $256\times 256$ Shepp-Logan phantom in MATLAB, see \hypref{Fig.}{fig:testp1:aa} and look at a single projection on $256$ detectors (angle $\theta = \pi / 2$). We take $x\in\mbbR^{256^2}$ to be the vector that represents the image and $y = Rx\in\mbbR^{256}$ its Radon transformation. The DPC-CT measurements are now given by $b = Dy$. In order to generate this data and avoid an inverse crime (see \cite{Siltanen2012}), we generated a weighted version of $b$ by taking
\begin{equation}\label{eq:modelerror}
\left\{\begin{aligned}
b_f &= (1 - \omega) D_fy + \omega D_cy\\
b_c &= \omega D_fy + (1 - \omega) D_cy
\end{aligned}\right.
\end{equation}
for $\omega = 0.2$. For this dataset this yields a model error of approximately $10\%$. The data, the measurements and the reconstructions can be seen in \hypref{Fig.}{fig:testp1}. The reconstructions in \hypref{Fig.}{fig:testp1:c} and \hypref{Fig.}{fig:testp1:d} were made with $256$ iterations of the LSQR algorithm. The results demonstrate that the central difference model is much more sensitive to the model errors we introduced by generating the data on a finer grid. Note that in this case both $D_f$ and $D_c$ are invertible, so the LSQR algorithm converged to $D_f^{-1}b$ and $D_c^{-1}b$ after the 256 iterations.

\begin{figure}
\subfloat[]{\raisebox{4mm}{\includegraphics[width=0.23\linewidth]{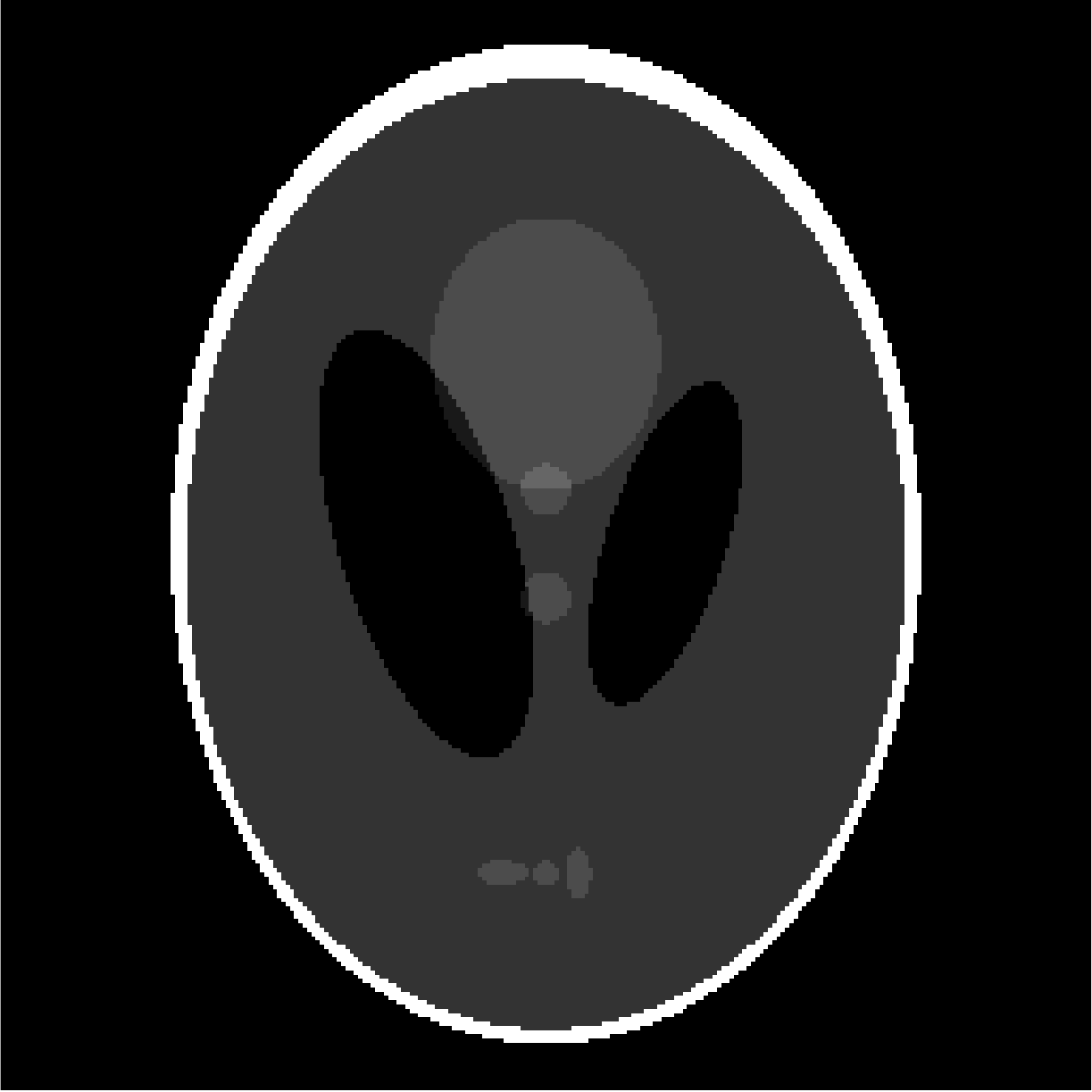}}\label{fig:testp1:aa}}\quad
\centering
\subfloat[]{\includegraphics[width=0.33\linewidth]{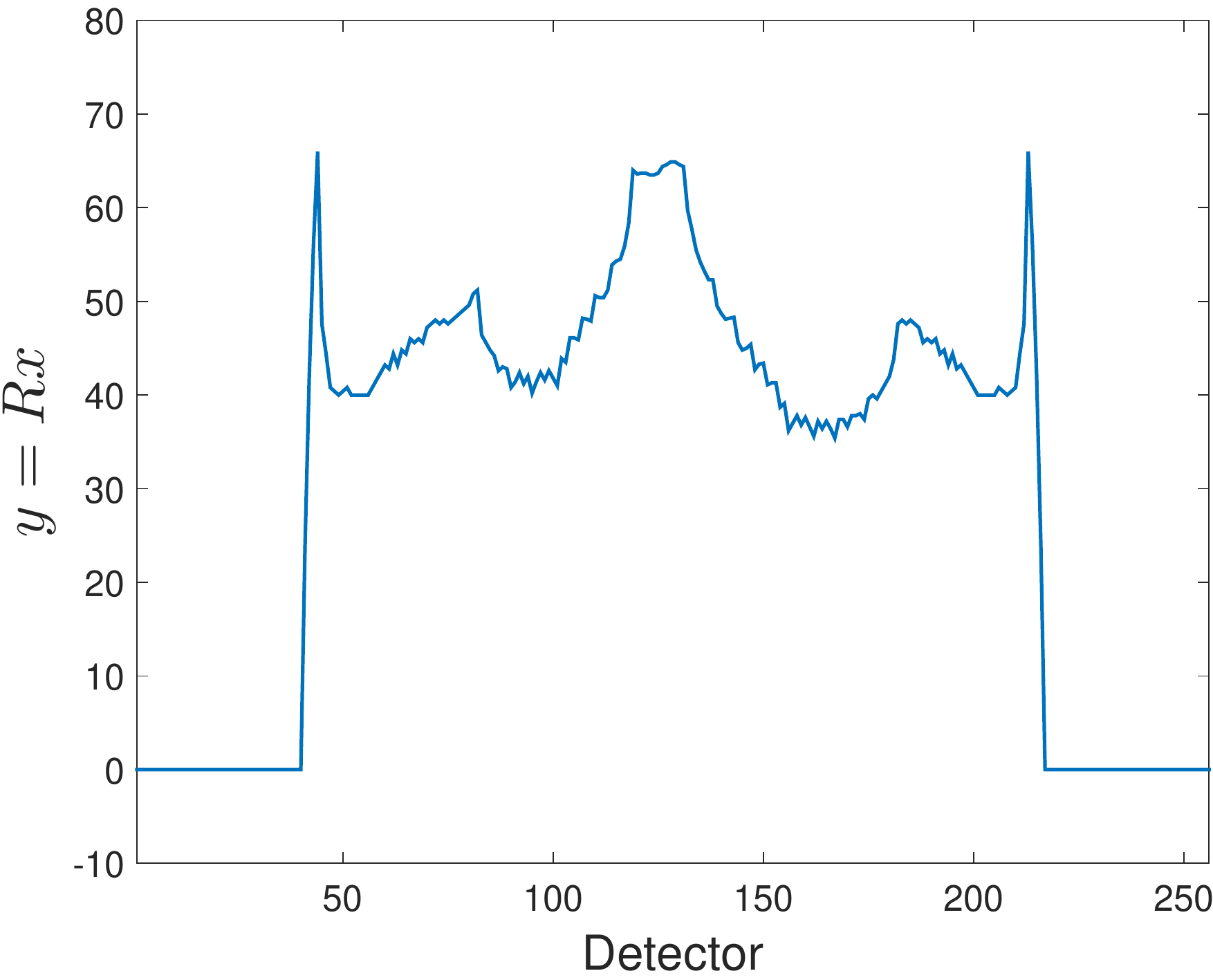}\label{fig:testp1:a}}\quad%\hspace*{\fill}
\subfloat[]{\includegraphics[width=0.33\linewidth]{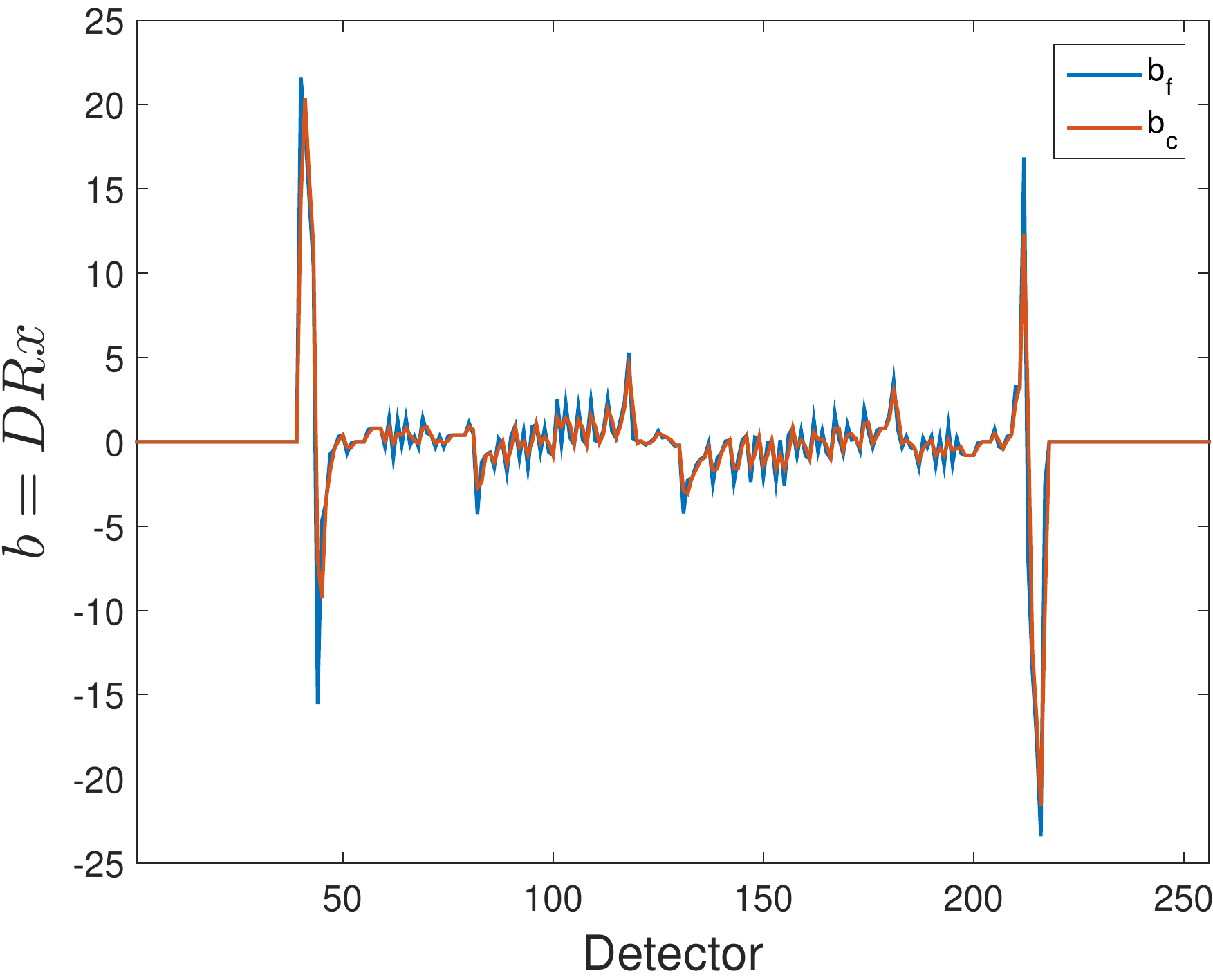}\label{fig:testp1:b}}\\%\hspace*{\fill}
\subfloat[]{\includegraphics[width=0.33\linewidth]{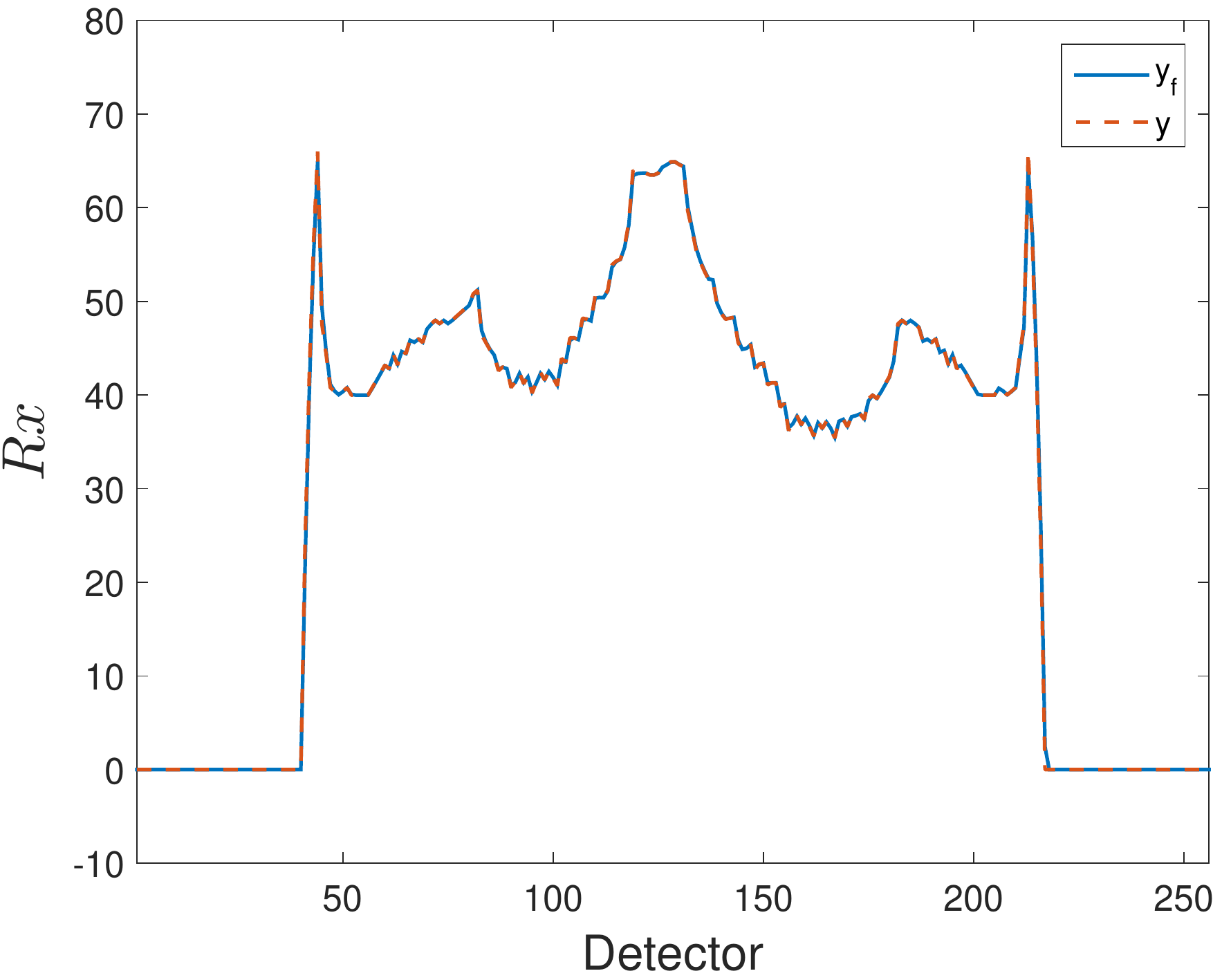}\label{fig:testp1:c}}\qquad%\hspace*{\fill}
\subfloat[]{\includegraphics[width=0.33\linewidth]{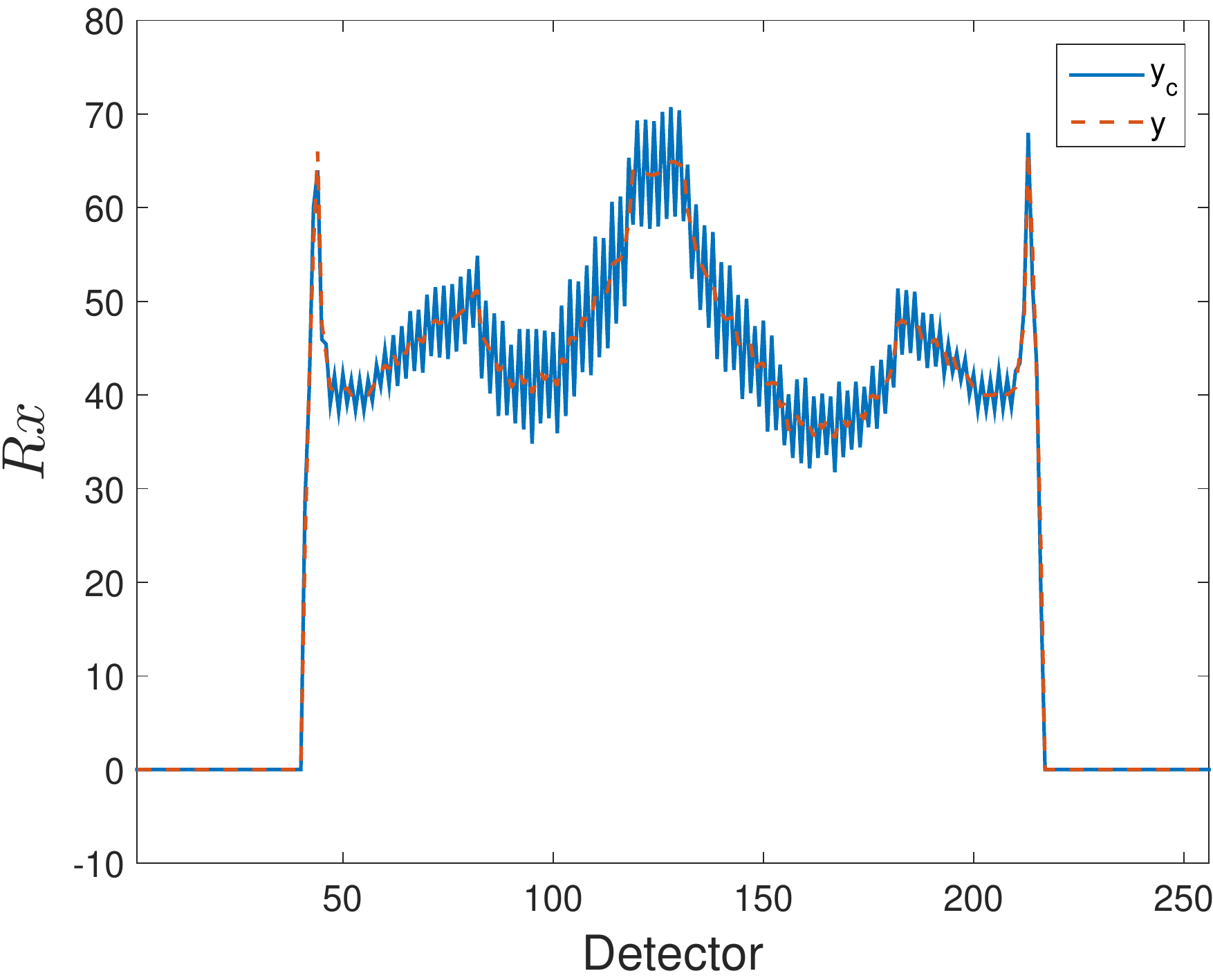}\label{fig:testp1:d}}
\caption{Data, measurements and reconstructions with model errors. \protect\subref{fig:testp1:aa} Shepp-Logan phantom. \protect\subref{fig:testp1:a} Radon projection of $x$ under a single projection angle ($\theta = \pi / 2$). \protect\subref{fig:testp1:b} The measurements $b$. \protect\subref{fig:testp1:c} LSQR reconstruction using the forward difference model. \protect\subref{fig:testp1:d} LSQR reconstruction using the central difference model.}
\label{fig:testp1}
\end{figure}

If, on the other hand, we assume that there are no model errors -- i.e. calculate $b$ without the weighting step -- but add noise instead, similar results can be observed. This can be seen in \hypref{Fig.}{fig:testp2} and \hypref{Fig.}{fig:testp2_reg}, where we added $10\%$ noise by generating a random normal vector $e\sim\mathcal{N}(0, 1)$ and setting
\begin{equation}\label{eq:addnoise}
b = \widetilde{b} + 0.10\frac{\|b\|}{\|e\|}e,
\end{equation}
for either the forward or the central difference model. Here, we also calculated the regularized solution with \hyperref[alg:GBiT]{GBiT} (using the maxcounter variable to perform 256 iterations).

 When we compare the different solutions, we see that the central difference model is, again, more sensitive to the errors in the measurements. Furthermore, as can be seen in \hypref{Fig.}{fig:testp2:c} and \hypref{Fig.}{fig:testp2:d}, regularization hardly improves the oscillatory nature of the solution of the central difference model. Also, because we assumed no prior knowledge of the solution and used $x_0 = 0$ and $L = I$, the Tikhonov solution given by \eqref{eq:TikReg} will be centred around $x_0$. This is because adding a constant function does not change the derivative and although these functions are not technically in the null space of either matrix, see \hypref{proposition}{prop1}, they only fail at one or both of the boundary points. This will have little to no weight in the regularization.

\begin{figure}
\centering
\subfloat[]{\includegraphics[width=0.33\linewidth]{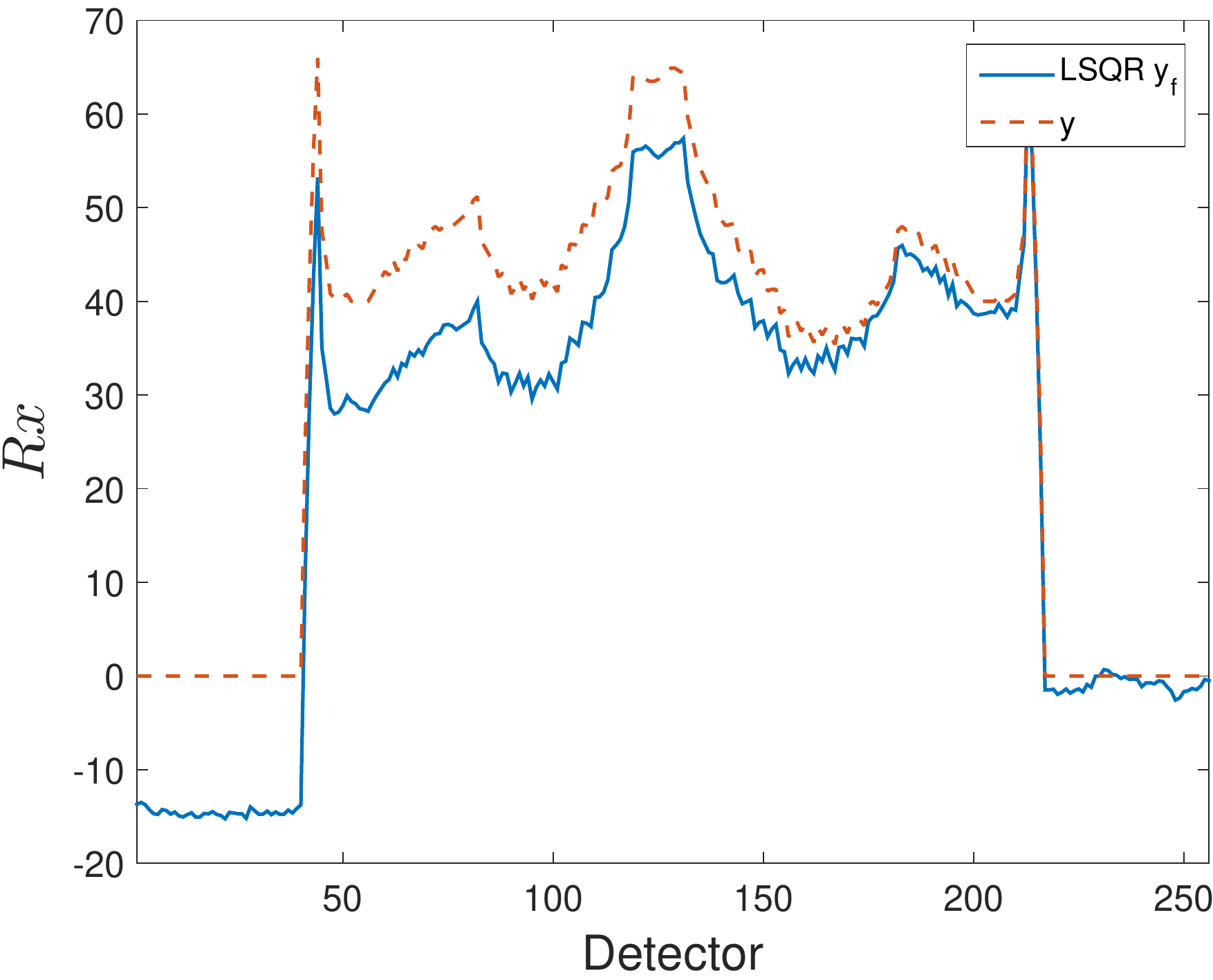}\label{fig:testp2:a}}\qquad%\hspace*{\fill} 
\subfloat[]{\includegraphics[width=0.33\linewidth]{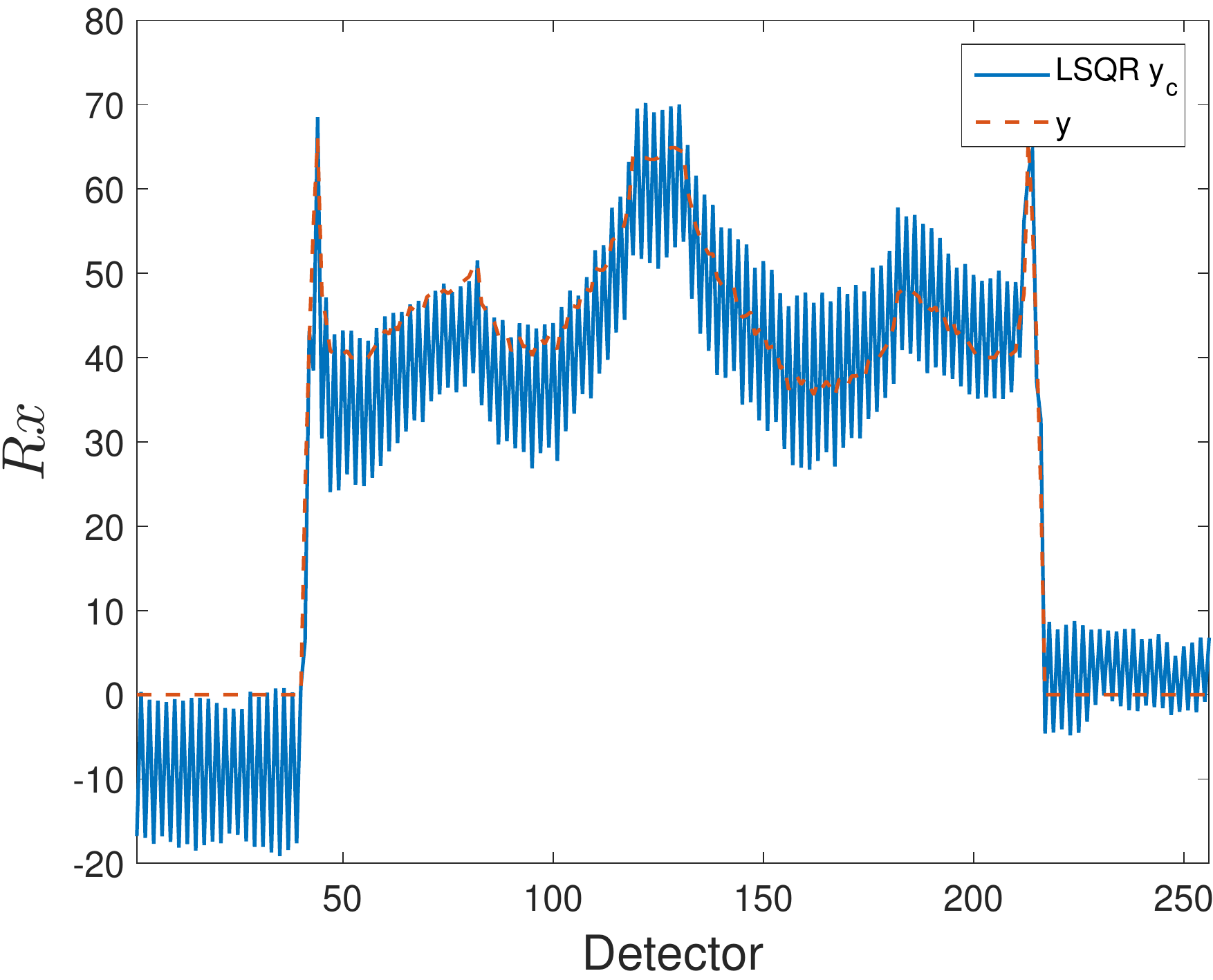}\label{fig:testp2:b}}\\%\hspace*{\fill}
\subfloat[]{\includegraphics[width=0.33\linewidth]{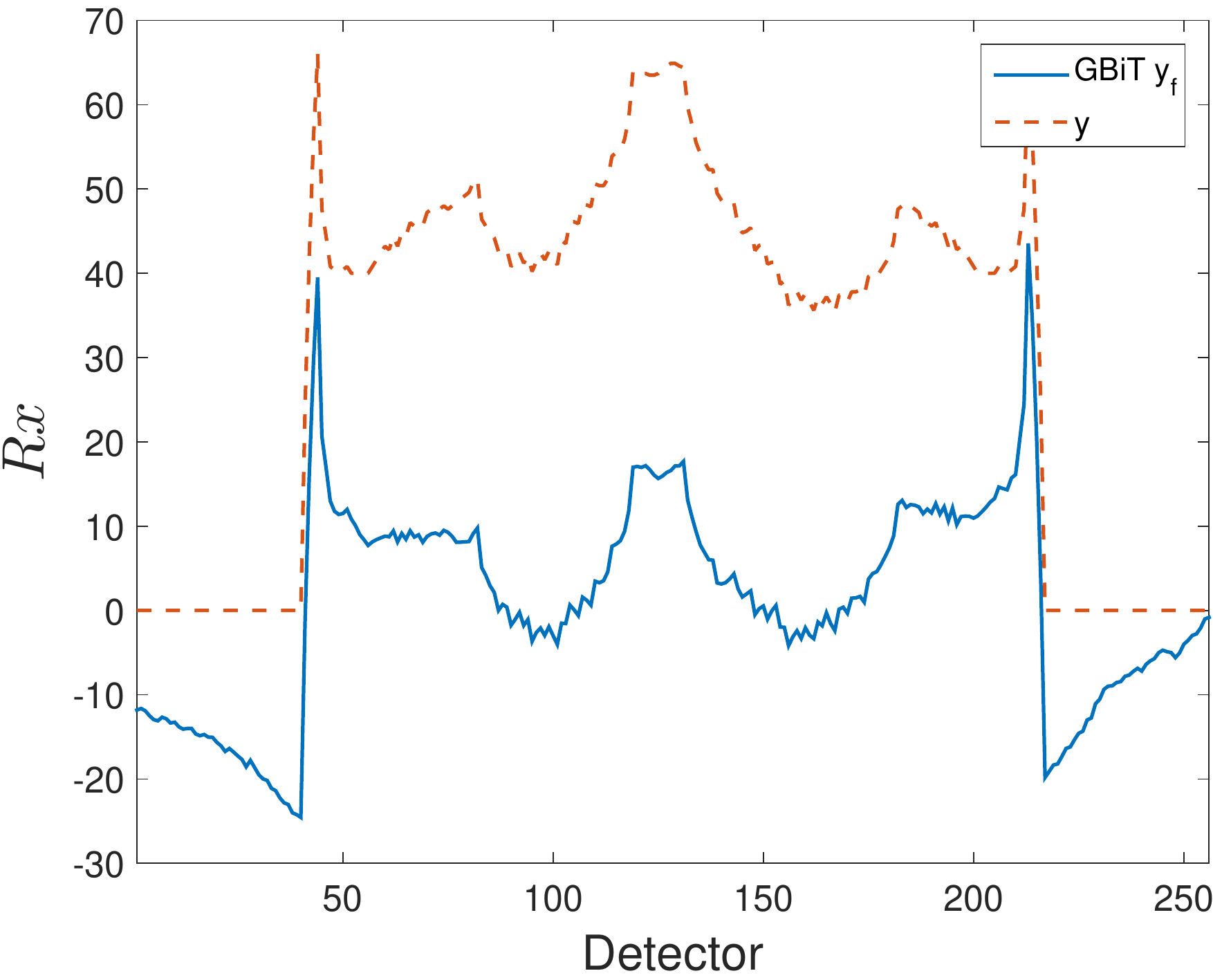}\label{fig:testp2:c}}\qquad%\hspace*{\fill}
\subfloat[]{\includegraphics[width=0.33\linewidth]{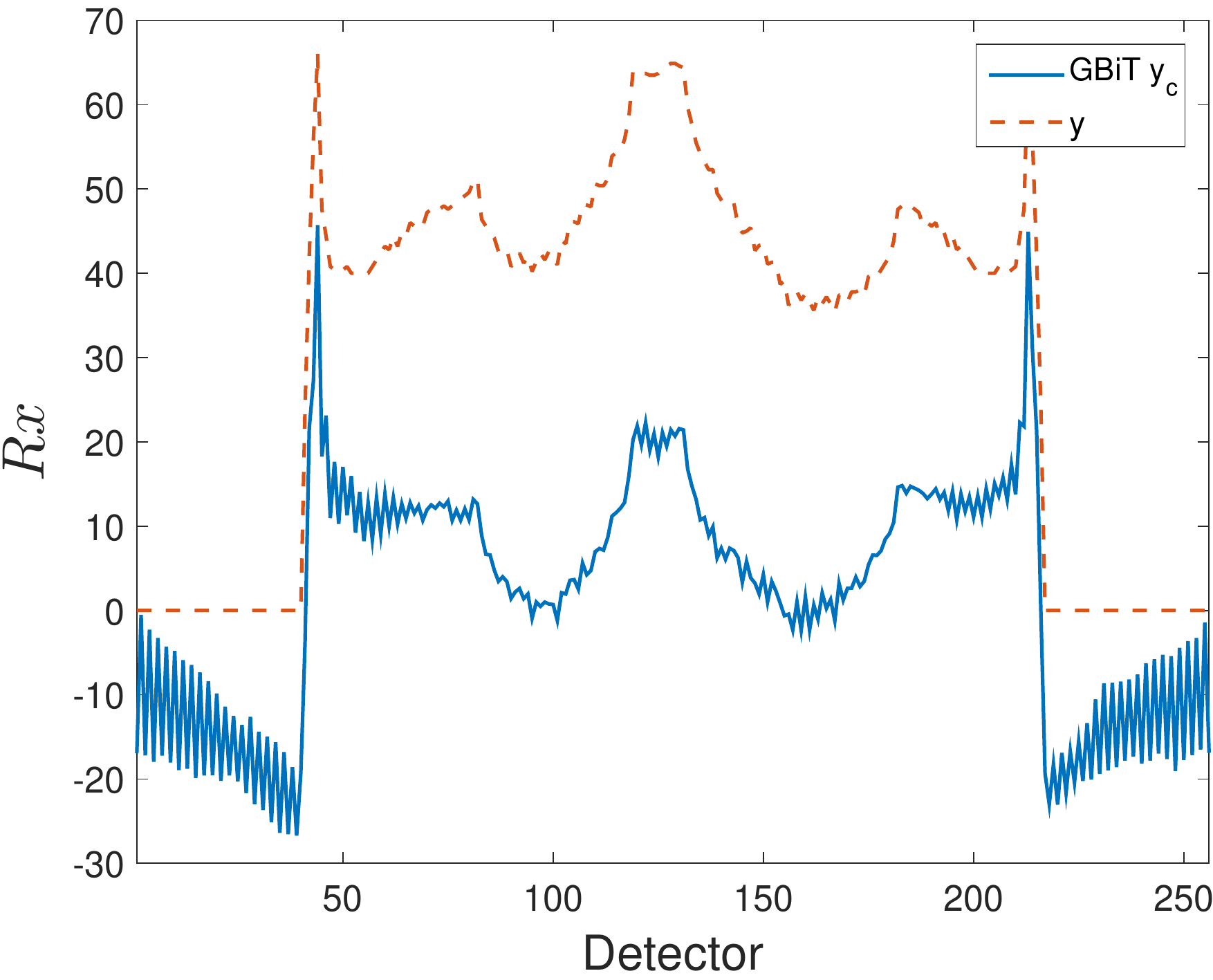}\label{fig:testp2:d}}
\caption{Reconstructions when noise is present. \protect\subref{fig:testp2:a} LSQR reconstruction using the forward difference model. \protect\subref{fig:testp2:b} LSQR reconstruction using the central difference model. \protect\subref{fig:testp2:c} GBiT reconstruction using the forward difference model. \protect\subref{fig:testp2:d} GBiT reconstruction using the central difference model.}
\label{fig:testp2}
\end{figure}

\begin{figure}
\centering
\subfloat[]{\includegraphics[width=0.3\linewidth]{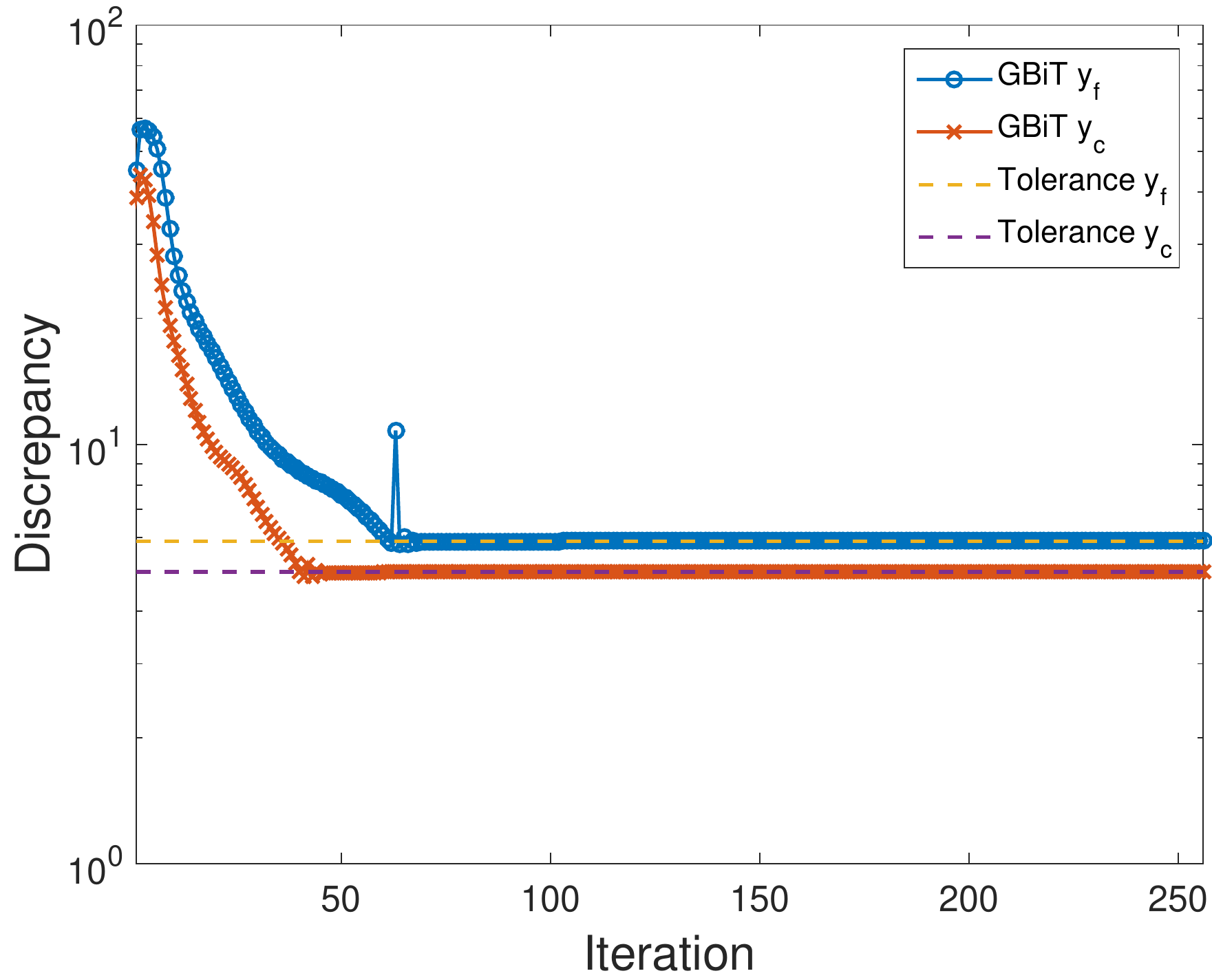}\label{fig:testp2_reg:a}}\quad
\subfloat[]{\includegraphics[width=0.3\linewidth]{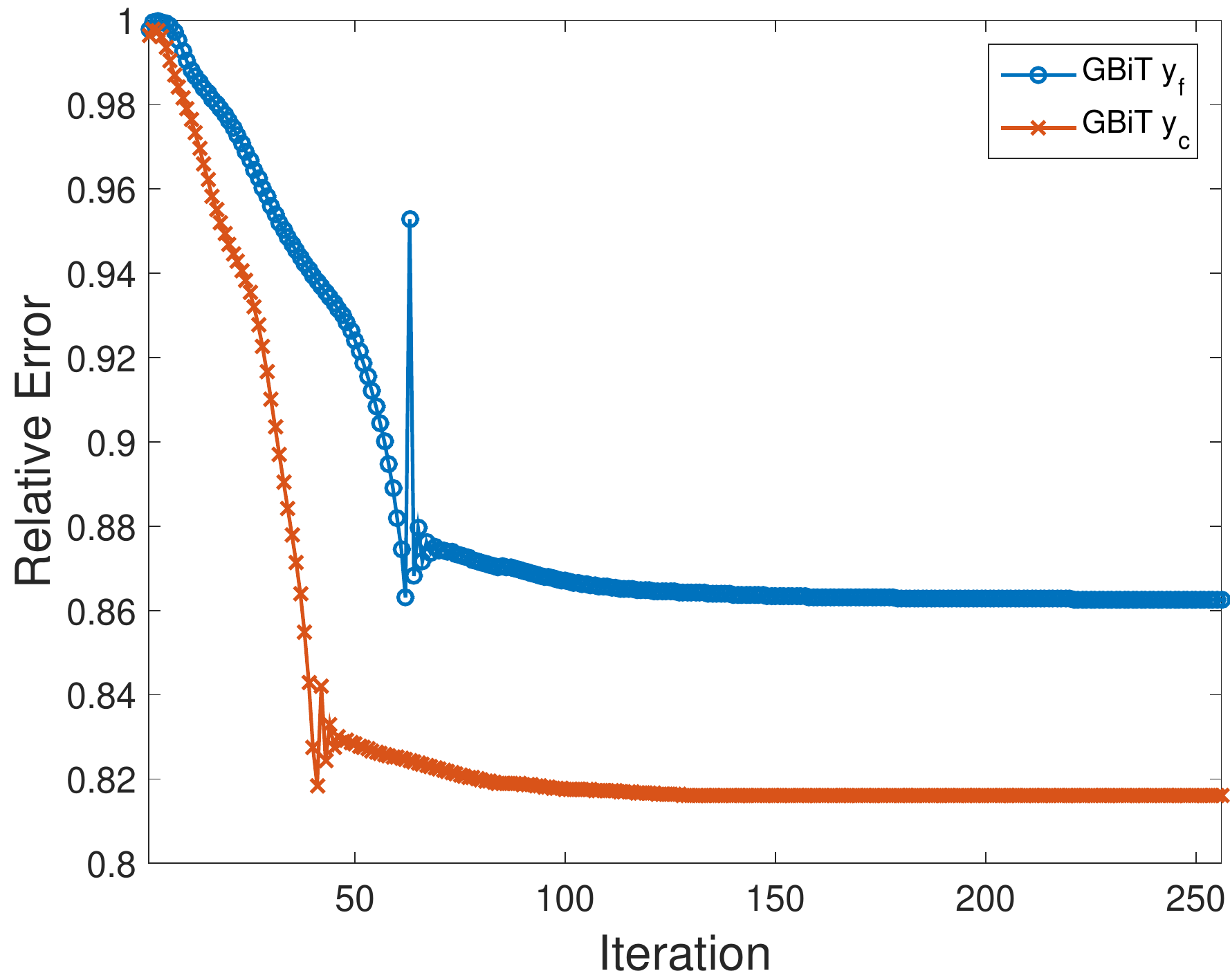}\label{fig:testp2_reg:b}}\quad
\subfloat[]{\includegraphics[width=0.3\linewidth]{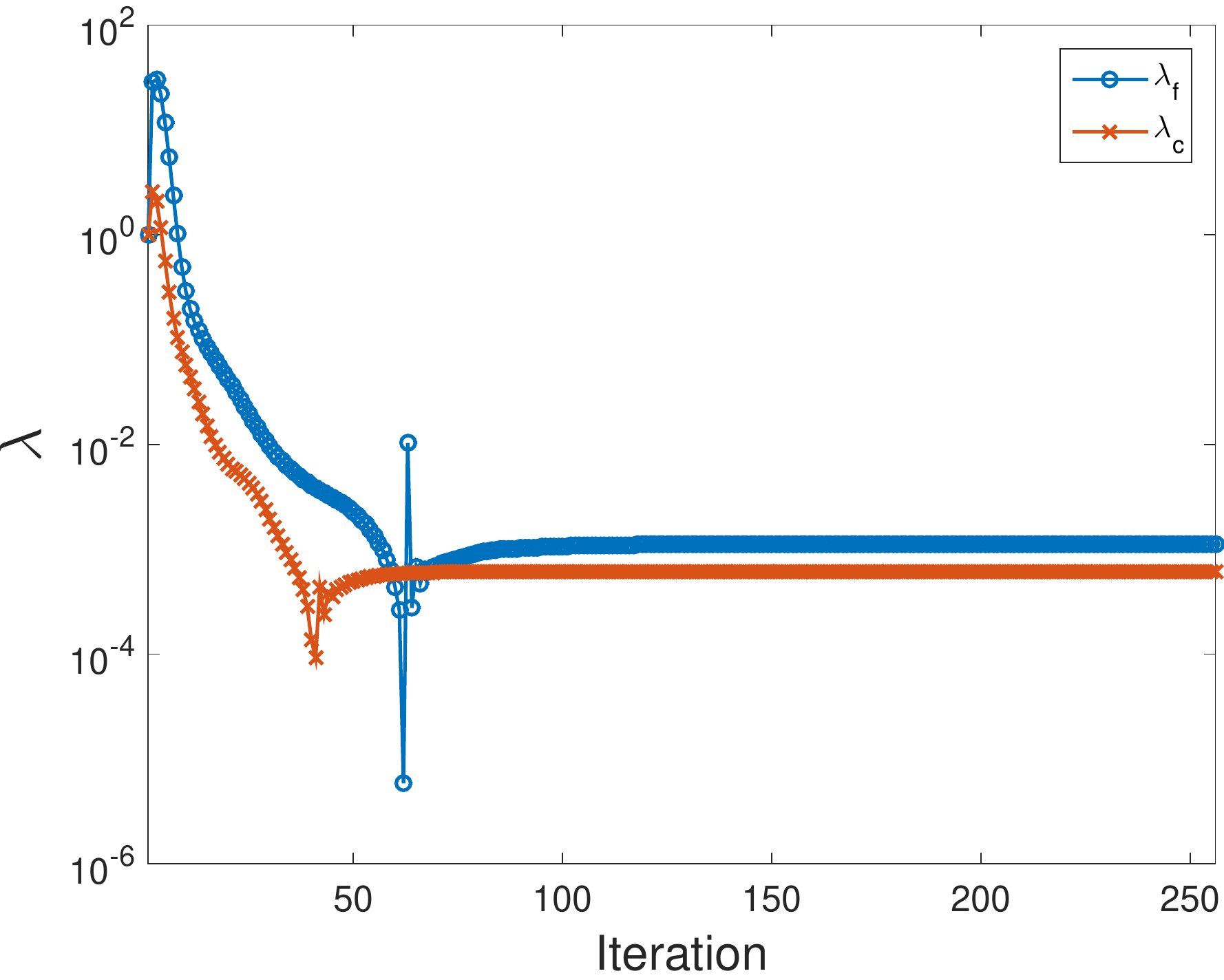}\label{fig:testp2_reg:c}}
\caption{Overview of the discrepancy \protect\subref{fig:testp2_reg:a}, relative error \protect\subref{fig:testp2_reg:b} and regularization parameter \protect\subref{fig:testp2_reg:c}
in each iteration of GBiT for the reconstructions from \hypref{Fig.}{fig:testp2:c} and \hypref{Fig.}{fig:testp2:d}. For the forward and the central difference model the discrepancy principle is satisfied after respectively $80$ and $43$ iterations. The oscillations in all three plots at these iterations motivates performing some extra iterations after the discrepancy principle is satisfied for the first time.}
\label{fig:testp2_reg} 
\end{figure}

Why the forward difference model gives much better reconstructions is explained by the inverse of the operators. If we write $b = \widetilde{b} + e$ for some error caused by either the model, the measurements or both, the solution of $D_fy = b$ is given by a simple back substitution:
\begin{equation}\label{eq:forwsmoothing}
y_i = - \sum_{j = i} ^n b_j = -\sum_{j = i}^n \widetilde{b}_j - \sum_{j=i}^n e_j\qquad(\forall i = 1, \ldots, m).
\end{equation}
So  the noise partially cancels itself out. This is seen in \hypref{Fig.}{fig:testp2:a} where the difference between the exact solution $y$ and the approximation $y_f$ does not increase monotonically from right to left. Important to note here is the assumption that the noise is centred around the measurements, i.e. $\mathbb{E}(b) = \widetilde{b}$ or $\mathbb{E}(e) = 0$ for whatever distribution is used to model the noise. If this is not the case, e.g. the noise has a constant offset w.r.t the measurements, this offset will add up in the back substitution. 

This is illustrated in \hypref{Fig.}{fig:testp3}, which is the same problem from \hypref{Fig.}{fig:testp2}, but where the noise was given a constant offset of 5. Here, we observe that the regularization suppresses this offset because -- as stated before -- it tries to centre the solution around $x_0 = 0$. From the reconstructions, however, we see that the central difference model does not have the same smoothing effect on the solution. This is because when $D_c$ is invertible (as is the case here), its inverse has a completely different structure\footnote{If $L$ is a strictly lower bidiagonal matrix consisting of alternating zeros and twos, this matrix can be seen as the sum of $L$ and $-L^T$.}: 
\[
D_c^{-1} = I_l\otimes\underbrace{\left(\begin{array}{rrrrrrr}0&-2&0&-2&0&-2&\cdots\\
\ 2&\ 0&\ 0&\ 0&\ 0&\ 0&\cdots\\
0&0&0&-2&0&-2&\cdots\\
2&0&2&0&0&0&\cdots\\
0&0&0&0&0&-2&\cdots\\
2&0&2&0&2&0&\cdots\\
\vdots&\vdots&\vdots&\vdots&\vdots&\vdots&\ddots\end{array}\right)}_{k\times k}\in\mbbR^{m\times m}
\]
If we write down the analogous formulation to \eqref{eq:forwsmoothing}, a first observation is that, because  the matrix consists of 0 and 2's, the noise components is multiplied by 2 as well. This implies that the part of the error that is not cancelled out doubles. Furthermore, since the rows have alternating plus and minus signs and a highly varying number of non-zero entries at alternating positions, the noise cancellation is less consistent compared to the forward difference model. This causes the heavy oscillations present in the non-regularized solutions.

\begin{figure}
\centering
\subfloat[]{\includegraphics[width=0.33\linewidth]{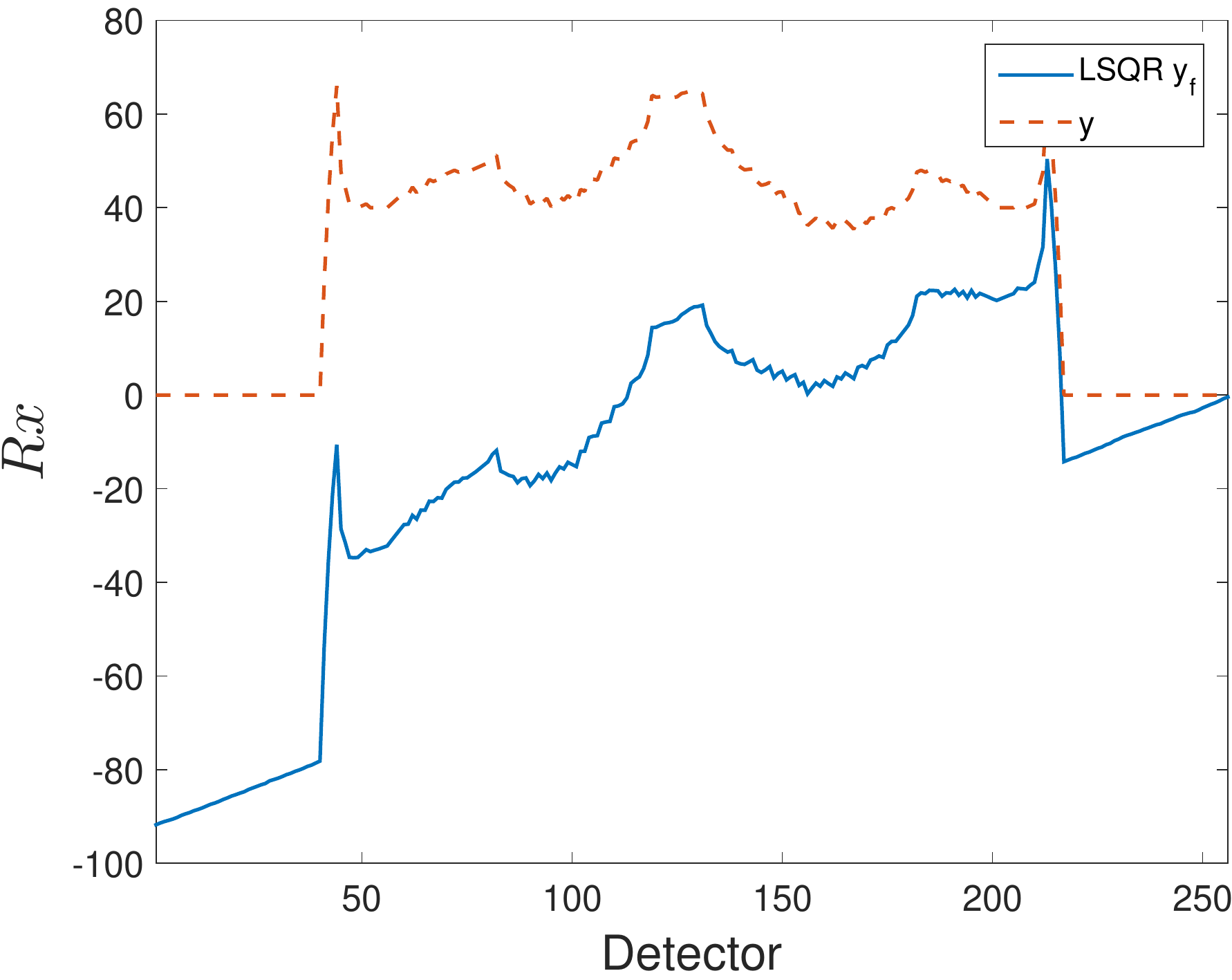}\label{fig:testp3:a}}\qquad%\hspace*{\fill}
\subfloat[]{\includegraphics[width=0.33\linewidth]{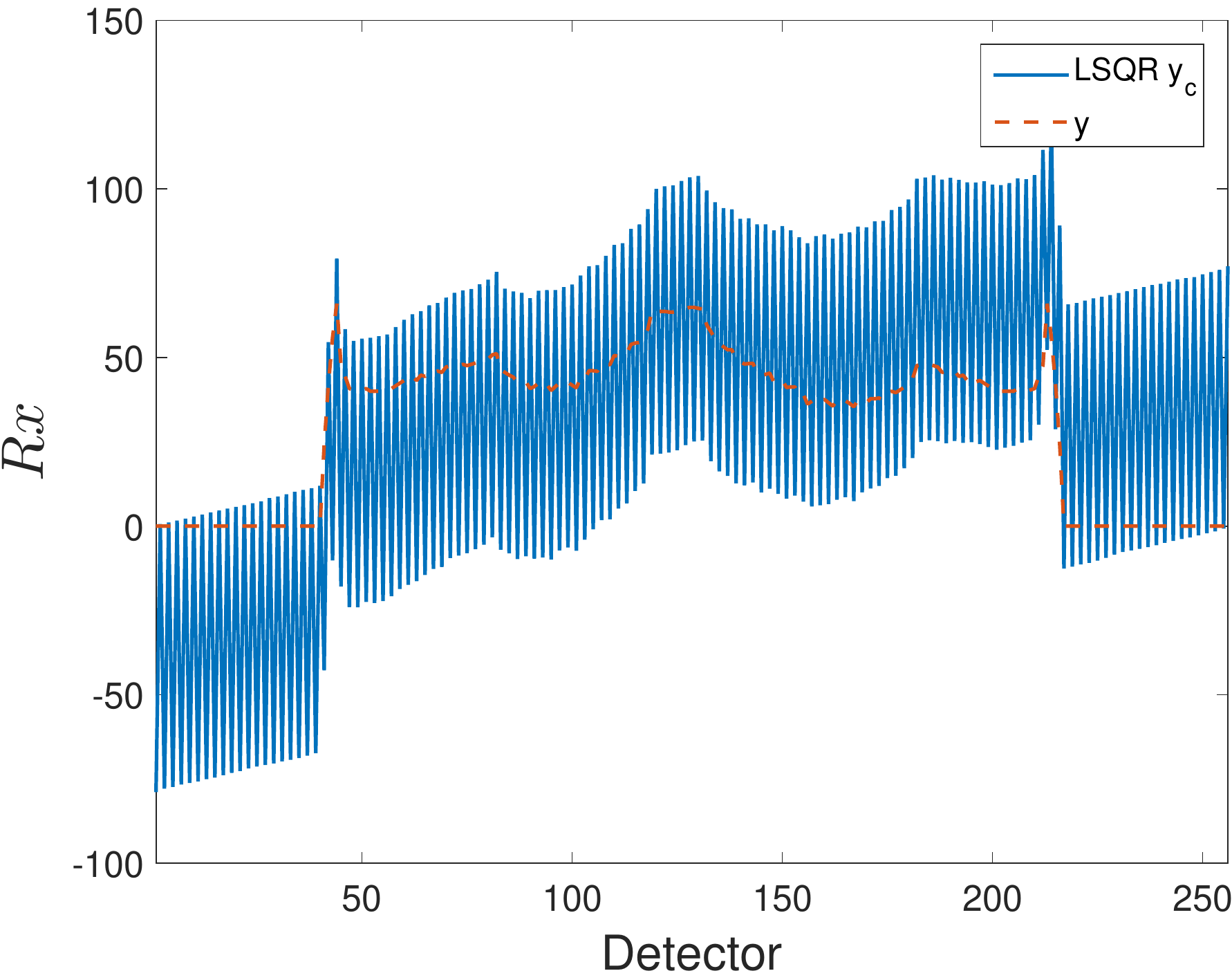}\label{fig:testp3:b}}\\%\hspace*{\fill}
\subfloat[]{\includegraphics[width=0.33\linewidth]{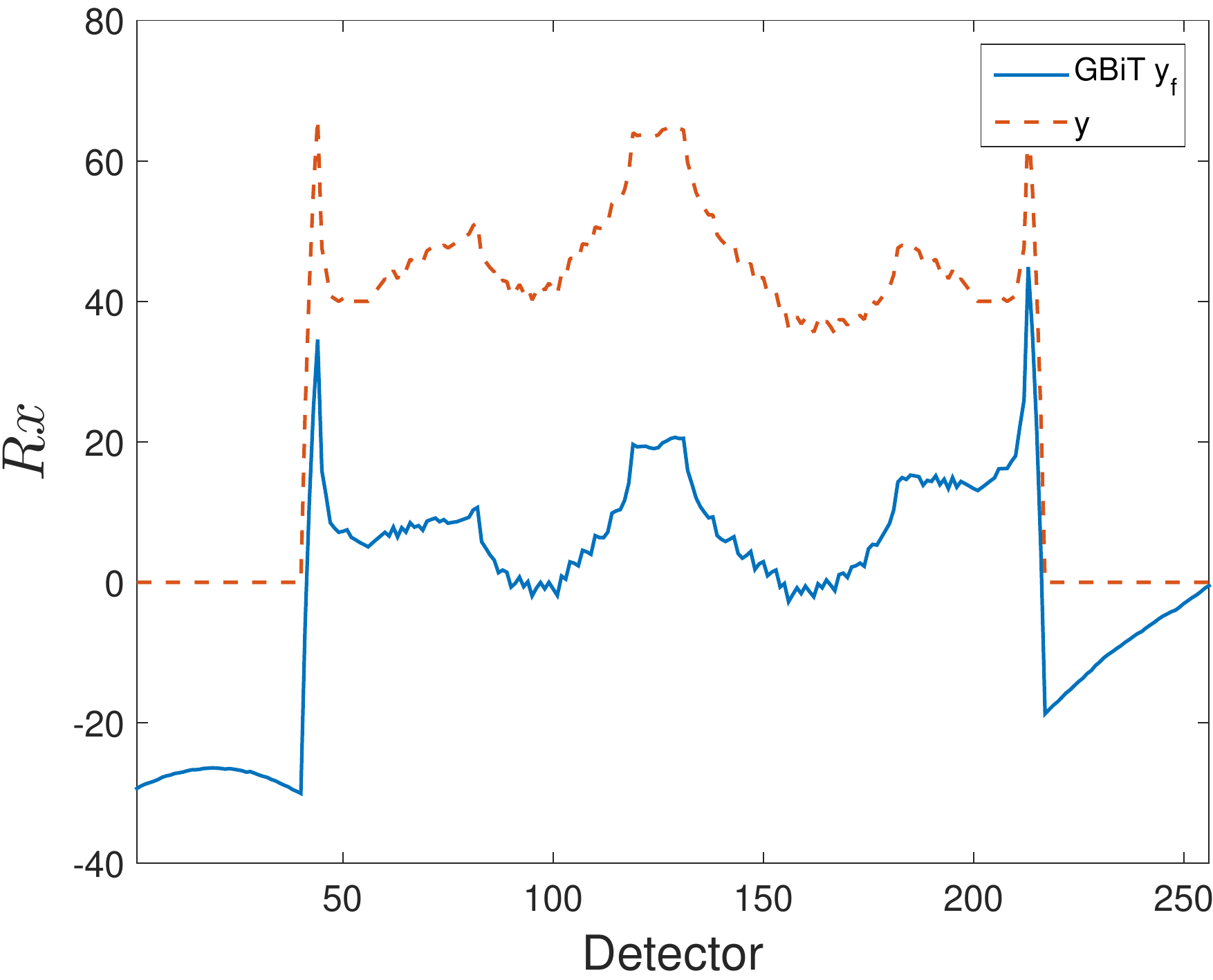}\label{fig:testp3:c}}\qquad%\hspace*{\fill}
\subfloat[]{\includegraphics[width=0.33\linewidth]{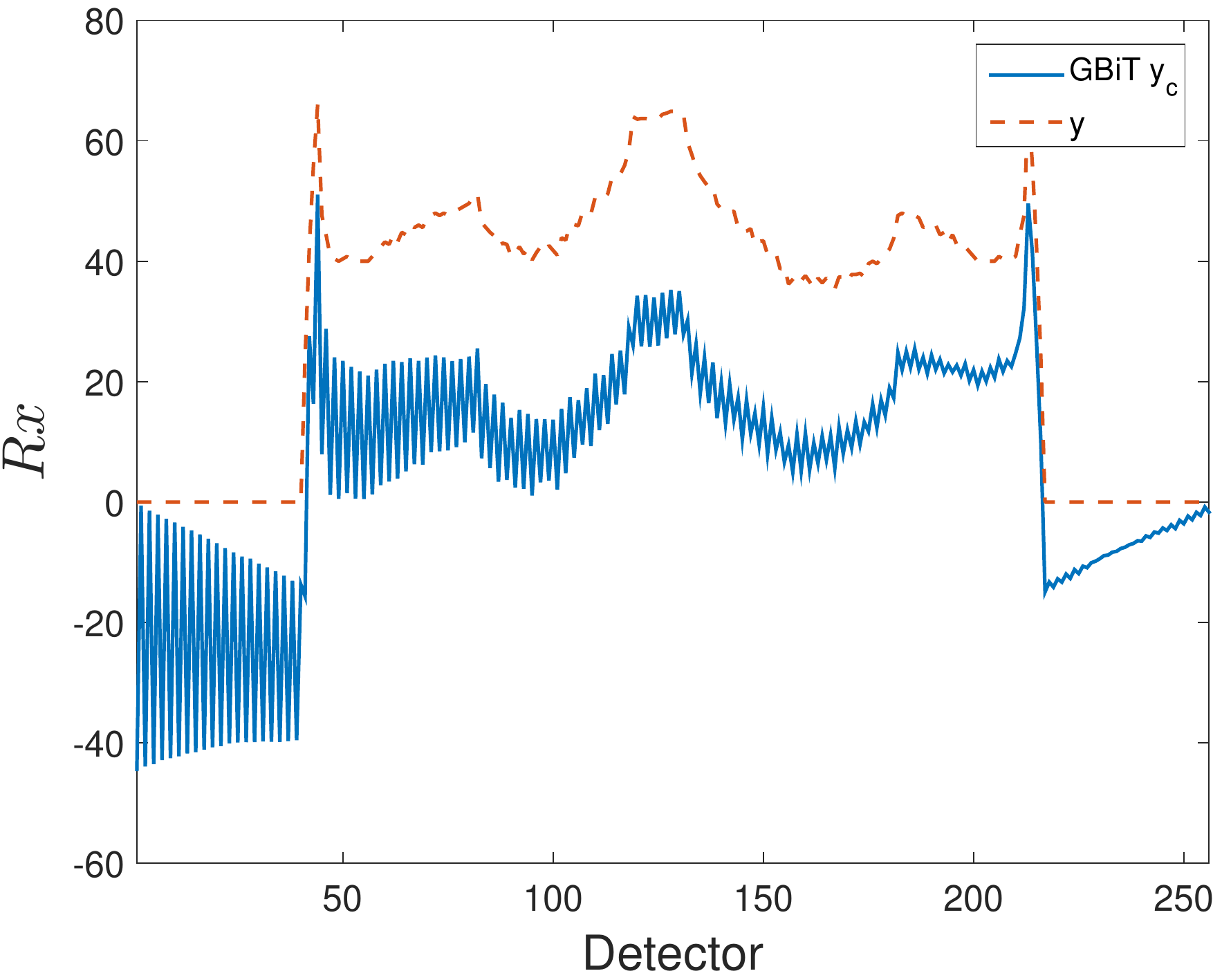}\label{fig:testp3:d}}
\caption{Reconstructions identical to \hypref{Fig.}{fig:testp2}, but with a constant offset added to the noise, i.e. $e$ is replaced by $e + 5$ in \eqref{eq:addnoise}.}
\label{fig:testp3}
\end{figure}

%%%%%%%%%%%%%%%%%%%%%%%% SECTION: Numerical Applications %%%%%%%%%%%%%%%%%%%%%%%%
\section{Application to phase contrast tomography}\label{sec:numapp}
The previous section looked at a single projection angle and only at the effect of the matrix $D_f$ and $D_c$. We will now look at the full tomography problem for multiple angles in $[0, \pi]$. We start with a test problem and end with an application to real phase contrast data.

% ------------ SUBSECTION  TESTPROBLEM  ---------------------
\subsection{PC-CT with simulated data}
We continue using the $256\times 256$ Shepp-logan phantom from \hypref{section}{sec:testp} as a test image and measurement data is generated from $360$ projection angles with $10\%$ model error using the same weighting method as in \eqref{eq:modelerror} and add $10\%$ noise to the resulting measurements as in \eqref{eq:addnoise}. Because the inverse of $D_f$ is given by a back substitution, we will also consider the reconstruction given by
\begin{equation}\label{eq:phr}
Rx=D_f^{-1}b,
\end{equation}
with and without regularization. This can be seen as two steps. First, reconstructing the Radon projections of $\delta$ for each projection angle from the measured phase contrast data and, second, approximate the image $x$ based on the models for classic absorption based CT. This reconstruction is shown in \hypref{Fig.}{fig:pcct:a}. However, since the classic tomography problem is know to be very ill-conditioned, this reconstruction is of low quality. Regularization, see \hypref{Fig.}{fig:pcct:d}, suppresses this but over-smooths the reconstruction compared to directly solving $DRx=b$. This is demonstrated by the value of the regularization parameter $\lambda$, shown in \hypref{Fig.}{fig:pcct2:c}, which is significantly higher compared to the other values. 

\hypref{Fig.}{fig:pcct:b} and \hypref{Fig.}{fig:pcct:c} demonstrate the results from the previous section, namely that the forward model is much less sensitive to the introduced measurement and model errors. This can be observed in \hypref{Fig.}{fig:pcct2:b}, which shows the relative error off all the reconstructions. Here it can be seen that the LSQR solution using the central model shows signs of semi-convergence after 13 iterations, whereas for the solution using the forward model this only happens after 44 iterations and at a slower rate. When we compute the reconstruction using regularization, \hypref{Fig.}{fig:pcct:e} and \hypref{Fig.}{fig:pcct:f}, there is no semi-convergence is, but the reconstruction with the central model has less contrast due to the fact that it is more oscillatory than the other regularized reconstructions. This can be seen in the horizontal slices of the images, see \hypref{Fig.}{fig:pcct:g}, \hypref{Fig.}{fig:pcct:h} and \hypref{Fig.}{fig:pcct:i}.

\begin{figure}
\centering
\subfloat[]{\includegraphics[width=0.28\linewidth]{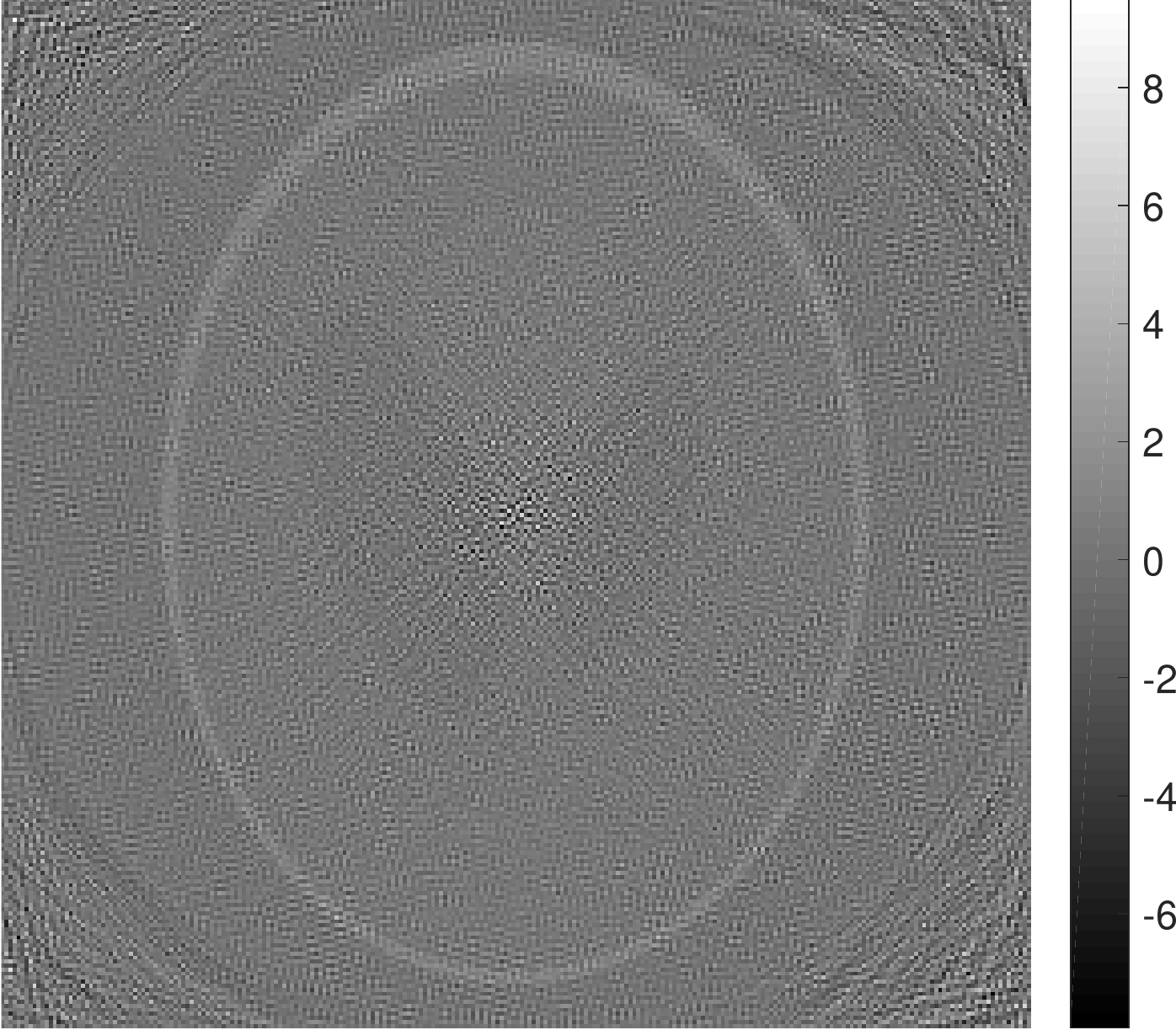}\label{fig:pcct:a}}\hspace*{\fill}
\subfloat[]{\includegraphics[width=0.28\linewidth]{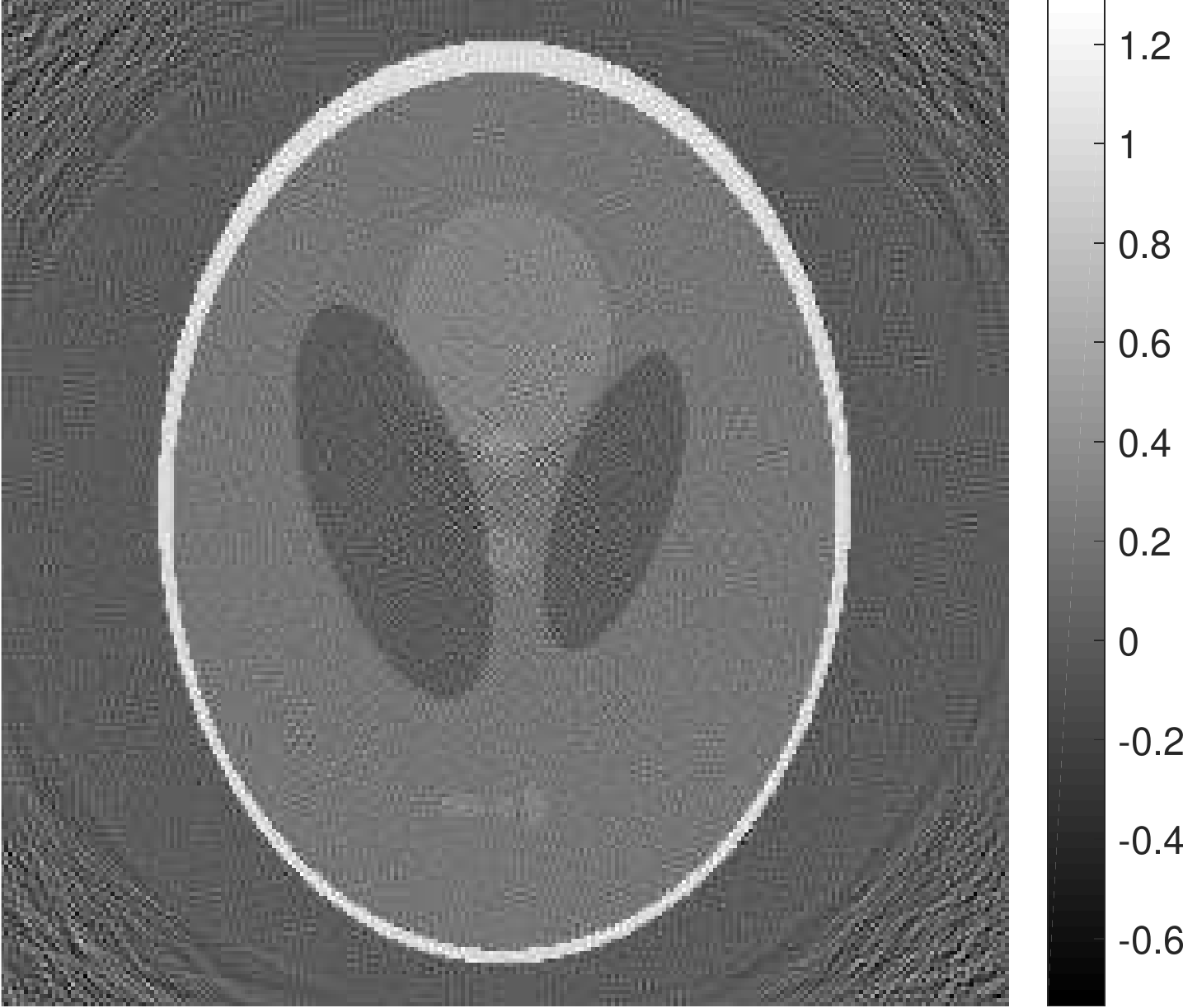}\label{fig:pcct:b}}\hspace*{\fill}
\subfloat[]{\includegraphics[width=0.28\linewidth]{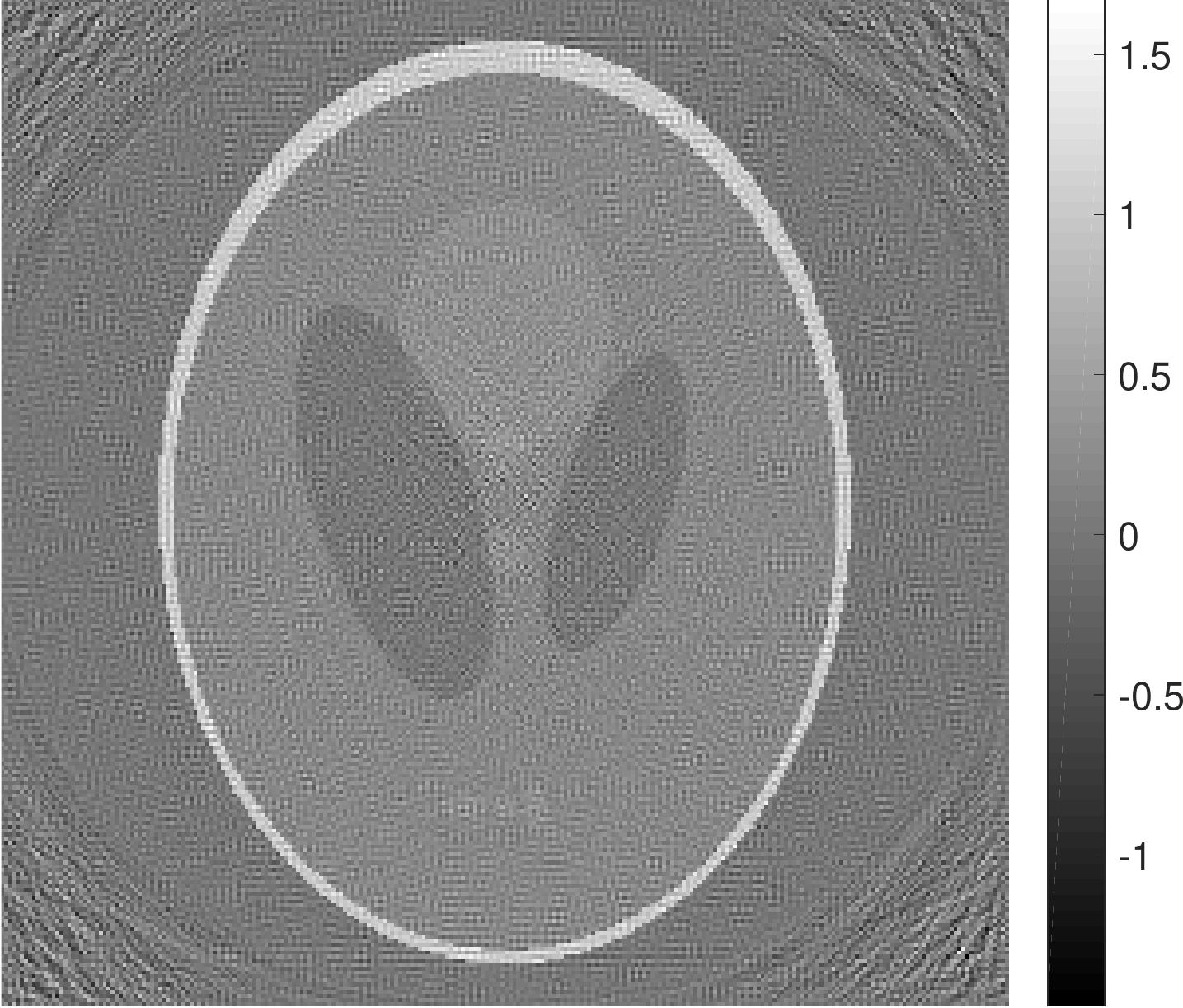}\label{fig:pcct:c}}\\
\subfloat[]{\includegraphics[width=0.28\linewidth]{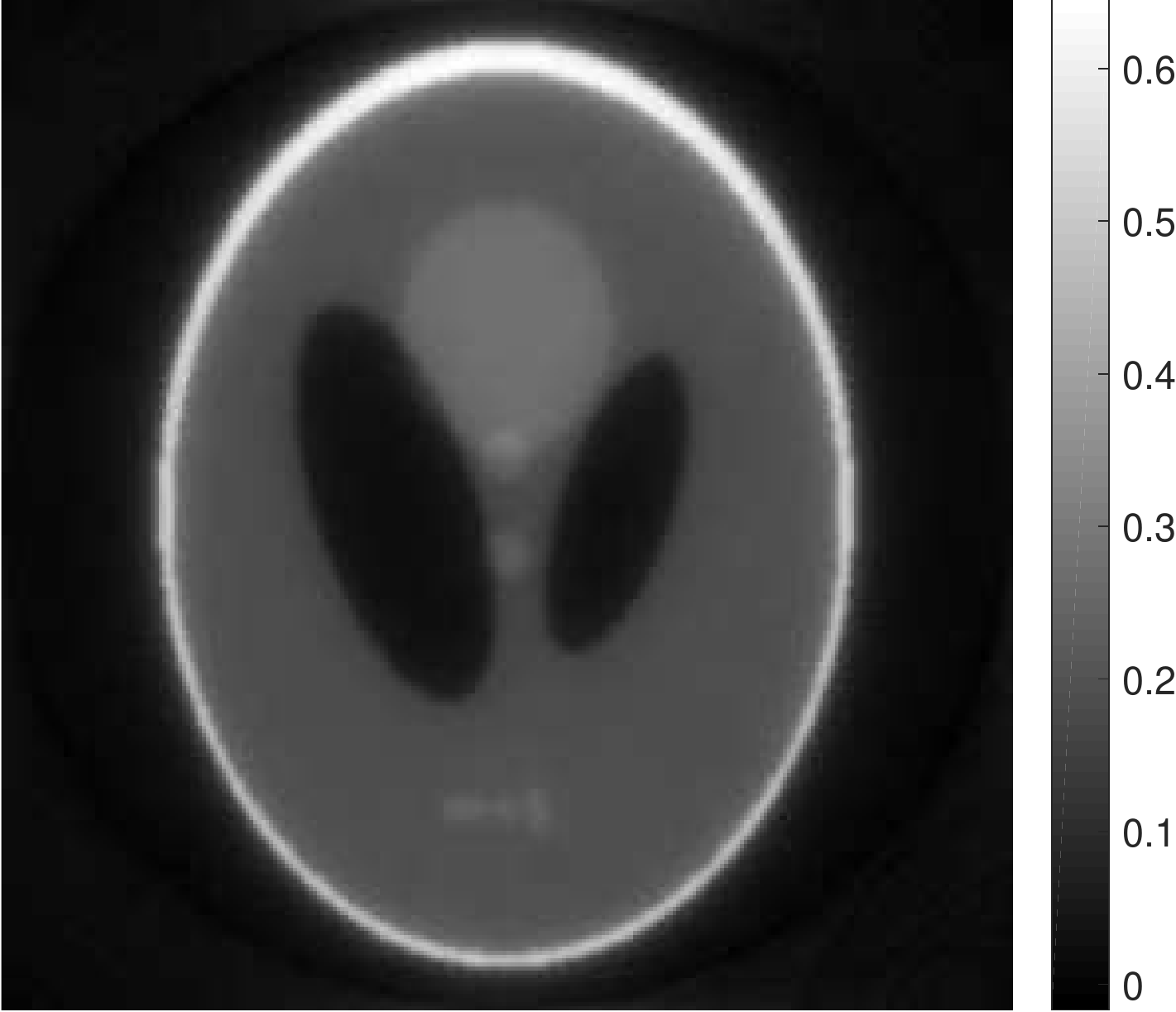}\label{fig:pcct:d}}\hspace*{\fill}
\subfloat[]{\includegraphics[width=0.28\linewidth]{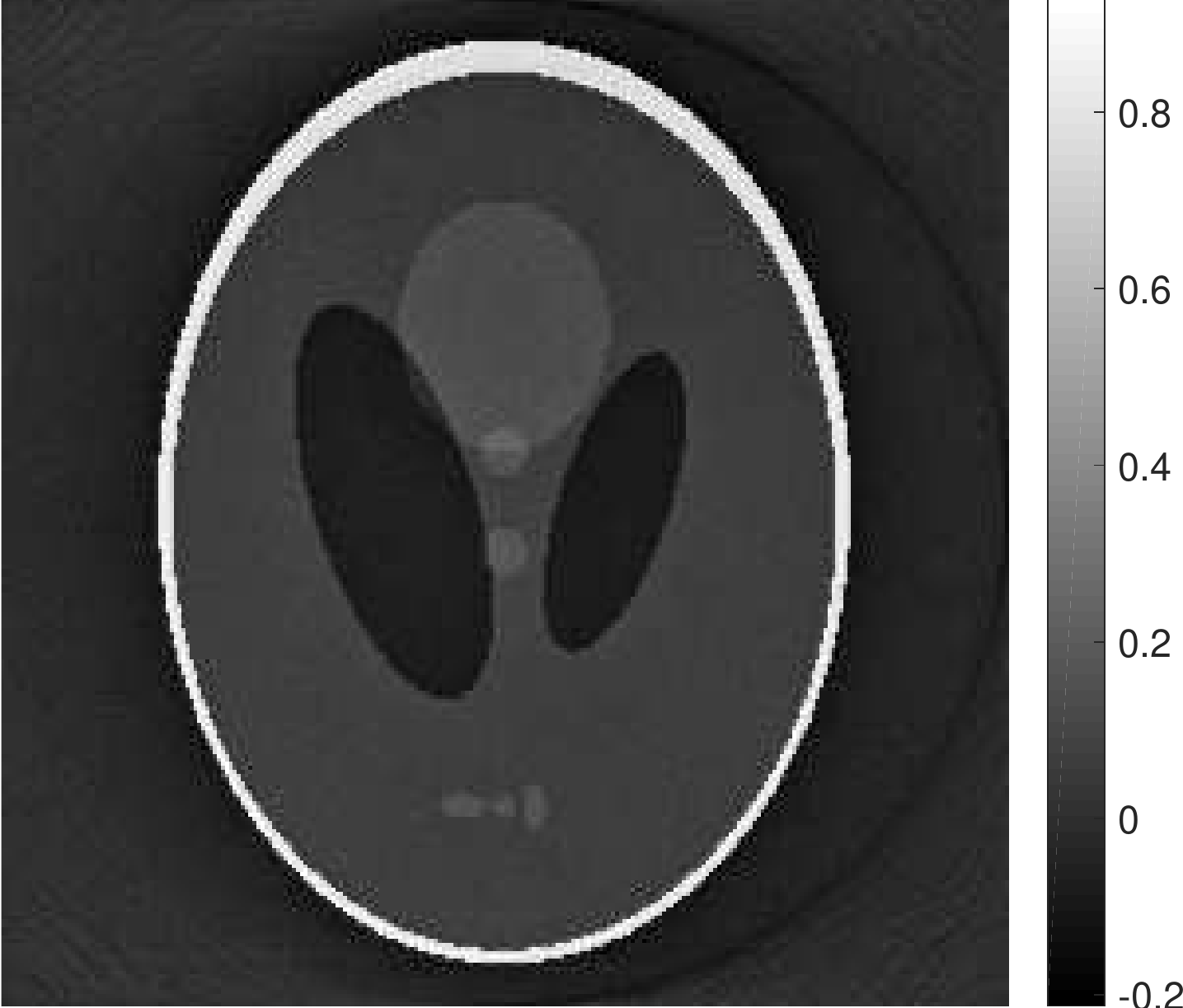}\label{fig:pcct:e}}\hspace*{\fill}
\subfloat[]{\includegraphics[width=0.28\linewidth]{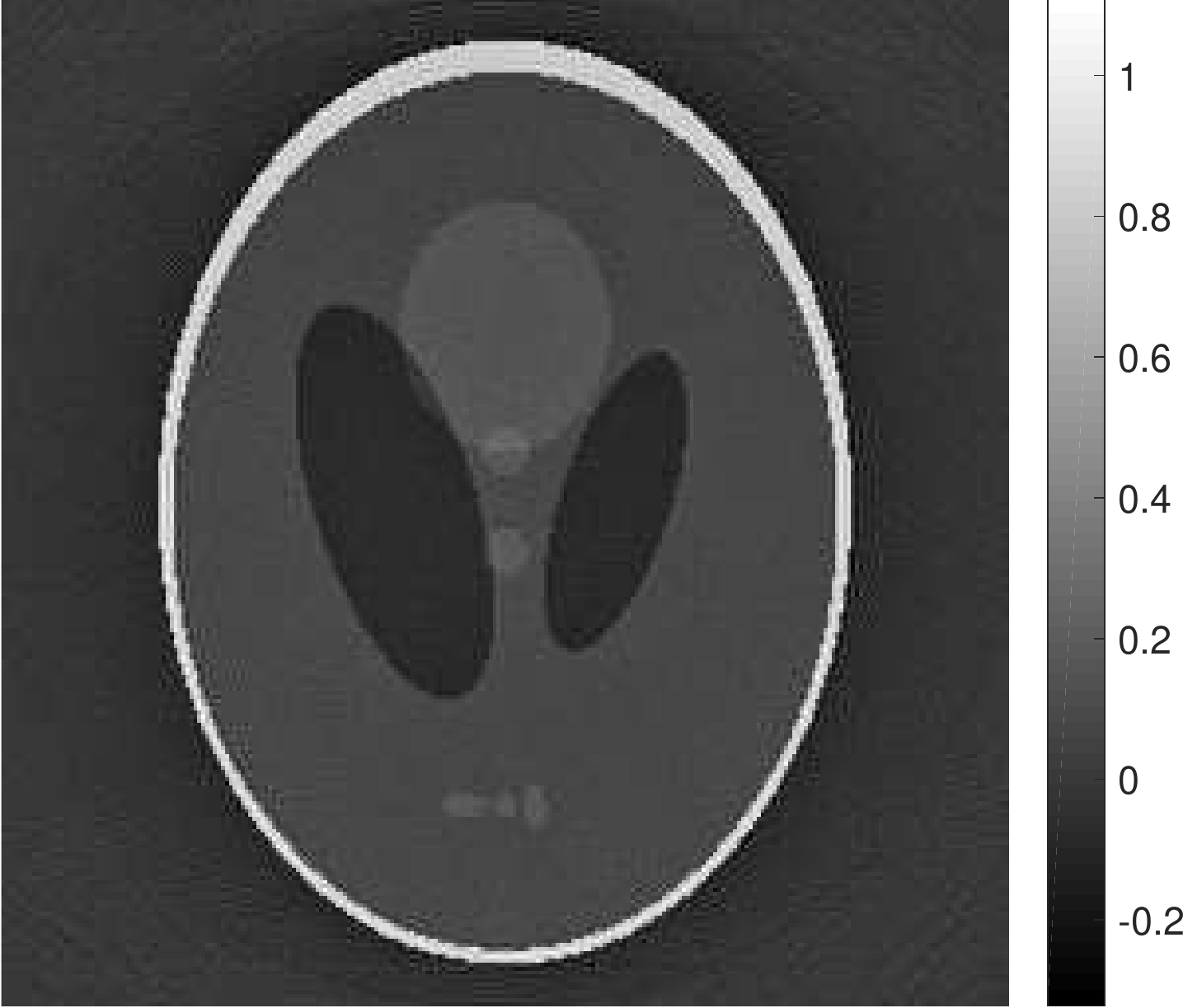}\label{fig:pcct:f}}\\
\subfloat[]{\includegraphics[width=0.28\linewidth]{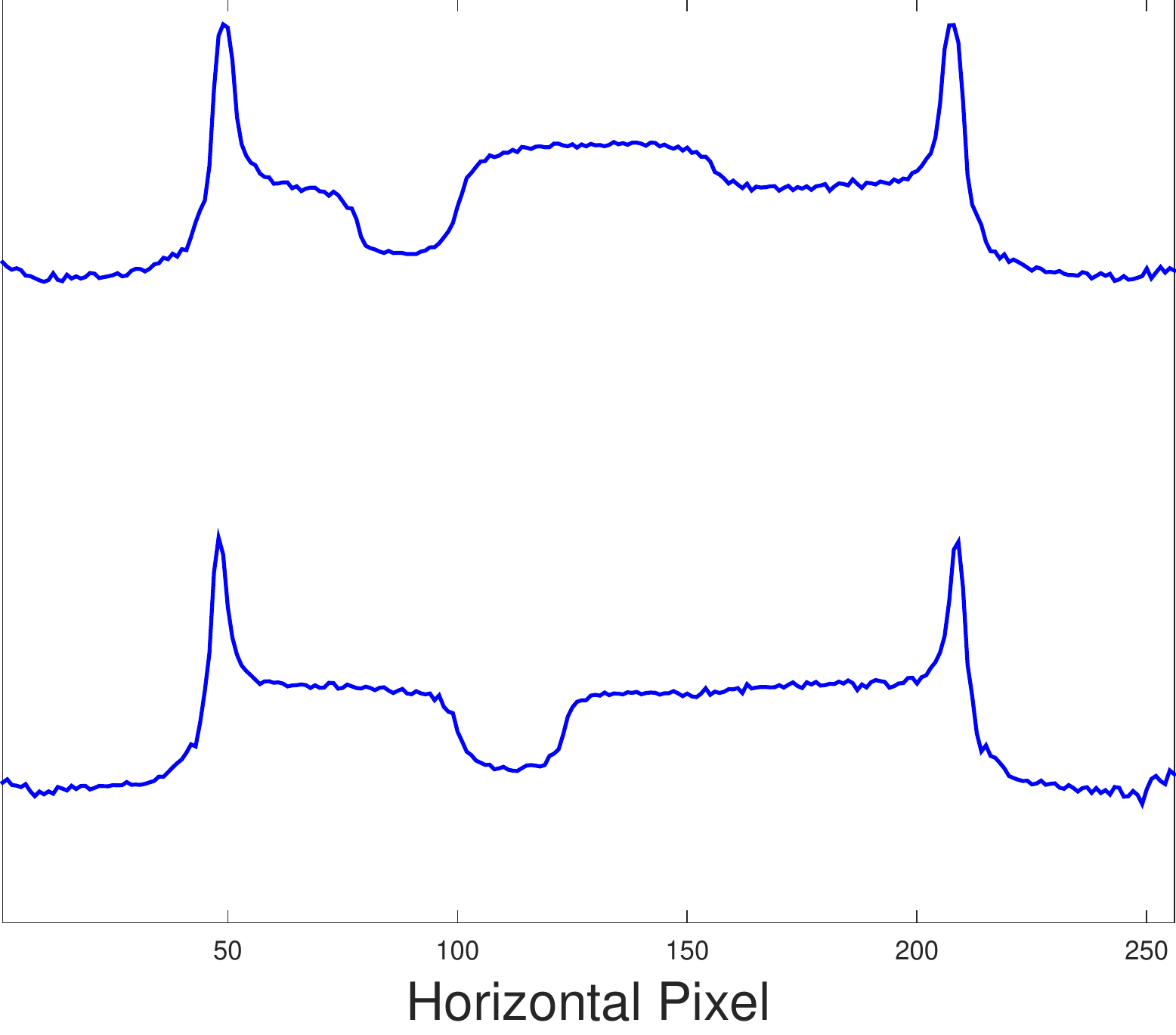}\label{fig:pcct:g}}\hspace*{\fill}
\subfloat[]{\includegraphics[width=0.28\linewidth]{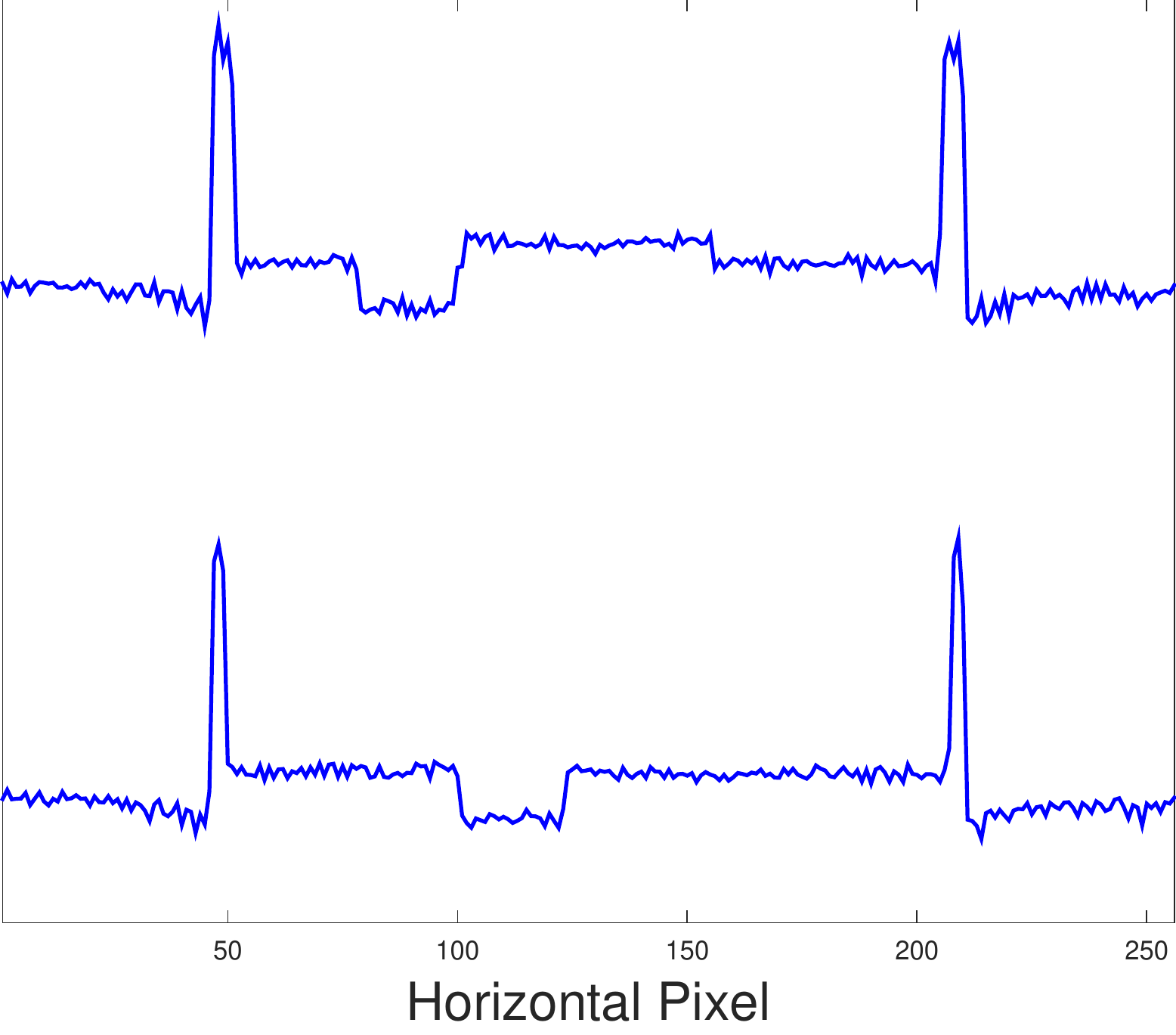}\label{fig:pcct:h}}\hspace*{\fill}
\subfloat[]{\includegraphics[width=0.28\linewidth]{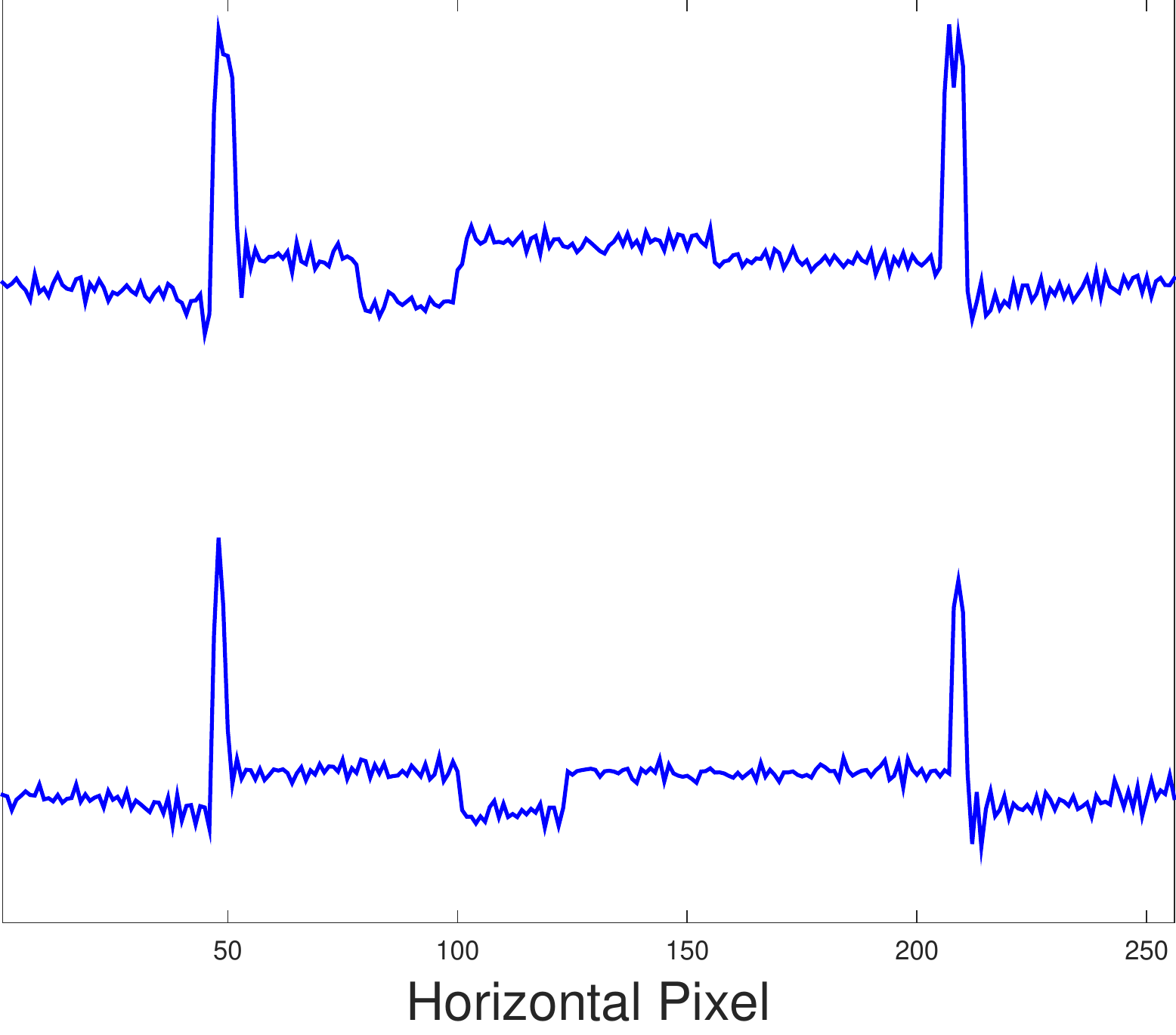}\label{fig:pcct:i}}
\caption{\protect\subref{fig:pcct:a}, \protect\subref{fig:pcct:d} Reconstruction based on the phase retrieval approach from \eqref{eq:phr} without and with regularization. \protect\subref{fig:pcct:b}, \protect\subref{fig:pcct:e} Reconstruction of $D_fRx=b$ without and with regularization. \protect\subref{fig:pcct:c}, \protect\subref{fig:pcct:f} Reconstruction of $D_cRx=b$ without and with regularization. Note the different scales of the reconstructions. \protect\subref{fig:pcct:g},\protect\subref{fig:pcct:h} \& \protect\subref{fig:pcct:i} Normalized horizontal slices of the regularized reconstructions above them at vertical pixels 85 and 171 (top to bottom).}
\label{fig:pcct}
\end{figure}

\begin{figure}
\centering
\subfloat[]{\includegraphics[width=0.45\linewidth]{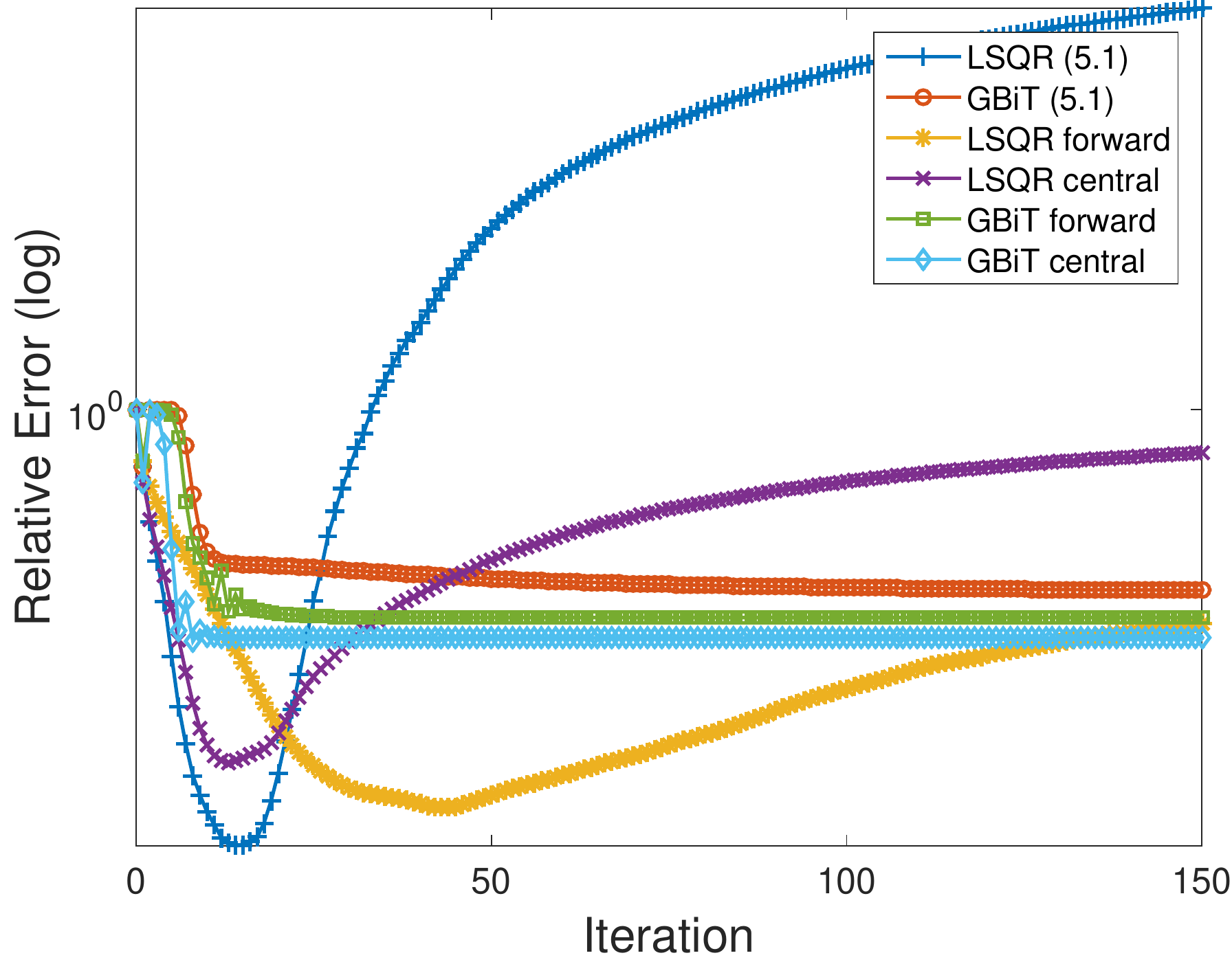}\label{fig:pcct2:b}}\quad
\subfloat[]{\includegraphics[width=0.45\linewidth]{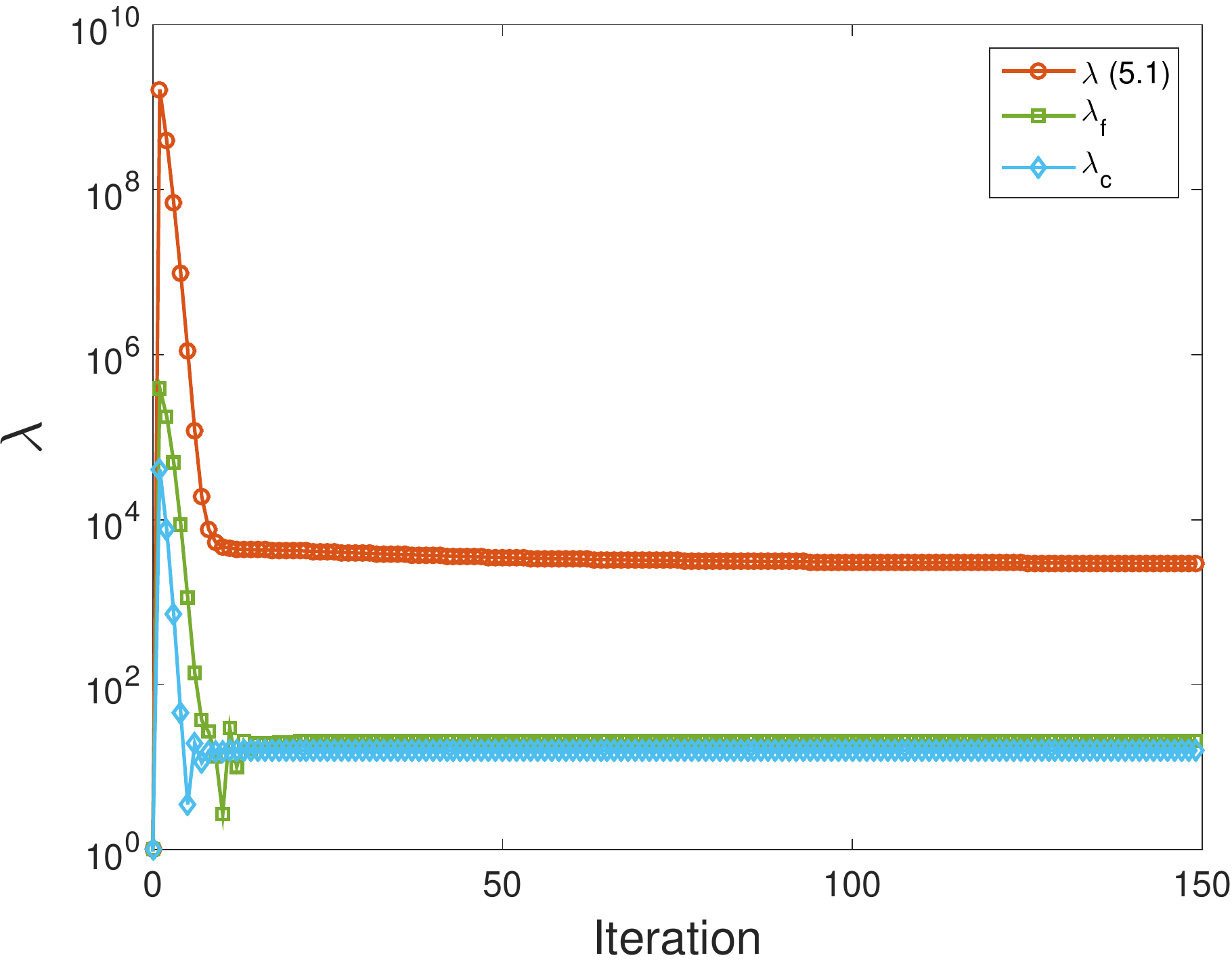}\label{fig:pcct2:c}}
\caption{Relative error \protect\subref{fig:pcct2:b} and regularization parameter \protect\subref{fig:pcct2:c} for the reconstructions from \hypref{Fig.}{fig:pcct}.}
\label{fig:pcct2}
\end{figure}

% ------------ SUBSECTION  DPC PROBLEM  ---------------------
\subsection{PC-CT with experimental data}
Finally, to illustrate the validity of our test problem and to illustrate that \hyperref[alg:GBiT]{GBiT} is a viable algorithm for real data sets, we apply the same calculations on a PC-CT data set from experimental measurements provided by Franz Pfeiffer \cite{Pfeiffer2007_2}. The data was taken over 2401 angles in $[0, \pi[$, with projections onto 1637 detectors. Because we have no prior knowledge about the noise, we use the alternative update scheme described in \hypref{section}{sec:altup}. As can be seen from \hypref{Fig.}{fig:real}, \hypref{Fig.}{fig:real2} and \hypref{Fig.}{fig:zoom}, the results for the real data set are similar to our test problem. Here, all reconstructions were made by calculating 20 iterations of the algorithms, but for the \hyperref[alg:GBiT]{GBiT} algorithm applied to the forward model, the central model and \eqref{eq:phr}, the discrepancy principle was satisfied after respectively 8, 9 and 7 iterations. 

A final note should be added here about the runtime of the algorithm. The ASTRA toolbox can be used as a black box for the multiplication with $R$ and $R^T$ in the algorithms. This can be done on either the CPU or GPU. These reconstructions were done using the CPU code, with an average runtime of 3 minutes per iterations for the \hyperref[alg:GBiT]{GBiT} method. Comparing this to the runtime filtered back projection (20 seconds) this is a very poor result. However, when using the GPU the runtime of our iterative reconstruction drops significantly and the full algebraic reconstruction can be done in approximately 60 seconds. The reason we did not do this is because the `line weights' we described in \hypref{section}{sec:disrec} aren't (yet) available for the GPU code. On the GPU only a `linear weighting' scheme has been implemented that takes into account the neighbouring pixels. Furthermore, the GPU implementation of this scheme uses some approximations in the multiplication with $R^T$, which can make the reconstruction differ from the CPU reconstruction with the same weights.

\begin{figure}
\centering
\subfloat[LSQR, \eqref{eq:phr}]{\includegraphics[width=0.4\linewidth]{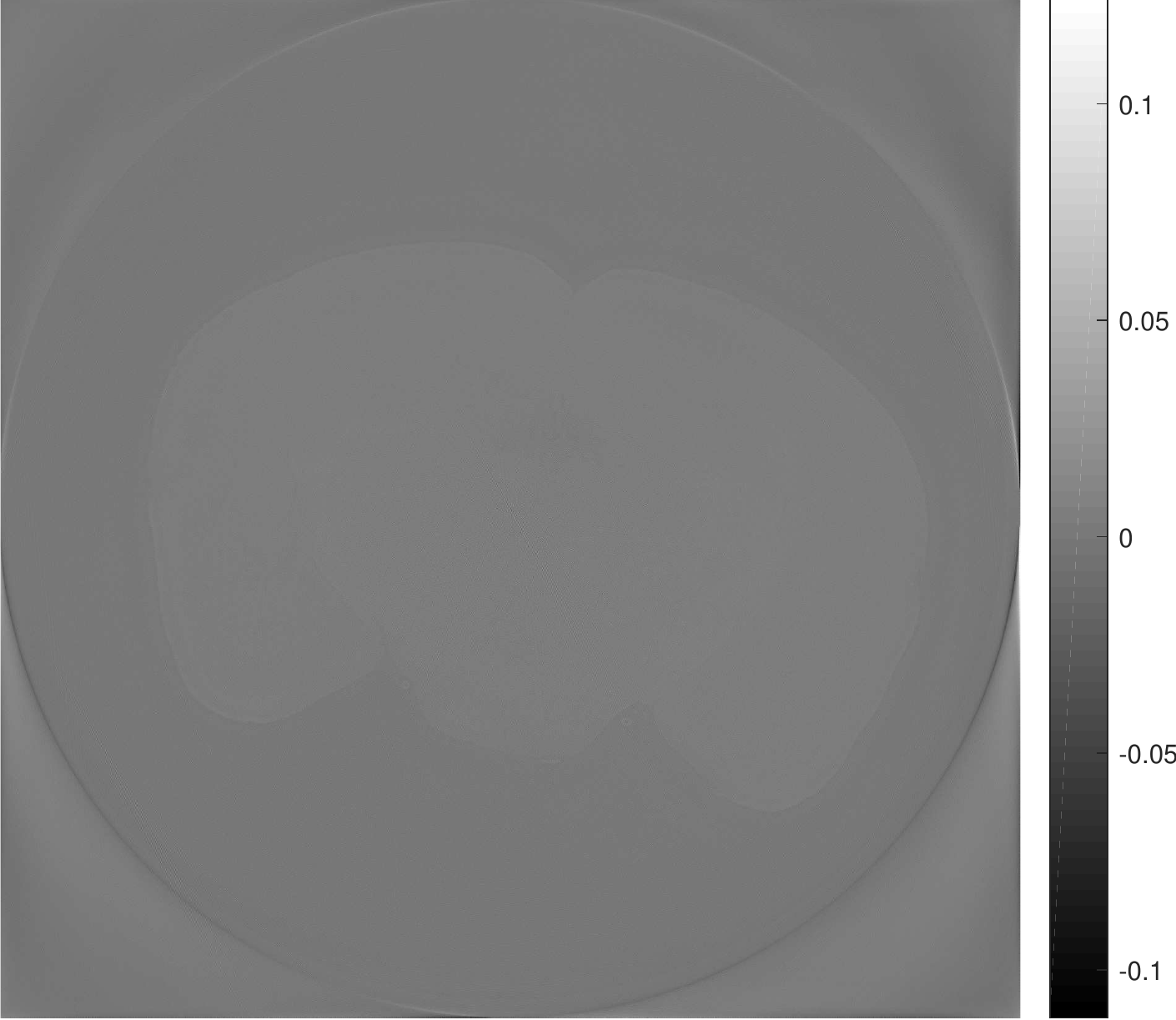}\label{fig:real:a}}\qquad\qquad
\subfloat[BGiT, \eqref{eq:phr}]{\includegraphics[width=0.4\linewidth]{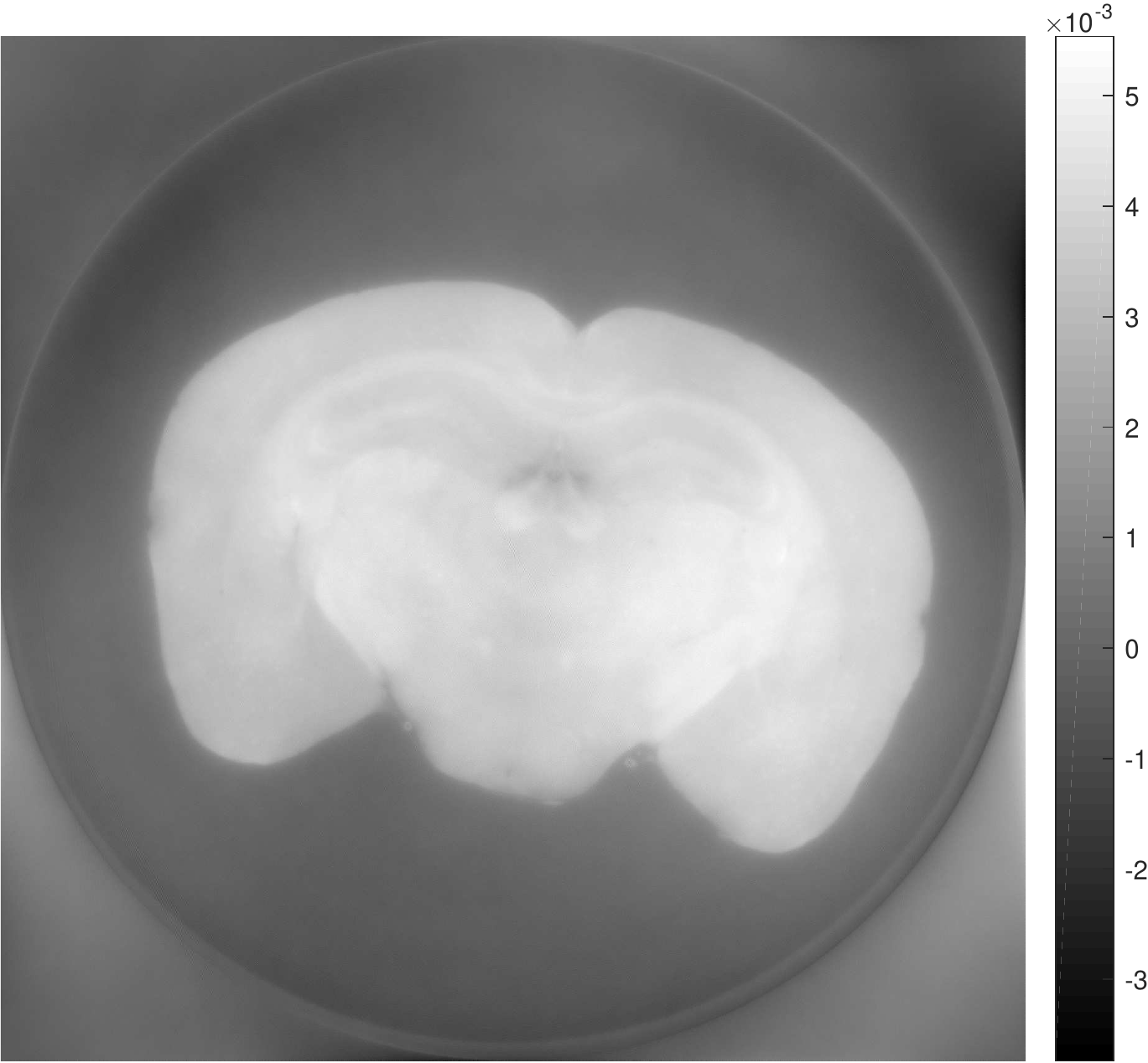}\label{fig:real:d}}\\
\subfloat[LSQR, forward model]{\includegraphics[width=0.4\linewidth]{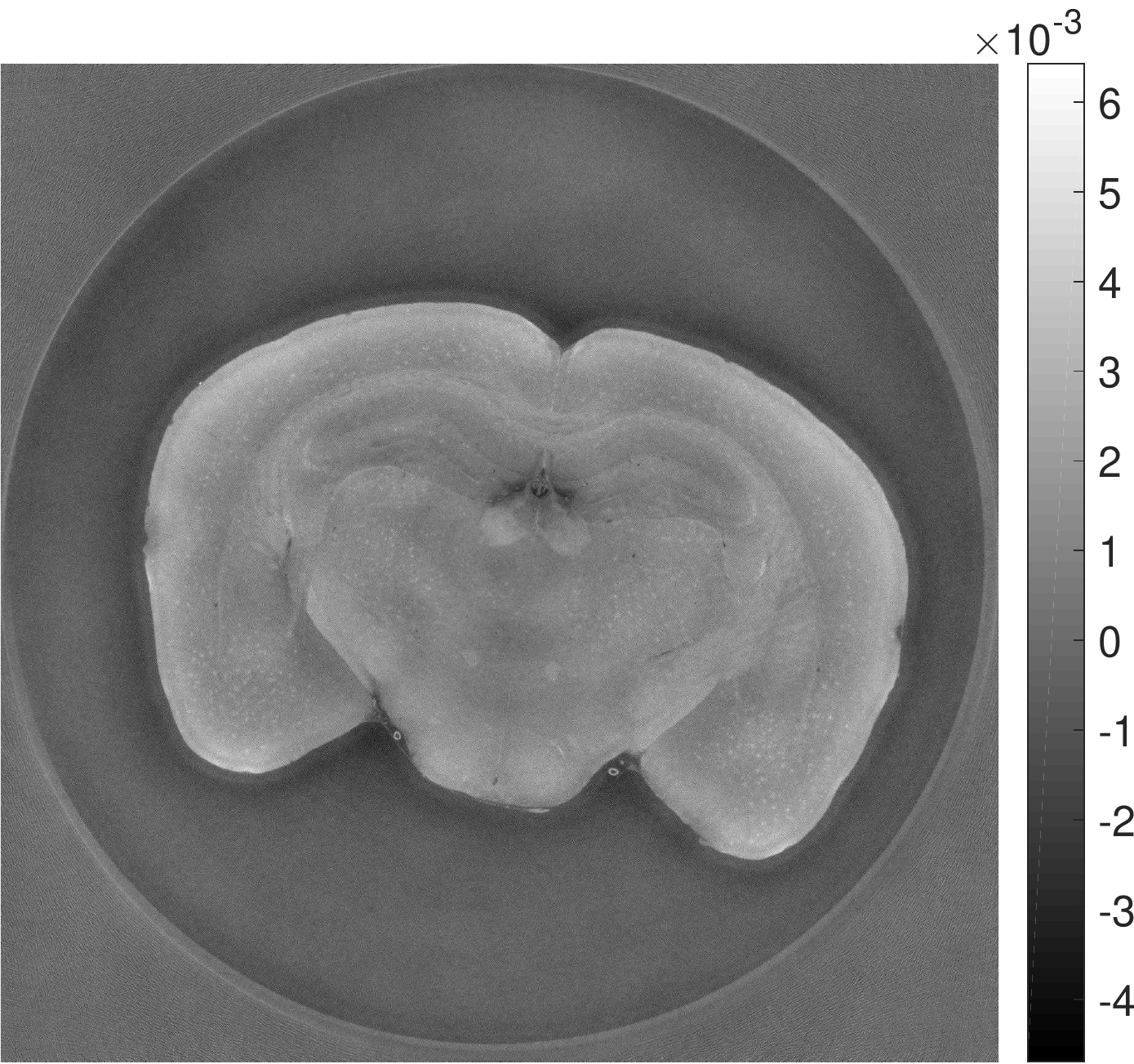}\label{fig:real:b}}\qquad\qquad
\subfloat[GBiT, forward model]{\includegraphics[width=0.4\linewidth]{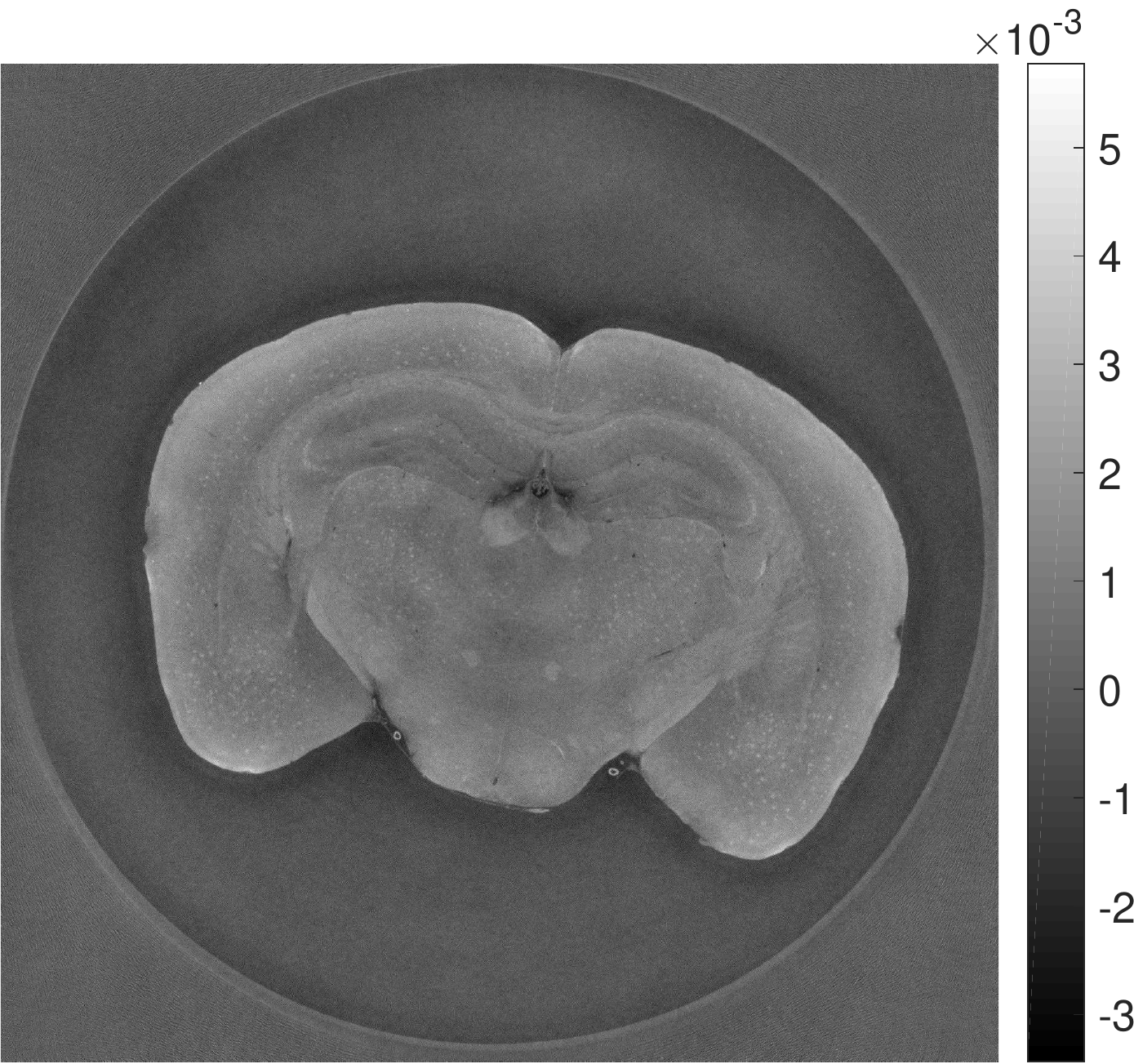}\label{fig:real:e}}\\
\subfloat[LSQR, central model]{\includegraphics[width=0.4\linewidth]{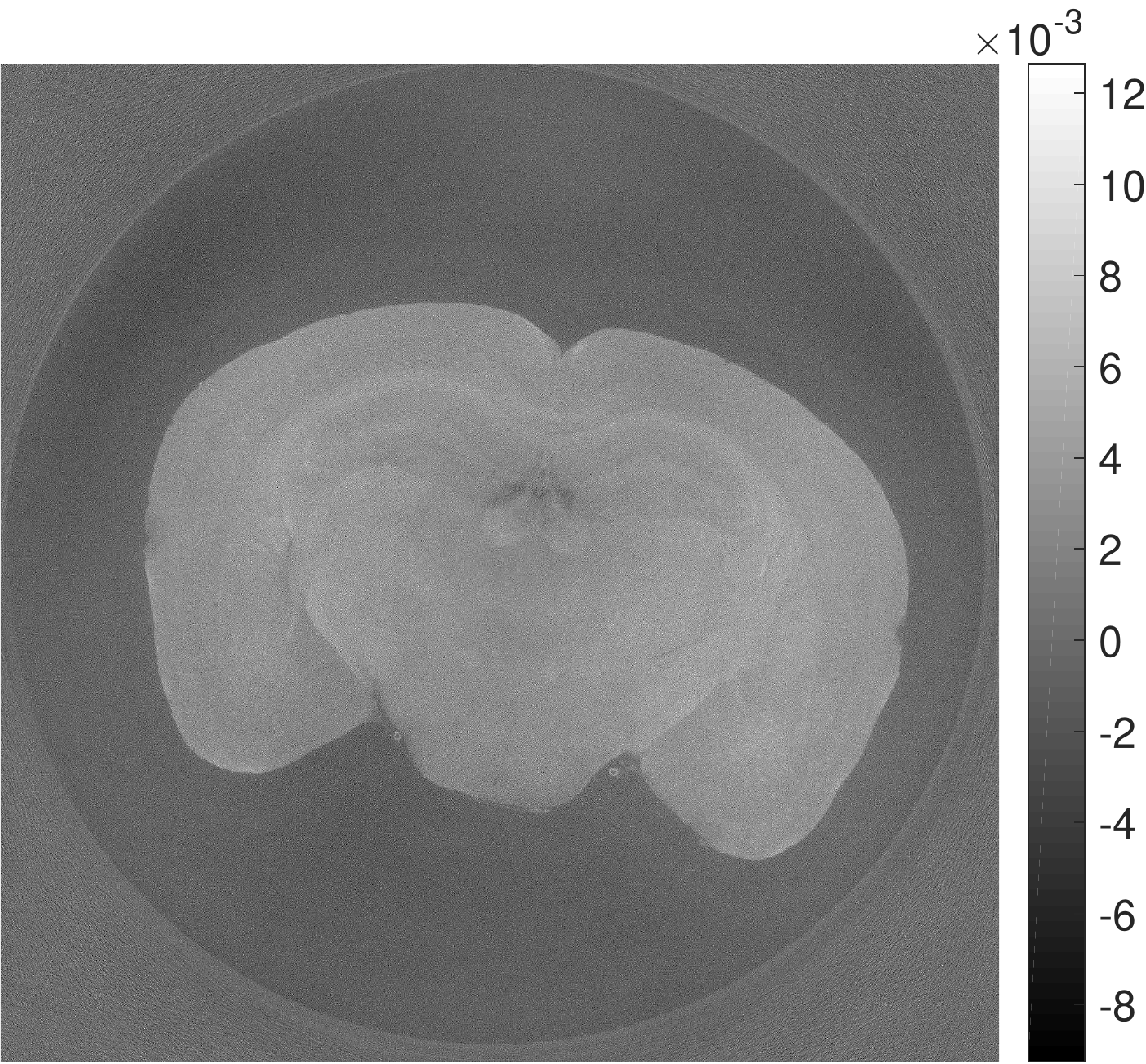}\label{fig:real:c}}\qquad\qquad
\subfloat[GBiT, central model]{\includegraphics[width=0.4\linewidth]{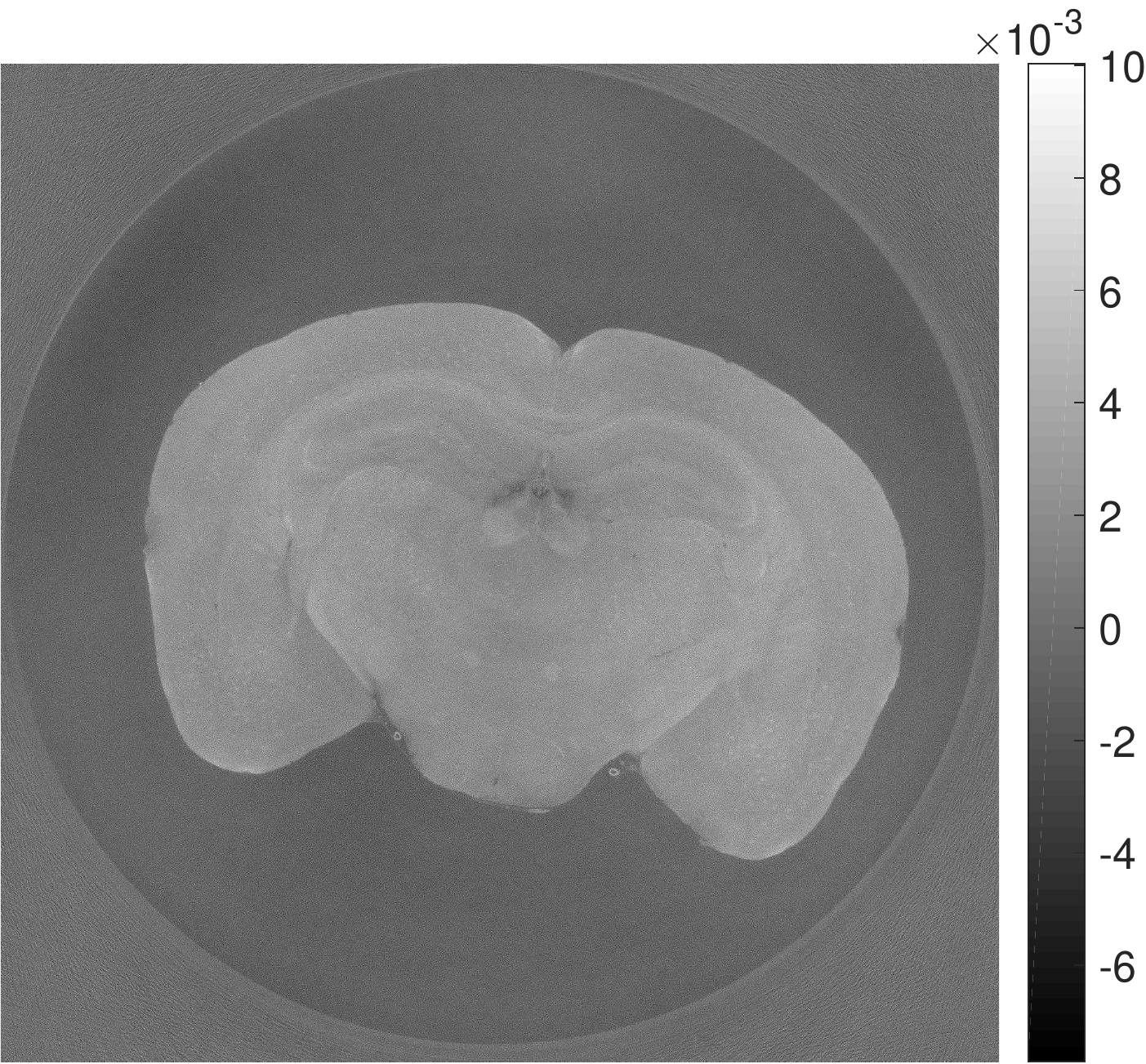}\label{fig:real:f}}
\caption{\protect\subref{fig:real:a}\protect\subref{fig:real:d} Reconstruction based on the phase retrieval approach from \eqref{eq:phr} without and with regularization. \protect\subref{fig:real:b}, \protect\subref{fig:real:e} Reconstruction of $D_fRx=b_n$ without and with regularization. \protect\subref{fig:real:c}, \protect\subref{fig:real:f} Reconstruction of $D_cRx=b_n$ without and with regularization.}
\label{fig:real}
\end{figure}

\begin{figure}
\centering
\subfloat[]{\includegraphics[width=0.45\linewidth]{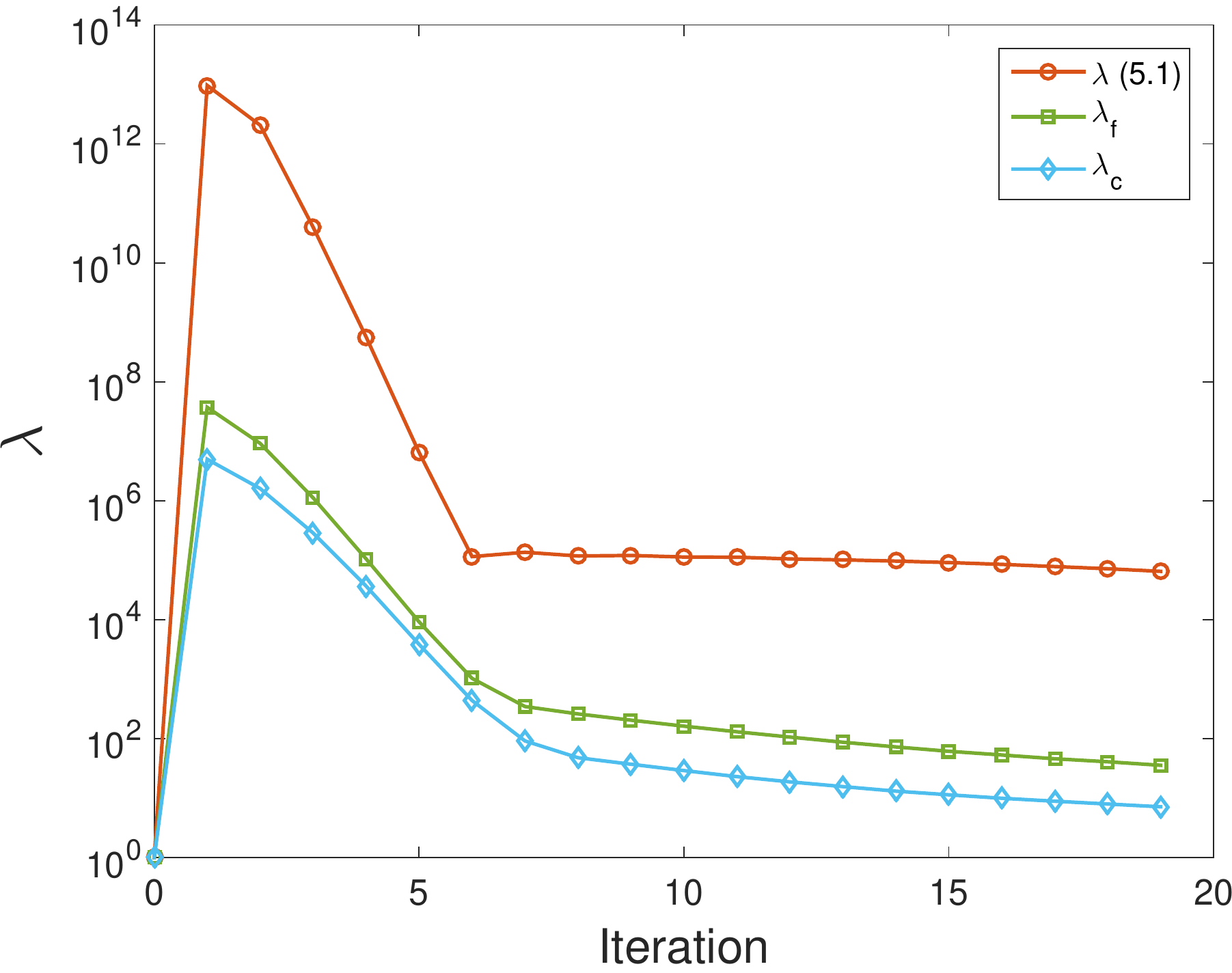}\label{fig:real2:a}}\\
\subfloat[GBiT, forward model]{\includegraphics[width=0.4\linewidth]{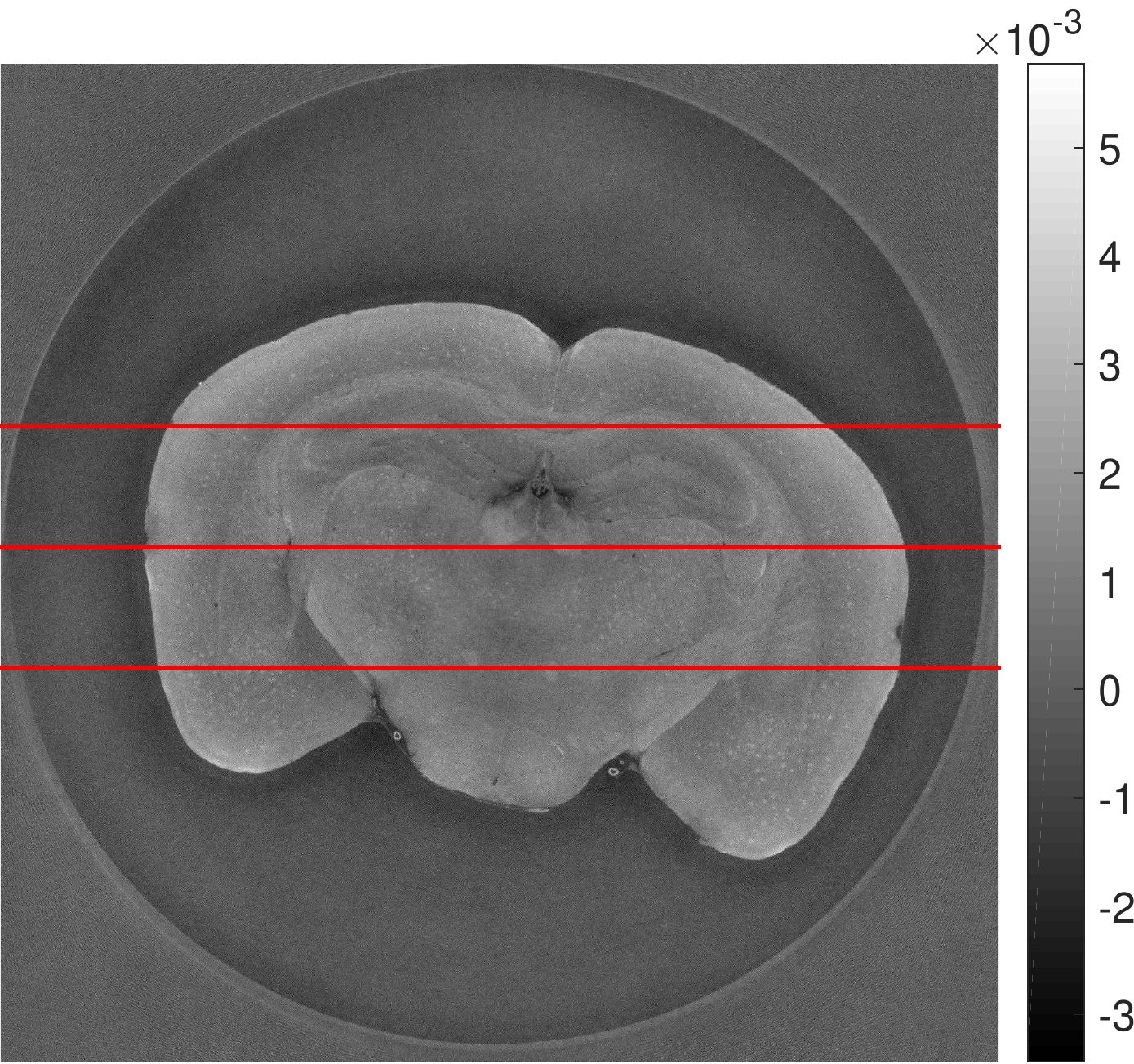}\label{fig:real2:b}}\qquad\qquad
\subfloat[GBiT, central model]{\includegraphics[width=0.4\linewidth]{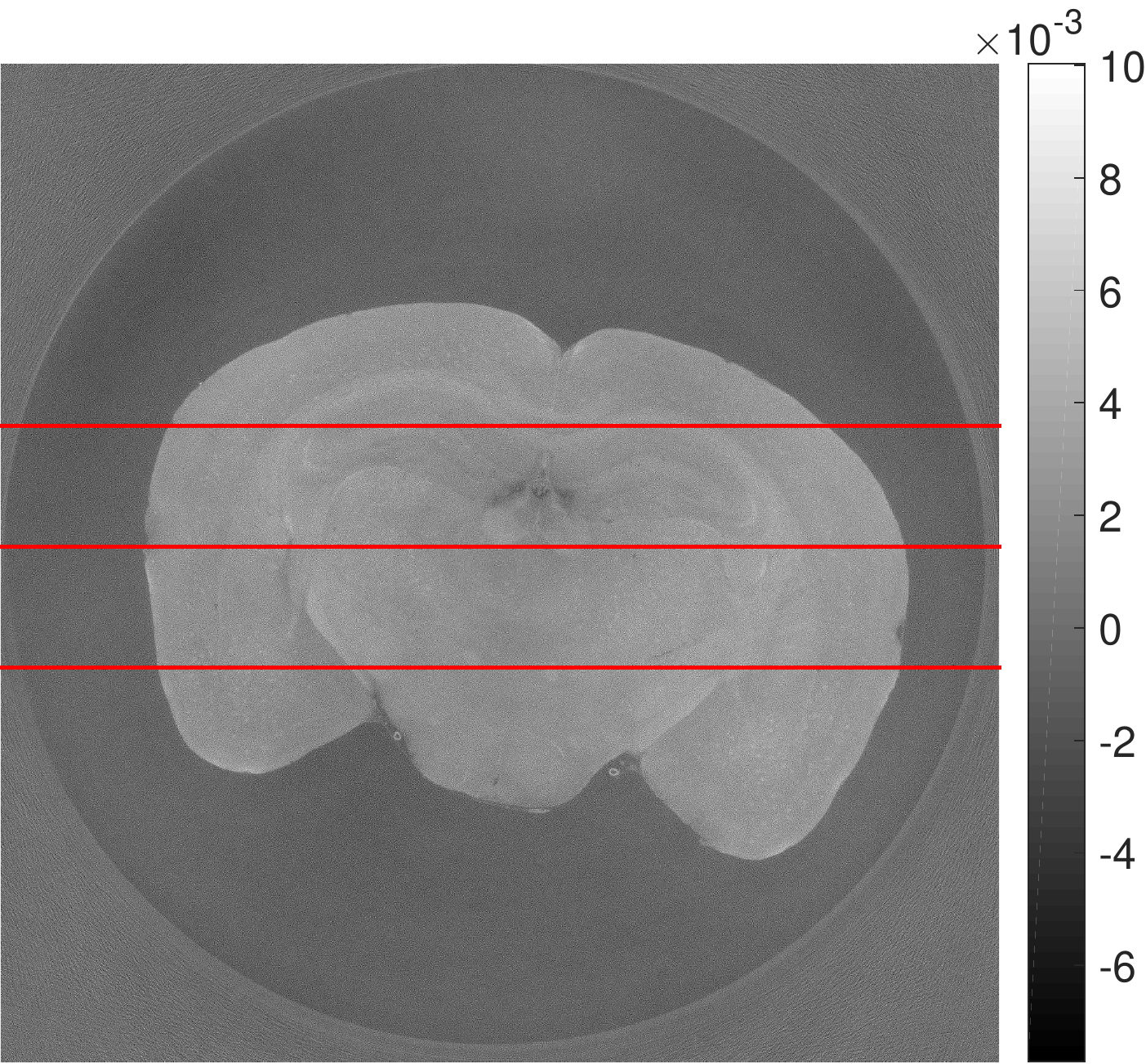}\label{fig:real2:c}}\\
\subfloat[]{\includegraphics[width=0.4\linewidth]{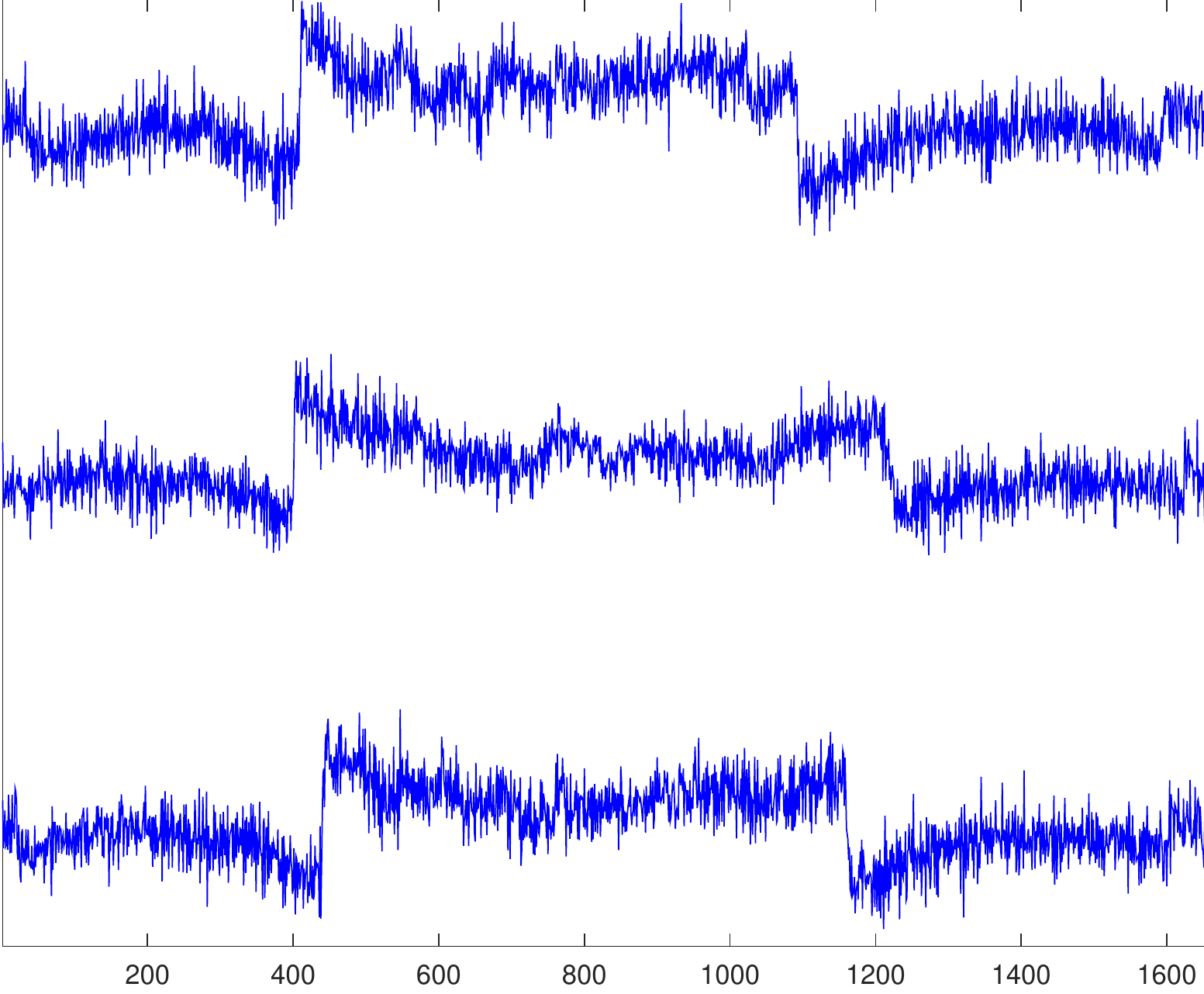}\label{fig:real2:d}}\qquad\qquad
\subfloat[]{\includegraphics[width=0.4\linewidth]{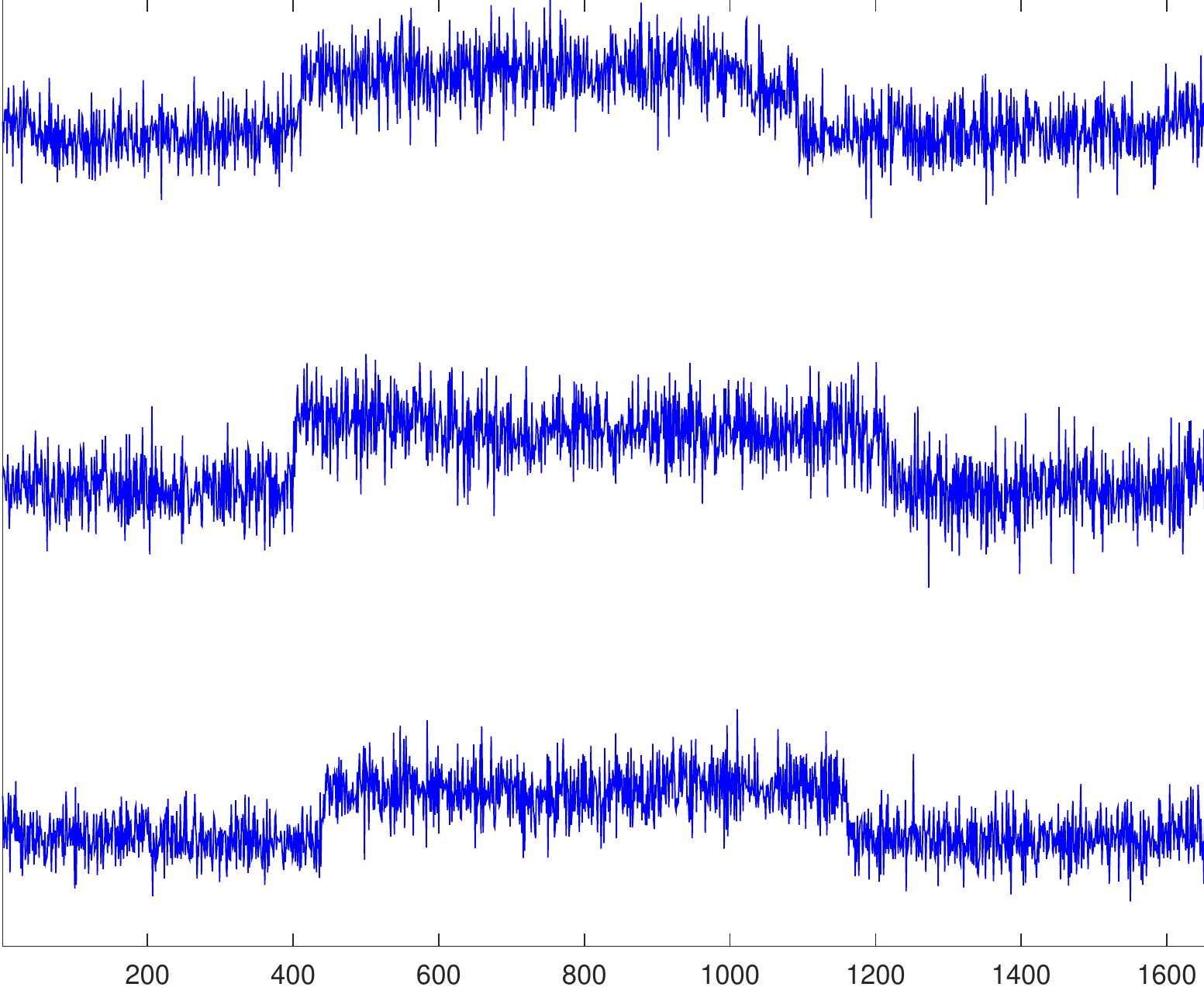}\label{fig:real2:e}}\\
\caption{\protect\subref{fig:real2:a} Regularization parameter $\lambda$. \protect\subref{fig:real2:b}-\protect\subref{fig:real2:e} Normalized horizontal slices of the regularized reconstructions above them at the denoted vertical positions. The higher oscillations for the the reconstruction based on the central model result in a lower contrast.}
\label{fig:real2}
\end{figure}

\begin{figure}
\centering
\subfloat[LSQR, forward model]{\includegraphics[width=0.4\linewidth]{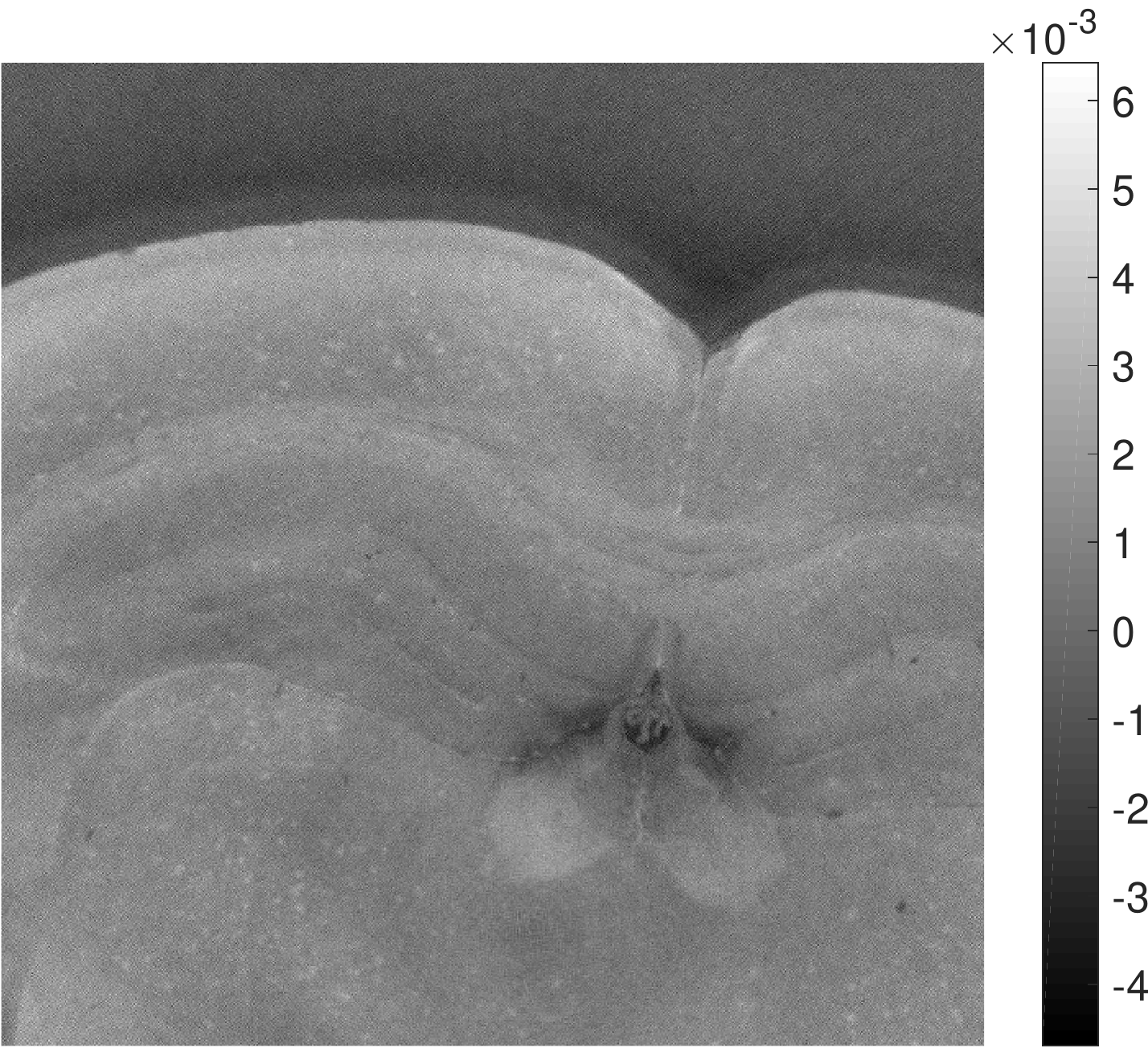}\label{fig:zoom:a}}\qquad\qquad
\subfloat[GBiT, forward model]{\includegraphics[width=0.4\linewidth]{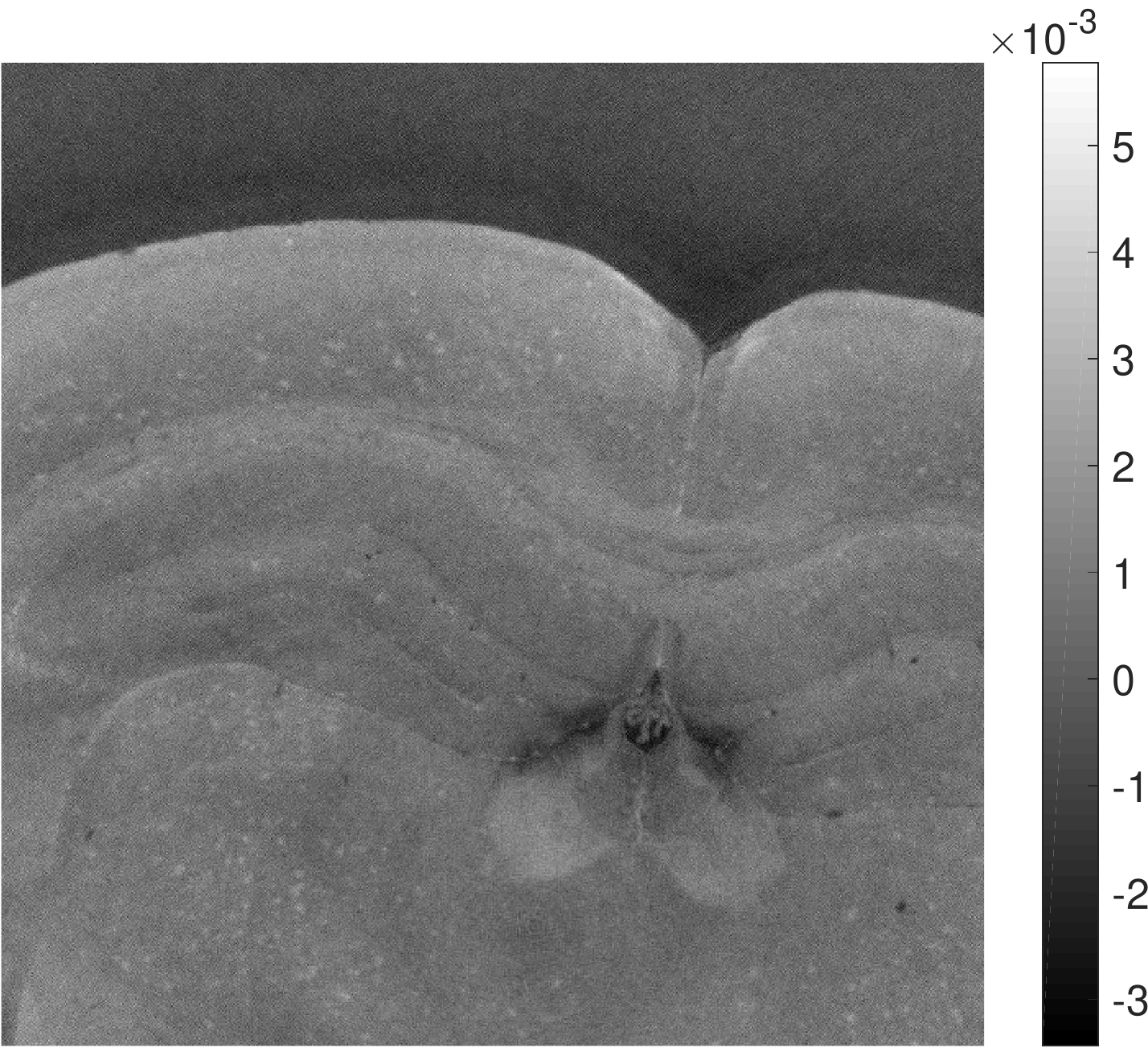}\label{fig:zoom:b}}\\
\subfloat[LSQR, central model]{\includegraphics[width=0.4\linewidth]{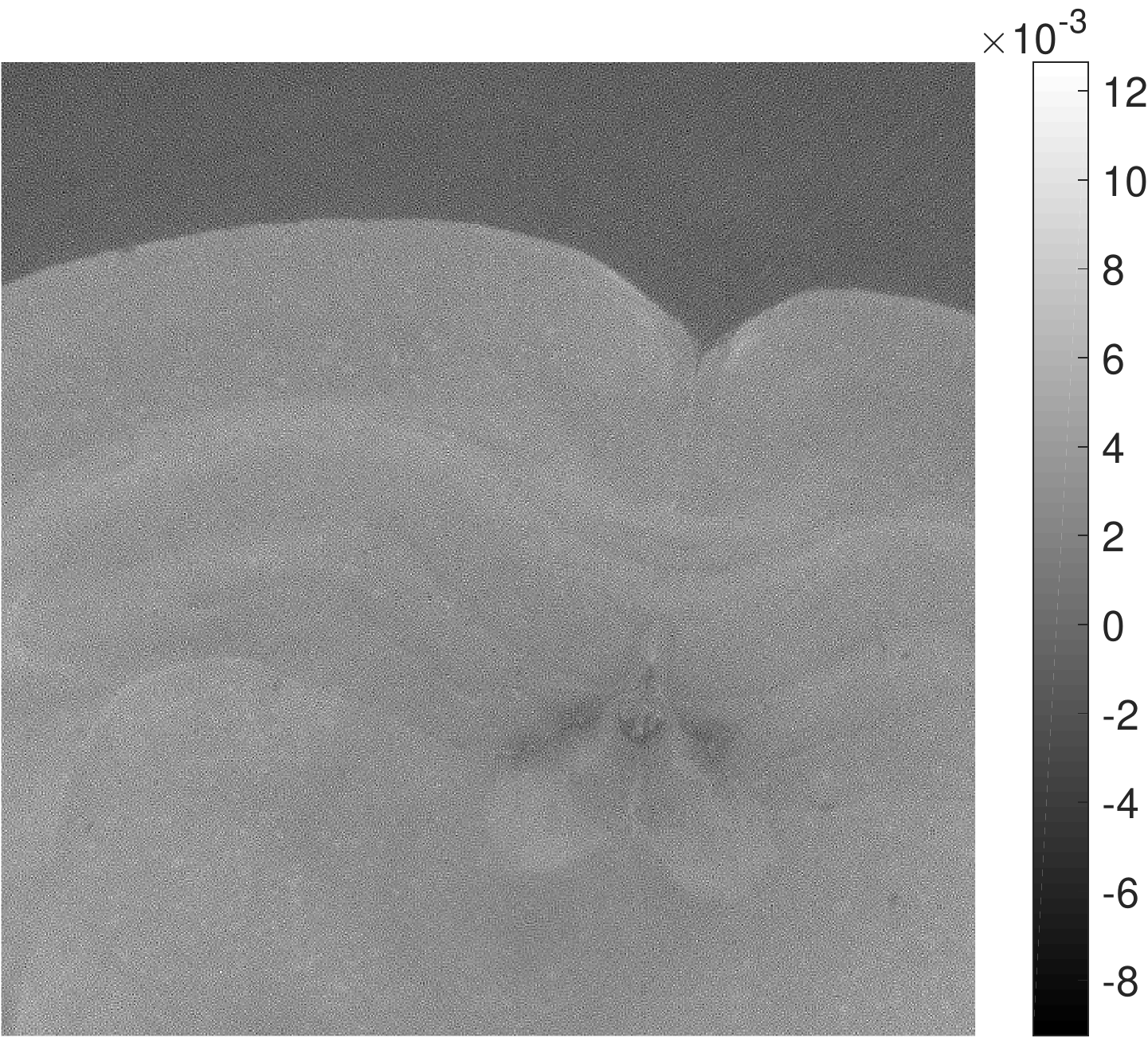}\label{fig:zoom:c}}\qquad\qquad
\subfloat[GBiT, central model]{\includegraphics[width=0.4\linewidth]{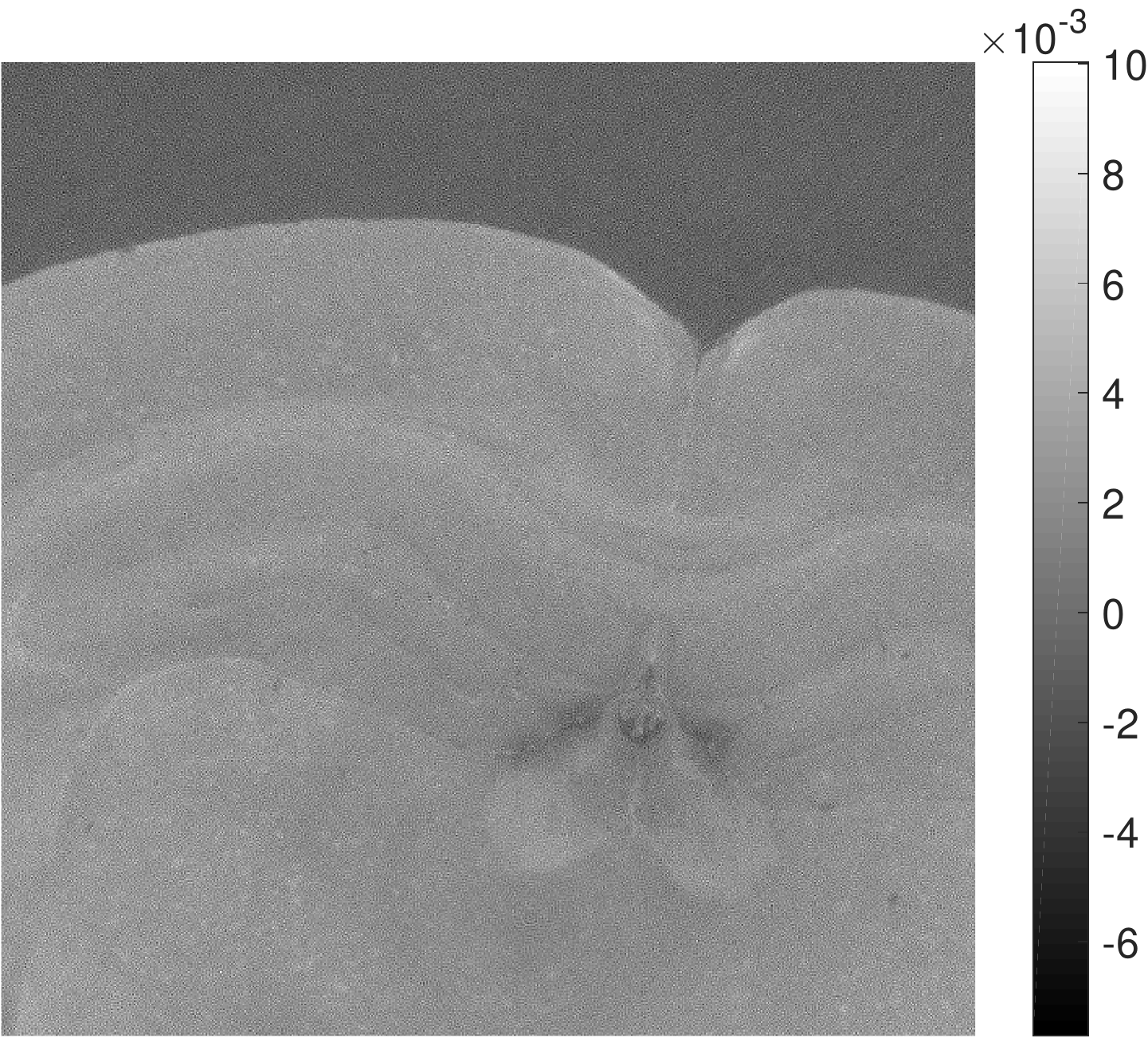}\label{fig:zoom:d}}
\caption{Details from the various reconstructions: \protect\subref{fig:zoom:a} detail from \ref{fig:real:b}, \protect\subref{fig:zoom:b} detail from \ref{fig:real:e}, \protect\subref{fig:zoom:c} detail from \ref{fig:real:c}, \protect\subref{fig:zoom:d} detail from \ref{fig:real:f}.}
\label{fig:zoom}
\end{figure}

%%%%%%%%%%%%%%%%%%%%%%%% SECTION: Discussion and Conclusions. %%%%%%%%%%%%%%%%%%%%%%%%
\section{Conclusions and Remarks}
In summary, we adapted the generalized Arnoldi-Tikhonov method so that it can be applied to the general Tikhonov problem for non-square matrices. This was done by using the bidiagonal decomposition in order to rewrite the corresponding linear system. We tested our method by applying it to the PC-CT problem, where we demonstrated that the forward difference approximation resulted in much better reconstructions than the central difference approximation. Combining our results we were able to find high-contrast reconstructions coming from experimental data using the ASTRA toolbox. Although our CPU implementation is significantly slower that the FBP algorithm, we can improve the reconstruction time significantly using a GPU implementation.

%%%%%%%%%%%%%%%%%%%%%%%% SECTION: Acknowledgements %%%%%%%%%%%%%%%%%%%%%%%%
\section*{Acknowledgements}
We acknowledge Franz Pfeiffer, TUM, for providing the experimental data. NS also thanks the Department of Mathematics and Computer Science, U. Antwerpen for financial support.  We also acknowledge fruitful discussions with Siegfried Cools.

%%%%%%%%%%%%%%%%%%%%%%%% BIBLIOGRAPHY %%%%%%%%%%%%%%%%%%%%%%%%
\bibliographystyle{abbrv}
\bibliography{References}

\end{document}